\documentclass[leqno,10pt]{amsart}
\usepackage{amsmath,amscd,amsthm,amsxtra}
\usepackage{epsfig,graphics,color,colortbl}
\usepackage{amssymb,latexsym}
\usepackage{mathabx}
\usepackage{mathrsfs,eucal,upgreek}
\usepackage[poly,all]{xy}
\usepackage{hyperref}
\usepackage{tikz-cd}
\usepackage{arydshln}
\usepackage{ytableau}
\usepackage[draft]{todonotes}

\newtheorem{thm}{\bf Theorem}[section]
\newtheorem{df}[thm]{\bf Definition}
\newtheorem{prop}[thm]{\bf Proposition}
\newtheorem{cor}[thm]{\bf Corollary}
\newtheorem{lem}[thm]{\bf Lemma}
\newtheorem{rem}[thm]{\bf Remark}
\newtheorem{ex}[thm]{\bf Example}

\numberwithin{equation}{section}

\newcommand{\A}{\mathcal{A}}

\newcommand{\W}{\mathcal{W}}
\newcommand{\cP}{\mathscr{P}}

\newcommand{\pf}{\noindent{\bfseries Proof. }}
\newcommand{\ov}{\overline}

\newcommand{\F}{\mathscr{F}}

\newcommand{\gl}{\mathfrak{gl}}
\newcommand{\Z}{\mathbb{Z}}
\newcommand{\C}{\mathbb{C}}

\newcommand{\te}{\widetilde{e}}
\newcommand{\tf}{\widetilde{f}}
\newcommand{\g}{\mathfrak{g}}

\newcommand{\mc}{\mathcal}
\newcommand{\mf}{\mathfrak}

\newcommand{\la}{\lambda}

\newcommand{\tl}[1]{\substack{\scalebox{0.85}{#1}}}

\newcommand{\llfloor}{\big\lfloor}
\newcommand{\rrfloor}{\big\rfloor}

\newcommand{\llceil}{\big\lceil}
\newcommand{\rrceil}{\big\rceil}

\newcommand{\Cequiv}{\stackrel{C}{\equiv}}

\begin{document}
\title[Combinatorial Howe duality of symplectic type]
{Combinatorial Howe duality of symplectic type}
\author{TAEHYEOK HEO}

\address{Department of Mathematical Sciences, Seoul National University, Seoul 08826, Korea}
\email{gjxogur123@snu.ac.kr}

\author{JAE-HOON KWON}

\address{Department of Mathematical Sciences and RIM, Seoul National University, Seoul 08826, Korea}
\email{jaehoonkw@snu.ac.kr}


\thanks{This work is supported by the National Research Foundation of Korea(NRF) grant funded by the Korea government(MSIT) (No.\,2019R1A2C1084833 and 2020R1A5A1016126).}

\begin{abstract}
We give a new combinatorial interpretation of Howe dual pairs of the form $(\g,{\rm Sp}_{2\ell})$, 
where $\g$ is a Lie (super)algebra of classical type. 
This is done by establishing a symplectic analogue of the RSK algorithm associated to this pair, 
in a uniform way which does not depend on $\g$. 
We introduce an analogue of jeu de taquin sliding for spinor model of irreducible characters of a Lie superalgebra $\g$ to define a $P$-tableau 
and show that the associated $Q$-tableau is given by a symplectic tableau due to King. 
\end{abstract}

\maketitle
\setcounter{tocdepth}{1}

\section{Introduction} \label{sec:intro}
The Robinson-Schensted-Knuth (simply RSK) algorithm or correspondence is one of the most fundamental results 
in the theory of symmetric functions with various applications and generalizations.
From a viewpoint of representation theory, it can be regarded as a combinatorial aspect of Howe duality \cite{H89, H95}, 
which is a more general principle with important applications in many areas of mathematics.  
Indeed, the dual pair $({\rm GL}_m,{\rm GL}_n)$, which acts on the symmetric or exterior algebra 
generated by $\C^m\otimes \C^n$ as mutual centralizers, yields the Cauchy identity or its dual as its associated character. 
It has been generalized to a pair $(\g,{\rm GL}_{n})$, where $\g$ is a Lie (super)algebra of type $A$ 
and an irreducible unitarizable highest weight $\g$-module (not necessarily integrable) appears in the duality (see \cite{CW03,Fr,H89,KacR2}). 
A uniform description of the RSK correspondence for this dual pair $(\g, {\rm GL}_n)$ 
together with a combinatorial model of the associated irreducible $\g$-modules is studied in \cite{K08}.

The purpose of this paper is to establish a symplectic analogue of the RSK algorithm associated to various Howe dual pairs, which include the symplectic group ${\rm Sp}_{2\ell}$.
The duality associated to a pair $(\g,{\rm Sp}_{2\ell})$, which we are interested in this paper, 
can be described as follows (see \cite{H89,H95,HTW,LZ,W99} and references therein). 
Let $\A$ be a $\Z_2$-graded linearly ordered set and 
let $\mathscr{E}_\A$ be the super exterior algebra generated by the superspace with a linear basis indexed by $\A$. 
Then $\F_\A=\mathscr{E}_\A^\ast \otimes \mathscr{E}_\A$ is a semisimple module over a classical Lie (super)algebra $\g_\A$, the type of which depends on $\A$, 
and the $\ell$-fold tensor power $\F^{\otimes \ell}_\A$ ($\ell\geq 1$) is a $(\g_\A,{\rm Sp}_{2\ell})$-module with the following multiplicity-free decomposition:
\begin{equation}\label{eq:Howe duality type C}
\F^{\otimes \ell}_\A \cong \bigoplus_{(\la, \ell) \in \mc{P}({\rm Sp})_\A}V_{\g_\A}(\la,\ell)\otimes V_{{\rm Sp}_{2\ell}}(\la),
\end{equation}
where the direct sum is over a set $\mc{P}({\rm Sp})_\A$ of pairs $(\la, \ell) \in \mathscr{P} \times \mathbb{Z}_+$ with $\ell(\la)\leq \ell$.
Here $V_{{\rm Sp}_{2\ell}}(\la)$ is the irreducible ${\rm Sp}_{2\ell}$-module corresponding to $\la$, 
and $V_{\g_\A}(\la,\ell)$ is the irreducible highest weight $\g_\mathcal{A}$-module 
corresponding to $V_{{\rm Sp}_{2\ell}}(\la)$ appearing in $\F^{\otimes \ell}_\A$ (see Remark \ref{rem:type Howe duality}).

In \cite{K15,K16}, the second author introduced a combinatorial object called a spinor model of type $BCD$, 
which gives the character of $V_{\g_\A}(\la,\ell)$ in \eqref{eq:Howe duality type C} in case of type $C$. 
As a set, the spinor model ${\bf T}_\A(\la,\ell)$ consists of sequences of usual semistandard tableaux of two-columned shapes 
with letters in $\A$, where two adjacent tableaux satisfy certain configuration. 
It can be viewed as a super analogue of Kashiwara-Nakashima (simply KN) tableaux of type $BCD$ \cite{KashNaka}, 
and has interesting applications including branching multiplicities for classical groups \cite{JK,K18-3}, 
crystal bases of quantum superalgebras of orthosymplectic type \cite{K15,K16}, and generalized exponents \cite{JK,LL}.

It is natural to ask whether we have an analogue of RSK algorithm for \eqref{eq:Howe duality type C} in terms of spinor model, and the main result in this paper is to construct an explicit bijection
\begin{equation}\label{eq:main result}
\xymatrixcolsep{4pc}\xymatrixrowsep{0pc}\xymatrix{
{\bf F}_\A^{\ell} 
\ar@{->}[r] & \ 
\displaystyle{\bigsqcup_{\la\in \mc{P}({\rm Sp})_\A}{\bf T}_\A(\la,\ell)\times {\bf K}(\la,\ell)}
}.
\end{equation}
Here ${\bf F}_\A^{\ell}$ is the set of $2\ell$-tuple of $\A$-semistandard tableaux of single-columned shapes with letters in $\A$, and ${\bf K}(\la,\ell)$ is the set of symplectic tableaux of shape $\la$ due to King \cite{King76}, which gives the character of $V_{{\rm Sp}_{2\ell}}(\la)$. Hence the bijection \eqref{eq:main result} yields the Cauchy type identity which follows from the decomposition \eqref{eq:Howe duality type C} for arbitrary $\A$. 

As a special case of our bijection \eqref{eq:main result} by putting $\A$ to be a finite set of $n$ elements with degree $0$, we recover the following well-known identity \cite{King75}
\begin{equation}\label{eq:dual spinors}
\prod_{i=1}^n \prod_{j=1}^\ell (x_i+x_i^{-1}+z_j+z_j^{-1}) = 
\sum_{\lambda \subseteq (n^\ell)}{sp_{\rho_n(\la, \ell)}(\mathbf{x})} sp_\lambda(\mathbf{z}),
\end{equation}
where $\rho_n(\la, \ell)$ is the transpose of the rectangular complement of $\la$ in $(n^\ell)$, 
and $sp_\la({\bf z})$ is the character of $V_{{\rm Sp}_{2\ell}}(\la)$ in $z_1^{\pm},\dots,z_\ell^{\pm}$ 
corresponding to $\la$.

On the other hand, by putting $\A$ to be a finite set of $n$ elements with degree 1 
and using the stability of ${\bf T}_\A(\la,\ell)$ for $\ell\geq n$ \cite{K18-3}, 
we also recover another well-known classical identity due to Littlewood \cite{Li} and Weyl \cite{Weyl},
\begin{equation}\label{eq:Littlewood}
\frac{1}{\prod_{i=1}^n \prod_{j=1}^\ell (1-x_iz_j)(1-x_iz_j^{-1})} = \sum_{\ell(\la)\leq n} sp_\lambda(\mathbf{z}) s_\lambda(\mathbf{x}) \prod_{1 \le i < j \le n}{(1-x_ix_j)^{-1}}, \\
\end{equation}
where $s_\la({\bf x})$ is the Schur polynomial in $x_1,\dots,x_n$, and $\ell(\la)$ is the length of $\la$. 

We remark that our bijection reduced to these cases is completely different from the ones in \cite{Tera} and \cite{Sun} for \eqref{eq:dual spinors} and \eqref{eq:Littlewood}, respectively, where the insertion algorithm in terms of the King tableaux is used.

To construct a bijection \eqref{eq:main result}, we first introduce an analogue of jeu de taquin sliding for spinor model of a skew shape, which plays a crucial role in this paper. 
If $\A$ is a finite set with degree $0$, then a spinor model of a skew shape is naturally in one-to-one correspondence 
with a set of symplectic tableaux of a skew shape in the sense of \cite{KashNaka}, 
and our sliding coincides with the symplectic sliding due to Sheats \cite{Sh} (see also \cite{Le02}).
A key observation to generalize the symplectic sliding algorithm to a spinor model with respect to arbitrary $\A$ 
is that in terms of spinor model, the Sheats' algorithm can be described as a sequence of Kashiwara operators 
with respect to an $\mf{sl}_{2\ell}$-crystal structure on ${\bf F}_\A^{\ell}$ defined by the jeu de taquin sliding of type $A$ \cite{La}. 
This implies that the sliding does {\em not} depend on the choice of $\A$, 
and enables us to define a symplectic analogue of RSK algorithm in a uniform way.
In this sense, the jeu de taquin sliding for spinor model introduced here is similar to the case of type $A$, 
where the related combinatorics depends not essentially on the choice of the set of letters $\A$ (cf.~\cite{BR,Remmel}).

Then we show that for ${\bf T}\in {\bf F}_\A^{\ell}$, 
there exists a unique tableau ${\tt P}({\bf T})\in {\bf T}_\A(\la,\ell)$ 
which can be obtained by applying our jeu de taquin sliding. 
We also define the recording tableau corresponding to ${\tt P}({\bf T})$, which is a sequence of oscillating tableaux, 
and then show that it corresponds bijectively to a tableau ${\tt Q}({\bf T})$ in ${\bf K}(\la,\ell)$ thanks to a recent result by Lee \cite{Lee}.

It would be interesting to find other applications of our RSK algorithm and new interpretation of symplectic sliding in terms of type $A$ crystals. We first expect an orthogonal analogue of RSK correspondence associated to Howe dual pairs $(\g,{\rm O}_n)$ by using the spinor model of type $B$ and $D$ \cite{K15,K16}. The case of type $D$ is more interesting, where a jeu de taquin sliding for KN tableaux in this case is not known yet \cite{Le03}. Also we may adopt the spinor model to realize the crystals of Kirillov-Reshetikhin crystals of classical affine type corresponding to fundamental weights to describe a combinatorial $R$ matrix and energy function on their tensor product, which is also closely related with Kazhdan-Lusztig polynomials (cf.~\cite{LOS12}).

The paper is organized as follows: after a brief review on necessary materials in Section \ref{sec:prel}, 
we recall the definition of spinor model of type $C$ in Section \ref{sec:spinor model}.  
In Section \ref{sec:KN tableaux}, we review the symplectic sliding on KN tableaux of a skew shape introduced in \cite{Sh} and further developed in \cite{Le02} by using crystals.
In Section \ref{sec:sliding for spinor}, we define an analogue of the jeu de taquin sliding for spinor model and then the $P$-tableau ${\tt P}({\bf T})$ for ${\bf T}\in {\bf F}^{\ell}_\A$.
In Section \ref{sec:recording for spinor}, we introduce the recording tableau or $Q$-tableaux for ${\tt P}({\bf T})$, say $Q({\bf T})$, which is a sequence of oscillating tableau, and show that it naturally corresponds to a King tableau $\texttt{Q}({\bf T})$.
In Section \ref{sec:RSK for spinor}, we show that the map ${\bf T}\mapsto ({\tt P}({\bf T}),{\tt Q}({\bf T}))$ gives a bijection in \eqref{eq:main result}.\newline

{\bf Acknowledgement}
The authors would like to thank S.-J. Lee for his kind explanation of his work \cite{Lee}.

\section{Preliminaries} \label{sec:prel}

\subsection{Notations} 
Let $\Z_+$ denote the set of non-negative integers.
Let $\cP$ be the set of partitions or Young diagrams $\la=(\la_1,\la_2,\dots)$. 
We denote by $\la'=(\la'_1,\la'_2,\dots)$ the conjugate of $\la$, 
and by $\lambda^\pi$ the skew Young diagram obtained by $180^{\circ}$-rotation of $\lambda$. 
For $n\geq 1$, let $\cP_{n}=\{\,\la\in\cP\,|\,\ell(\la)\leq n\,\}$, 
where $\ell(\la)$ is the length of $\la$. 

Let $\mc{A}$ be a linearly ordered set with a $\mathbb{Z}_2$-grading $\mc{A}=\mc{A}_0\sqcup\mc{A}_1$.
For $n\geq 1$, we let 
\begin{equation*}
\begin{split}
& [n] = \{ 1 < 2 < \cdots < n \}, \quad [\bar{n}] = \{ \bar{n} < \overline{n-1} < \cdots < \bar{1} \},\\
& \mc{I}_n = [n]\cup [\ov{n}] = \{ 1 < \cdots < n < \bar{n} < \cdots < \bar{1} \},
\end{split}
\end{equation*}
where we assume that all the entries of these sets are assumed to be of degree $0$, and
\[ [n]' = \{ 1' < 2' < \cdots < n' \}, \]
where we assume that all the entries are assumed to be of degree $1$.
For positive integers $m$ and $n$, let
\[ \mathbb{I}_{m|n} = \{ 1 < 2 < \cdots < m < 1' < 2' < \cdots < n' \} \]
with $(\mathbb{I}_{m|n})_0 = [m]$ and $(\mathbb{I}_{m|n})_1 = [n]'$.

For a skew Young diagram $\lambda/\mu$, 
let $SST_\A(\lambda/\mu)$ be the set of {\it $\mathcal{A}$-semistandard} (or semistandard) tableaux of shape $\lambda/\mu$, 
that is, tableaux with entries in $\A$ such that 
(1) the entries in each row (resp. column) are weakly increasing from left to right (resp. from top to bottom), 
(2) the entries in $\mc{A}_0$ (resp. $\mc{A}_1$) are strictly increasing in each column (resp. row).
 
For $T\in SST_\A(\lambda/\mu)$, let $w(T)$ be a word obtained by reading the entries of $T$ column by column 
from right to left, and from top to bottom in each column (cf.~\cite{Ful}).
For $T\in SST_\A(\lambda)$ and $a\in \A$, 
we denote by $a \rightarrow T$ the tableau obtained by the column insertion of $a$ into $T$ (cf. \cite{BR, Ful}). 
For a word $w=w_1\cdots w_r$, we define $(w\rightarrow T)=(w_r\rightarrow (\cdots \rightarrow(w_1\rightarrow T)))$. 
For a semistandard tableau $S$, we define $(S\rightarrow T)=(w(S)\rightarrow T)$.

For $a,b,c\in \Z_+$, let $\lambda(a, b, c)=(2^{b+c}, 1^a)/(1^b)$
be a skew Young diagram with two columns. 
$$
\lambda(2,1,3)\ =\
\resizebox{.055\hsize}{!}{$
{\def\lr#1{\multicolumn{1}{|@{\hspace{.75ex}}c@{\hspace{.75ex}}|}{\raisebox{-.04ex}{$#1$}}}\raisebox{-.6ex}
{$\begin{array}{cccc}
\cline{2-2}
\cdot &\lr{}\\
\cline{1-1}\cline{2-2}
\lr{} &\lr{}\\
\cline{1-1}\cline{2-2}
\lr{} &\lr{\ \ }\\
\cline{1-1}\cline{2-2}
\lr{}&\lr{}\\
\cline{1-1}\cline{2-2}
\lr{\ \ }& \cdot \\
\cline{1-1}
\lr{} & \cdot \\
\cline{1-1}
\end{array}$}}$}$$
For $T\in SST_\A(\la(a, b, c))$, let $T^{\texttt{L}}$ and $T^{\texttt{R}}$ denote the left and right columns of $T$, respectively. 

If necessary, we assume that a tableau is placed on the plane $\mathbb{P}_L$ with a horizontal line $L$. 
Let $U_1,\dots, U_r$ be column tableaux (that is, tableaux of single-columned shapes), 
which are $\A$-semistandard. 
Let $\llfloor U_1,\dots,U_r \rrfloor$ denote the tableau in $\mathbb{P}_L$ (not necessarily of partition shape), 
where the $i$-th column from the left is $U_i$ and its bottom edge lies on $L$. 
Similarly, let $\llceil U_1,\dots,U_r \rrceil$ denote the tableau, 
where the $i$-th column from the left is $U_i$ and its top edge lies on $L$.

For $(u_1,\dots, u_r) \in \Z_+^r$, let 
\begin{equation*}
\llfloor U_1,\dots,U_r \rrfloor_{(u_1,\dots,u_r)},\quad 
\llceil U_1,\dots,U_r \rrceil^{(u_1,\dots,u_r)}
\end{equation*}
be the tableaux obtained from $\llfloor U_1,\dots,U_r \rrfloor$ and $\llceil U_1,\dots,U_r \rrceil$ 
by sliding each $U_i$ by $u_i$ positions up and down, respectively.

\begin{ex}{\rm
\[ \begin{tikzpicture}
\ytableausetup {mathmode, boxsize=1.2em}
\node [above right] at (0, 0) {
\begin{ytableau}
\none & \none & \tl 2 & \tl{$2'$} \\
\tl 1 & \tl 1 & \tl{$1'$} & \tl{$2'$} \\
\tl 2 & \tl 3 & \tl{$1'$} \\ 
\tl{$1'$} & \tl{$2'$} & \tl{$3'$} \\ 
\tl{$3'$}
\end{ytableau} };
\node [above right] at (7, 0) {
\begin{ytableau} 
\none & \none & \tl 3 \\ 
\none & \tl 3 & \bar{\tl 5} \\ 
\tl 3 & \bar{\tl 5} & \bar{\tl 3} \\ 
\tl 5 & \bar{\tl 4} \\
\bar{\tl 3} & \bar{\tl 2} \\
\bar{\tl 1}
\end{ytableau}};

\draw [dotted] (-0.3, 0.362em) -- (3, 0.362em);
\node at (-0.5, 0.4em) {\small $L$};
\draw [dotted] (6.8, 7.83em) -- (9, 7.83em);
\node at (9.3, 7.83em) {\small $L$};
\node [below] at (1.8, -0.2) {$\llfloor U_1, U_2, U_3, U_4 \rrfloor_{(0, 1, 1, 3)}$};
\node [below] at (1.8, -1) {\small $\in SST_{\mathbb{I}_{4|3}}((4, 4, 3, 3, 1)/(2))$};
\node [below] at (8.25, -0.2) {$\llceil U_1, U_2, U_3 \rrceil^{(2, 1, 0)}$};
\node [below] at (8.25, -1) {\small$\in SST_{\mathcal{I}_5}((3, 3, 3, 2, 2, 1)/(2, 1))$};
\end{tikzpicture} \]
} \end{ex}


\subsection{Crystal and Sch\"{u}tzenberger's jeu de taquin} 
Let us briefly recall the notion of crystals (see \cite{HK,Kas95} for more details).
Let $\mathfrak{g}$ be a symmetrizable Kac-Moody algebra associated to a generalized Cartan matrix $A = (a_{ij})_{i, j \in I}$ indexed by $I$.
A {\em $\g$-crystal} is a set $B$ together with the maps ${\rm wt} : B \rightarrow P$, $\varepsilon_i, \varphi_i: B \rightarrow \mathbb{Z}\cup\{-\infty\}$ 
and $\te_i, \tf_i: B \rightarrow B\cup\{{\bf 0}\}$ for $i\in I$ satisfying certain axioms, 
where ${\bf 0}$ is a formal symbol and $P$ is the weight lattice of $\mathfrak{g}$. 
For a dominant integral weight $\Lambda$, we denote by $B(\Lambda)$ the crystal associated to an irreducible highest weight $\g$-module with highest weight $\Lambda$. 
A crystal $B$ is called regular if it is isomorphic to a disjoint union of $B(\Lambda)$'s. In this case, we have 
$\varepsilon_i(b)=\max\{\,k\,|\,\te_i^kb\neq {\bf 0}\,\}$ and 
$\varphi_i(b)=\max\{\,k\,|\,\tf_i^kb\neq {\bf 0}\,\}$ for $b\in B$ and $i\in I$. 
For example, if $\A=[r]$ with $r\geq 2$, then $SST_{[r]}(\la)$ is a connected regular $\mf{gl}_r$-crystal for $\la\in \cP_r$ \cite{KashNaka}.

Let $\mc{A}$ be a $\mathbb{Z}_2$-graded linearly ordered set.
The Sch\"{u}tzenberger's jeu de taquin is also available for $\A$-semistandard tableaux (cf. \cite{BR, Remmel}), 
where we apply sliding process provided that the resulting tableau remains semistandard.
For example, when $\mathcal{A} = \mathbb{I}_{2|2}$, we have    
  \[
    \begin{ytableau} \none & 1 \\ 1 & 2 \end{ytableau}
    \ \longrightarrow \
    \begin{ytableau} 1 & 1 \\ 2 \end{ytableau}\ , \quad 
    \begin{ytableau} \none & \tl{$1'$} \\ \tl{$1'$} & \tl{$2'$} \end{ytableau}
    \ \longrightarrow \
    \begin{ytableau} \tl{$1'$} & \tl{$2'$} \\ \tl{$1'$} \end{ytableau}.
  \]

We use this algorithm in terms of crystal operator $\mc E$ and $\mc F$ for $\mf{sl}_2$ (cf. \cite{La}), 
which plays an important role in this paper.

For $T \in SST_\A(\lambda(a, b, c))$ and $0 \le k \le \min\{ a, b \}$, consider a tableau $T'$ of shape $\lambda(a-k, b-k, c+k)$ by sliding down $T^{\texttt{R}}$ by $k$ positions. Let $\mathfrak{r}_T$ be the maximal $k$ such that $T'$ is semistandard.

For $T\in SST_\A(\la(a,b,c))$ with ${\mf r}_T=0$, we define  
\begin{itemize}
\item[(1)] $\mc E T$ to be the tableau in $SST_\A(\la(a-1,b+1,c))$ obtained from $T$ by applying jeu de taquin sliding to the position below the bottom of $T^{\texttt{R}}$ when $a>0$,

\item[(2)] $\mc F T$ to be the tableau in $SST_\A(\la(a+1,b-1,c))$ obtained from $T$ by applying jeu de taquin sliding to the position above the top of $T^{\texttt{L}}$ when $b>0$.

\end{itemize}
Here we assume that $\mathcal{E}T = \mathbf{0}$ and $\mathcal{F}T = \mathbf{0}$ when $a=0$ and $b=0$, respectively, where ${\bf 0}$ is a formal symbol. 

\begin{rem} \label{rem:crystalstr}{\rm 
(1) When we apply jeu de taquin sliding to a corner of $T$, we have $\mathfrak{r}_T >0$ if and only if a vertical move occurs by \cite[Lemma 6.2]{K15}. Thus, we have ${\mf r}_{\mc E T}=0$ and ${\mf r}_{\mc F T}=0$ whenever ${\mc E}T$ and ${\mc F}T$ are defined.

(2) Let $\varepsilon(T)=\max\{\,k\,|\,\mathcal{E}^kT \ne 0\,\}$, $\varphi(T)=\max\{\,k\,|\,\mathcal{F}^kT \ne 0\,\}$, and let ${\rm wt}(T)=\varphi(T) - \varepsilon(T)$.
Then $B=\{\,{\mc E}^k T\,|\,0\leq k\leq a\,\}\cup \{\,{\mc F}^l T\,|\,0\leq l\leq b\,\}$ 
forms a regular $\mf{sl}_2$-crystal with respect to $\mc E$ and $\mc F$, where we identify the weight lattice of $\mf{sl}_2$ with $\Z$.
Indeed, if we let $P=\Z\epsilon_1\oplus\Z\epsilon_2$ be the weight lattice of $\mf{gl}_2$ with $\epsilon_1-\epsilon_2$ the simple root and define ${\rm wt}(T) = m_1\epsilon_1+m_2\epsilon_2$, 
where $m_1$ (resp. $m_2$) is the number of boxes of $T^{\tt R}$ (resp. $T^{\tt L}$), 
then we may regard $B$ as a $\mathfrak{gl}_2$-crystal.}
\end{rem}

\begin{ex} 
{\rm 

Suppose that $\mathcal{A} = \mathbb{I}_{4|3}$.

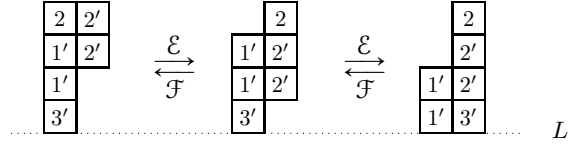
\[ \begin{tikzpicture}
\node [above right] at (0, 0) { 
\begin{ytableau} 
\tl 2 & \tl{$2'$} \\ 
\tl{$1'$} & \tl{$2'$} \\ 
\tl{$1'$} \\ 
\tl{$3'$}
\end{ytableau}};
\node [above right] at (2.5, 0) {
\begin{ytableau} 
\none & \tl 2 \\ 
\tl{$1'$} & \tl{$2'$} \\ 
\tl{$1'$} & \tl{$2'$} \\ 
\tl{$3'$} 
\end{ytableau}};
\node [above right] at (5, 0) {
\begin{ytableau} 
\none & \tl 2 \\ 
\none & \tl{$2'$} \\ 
\tl{$1'$} & \tl{$2'$} \\ 
\tl{$1'$} & \tl{$3'$} 
\end{ytableau}};

\node at (1.87, 3em) {$\longrightarrow$}; 
\node at (4.4, 3em) {$\longrightarrow$};
\node [above] at (1.87, 3em) {$\mc E$}; 
\node [above] at (4.4, 3em) {$\mc E$};
\node at (1.87, 2.6em) {$\longleftarrow$}; 
\node at (4.4, 2.6em) {$\longleftarrow$};
\node [above] at (1.87, 1.35em) {$\mc F$}; 
\node [above] at (4.4, 1.35em) {$\mc F$};

\draw [dotted] (-0.3, 0.362em) -- (6.5, 0.362em);
\node [below] at (7, 0.4) {\small $L$};
\end{tikzpicture} \]
}
\end{ex}

Now, for $(U,V)\in SST_{\A}((1^u))\times SST_{\A}((1^v))$ ($u,v\in\Z_+$), we define 
\begin{equation}\label{eq:E and F via jdt}
\begin{split}
\mc{X}(U,V)&=
\begin{cases}
((\mc{X}T)^{\tt L},(\mc{X}T)^{\tt R}) &\text{if $\mc{X}T\neq {\bf 0}$},\\
{\bf 0} &\text{if $\mc{X}T={\bf 0}$},
\end{cases} \quad\quad (\mc{X}=\mc{E}, \mc{F}),
\end{split}
\end{equation}
where $T$ is the unique tableau in $SST_{\A}(\la(u-k,v-k,k))$ for some $0\leq k\leq \min\{u,v\}$ such that $(T^{\tt L},T^{\tt R})=(U,V)$ and ${\mf r}_T=0$. 


\subsection{Crystal and RSK correspondence} 
Let $\mc{A}$ be a $\mathbb{Z}_2$-graded linearly ordered set.
For $r\geq 2$, let 
\begin{equation}\label{eq:Fock space}
{\bf E}_\A^r=\bigsqcup_{(u_r,\dots,u_1)\in \Z_+^r} 
SST_{\A}((1^{u_r}))\times \dots\times SST_{\A}((1^{u_1})).
\end{equation}
For $(U_r,\dots,U_1)\in {\bf E}^r_\A$ and $1\leq i\leq r-1$, we define
\begin{equation}
\begin{split}
\mc{X}_i(U_r,\dots,U_1) &= 
\begin{cases}
(U_r,\dots, \mc{X}(U_{i+1},U_i),\dots,U_1) & \text{if $\mc{X}(U_{i+1},U_i)\neq {\bf 0}$},\\
{\bf 0} &\text{if $\mc{X}(U_{i+1},U_i)={\bf 0}$},
\end{cases}
\end{split}
\end{equation} 
where $\mc{X}(\, \cdot\, ,\, \cdot \, )$ for $\mc{X}=\mc{E},\mc{F}$ is defined in \eqref{eq:E and F via jdt}.

Let $P=\bigoplus_{i=1}^r\Z\epsilon_i$ be the weight lattice of $\mf{gl}_r$, where $\{\,\alpha_i=\epsilon_i-\epsilon_{i+1}\,|\,i=1,\dots,r-1\,\}$ is the set of simple roots.
Given $T = (U_r, \dots, U_1) \in \mathbf{E}_\mathcal{A}^r$, let ${\rm wt}(T) = \sum_{i=1}^r m_i\epsilon_i$, where $m_i$ is the number of boxes of $U_i$.

\begin{lem}\label{lem:regularity of E}
${\bf E}^r_\A$ is a regular $\mf{gl}_r$-crystal with respect to ${\rm wt}$, $\mc{E}_i$ and $\mc{F}_i$ for $1\leq i\leq r-1$.
\end{lem}
\pf It can be proved by similar arguments as in \cite{K09,La}.
Let ${\bf M}_{\A\times [r]}$ be the set of matrices $\mathbf{m}= (m_{ab})$ with non-negative integral entries $(a\in \A,\ b \in [r])$ satisfying
(1) $m_{a,b} \in \{0,1\}$ if $a\in \A_0$, (2) $\sum_{ab} m_{ab} < \infty$.
There is a natural bijection from ${\bf E}^r_\A$ to ${\bf M}_{\A\times [r]}$, 
where $(U_r,\dots,U_1) \in \mathbf{E}_\mathcal{A}^r$ is sent to $\mathbf{m}= (m_{ab})$ 
such that $m_{ab}$ is the number of occurrences of $a$ in $U_b$.

Let ${\bf m}=(m_{ab}) \in {\bf M}_{\A\times [r]}$ be given.
For $a\in \A$, we may identify the $a$-th row of ${\bf m}$ with a unique tableau $T^{(a)}$ 
in $SST_{[r]}((u))$ (resp. $SST_{[r]}((1^u))$) if $a \in \mathcal{A}_0$ (resp. $a \in \mathcal{A}_1$), where $u=\sum_{b}m_{ab}$ and $m_{ab}$ is the number of occurrences of $b$ in $T^{(a)}$.
We may define a regular $\mf{gl}_r$-crystal structure on ${\bf M}_{\A\times [r]}$ by regarding ${\bf m}$ as $\overset{\longrightarrow}{\bigotimes}_{a\in \A}T^{(a)}$.
Then we can check that the associated operators $\te_i$ and $\tf_i$ for $1\leq i\leq r-1$ coincide with $\mc{E}_i$ and $\mc{F}_i$, and the $\mf{gl}_r$-weight is equal to ${\rm wt}$.
Hence $\mathbf{M}_{\mathcal{A} \times [r]}$ is a regular $\mf{gl}_r$-crystal since it is a disjoint union of tensor products of regular $\mf{gl}_r$-crystals $SST_{[r]}((u))$ and $SST_{[r]}((1^u))$. \qed
\vskip 2mm

Let us recall the RSK correspondence, which explains the decomposition of the $\mf{gl}_r$-crystal ${\bf E}^r_\A$ into its connected components. 
Let ${\bf U}=(U_r,\dots,U_1)\in {\bf E}_\A^r$ given. 
Let $P({\bf U})=(U_r\rightarrow (\cdots\rightarrow(U_2\rightarrow U_1) \cdots))$ 
and $Q({\bf U})$ be the corresponding recording tableau, that is, 
if ${\rm sh}(P({\bf U}))=\la$ (the shape of $P({\bf U})$) and 
$P_i=(U_i\rightarrow (\cdots(U_2\rightarrow U_1)\cdots))$ for $1\leq i\leq r$, then $Q({\bf U})$ is the unique tableau in $SST_{[r]}(\la')$ such that its subtableau ${\rm sh}(P_{i})'/{\rm sh}(P_{i-1})'$ is a horizontal strip filled with $i$, where we assume that $P_0$ is the empty tableau.
Then we have a bijection
\begin{equation}\label{eq:RSK}
\xymatrixcolsep{4pc}\xymatrixrowsep{0.5pc}\xymatrix{
{\bf E}^r_\A  \ar@{->}^{\hskip -20mm \kappa_\A}[r] & \  \displaystyle{\bigsqcup_{\la'\in \cP_r}} SST_\A(\la)\times SST_{[r]}(\la') \\
{\bf U}  \ar@{|->}[r] & (P({\bf U}), Q({\bf U}))}.
\end{equation}

\begin{lem} \label{lem:sl crystal}
The bijection $\kappa_\A$ is an isomorphism of $\mf{gl}_r$-crystals, where the right-hand side is a regular $\mf{gl}_r$-crystal with respect to the second component.
\end{lem}
\pf It follows from the argument in the proof of Lemma \ref{lem:regularity of E} and the symmetry of RSK correspondence.
\qed

\begin{rem}{\rm
When $\mathcal{A} = \mathbb{I}_{m|n}$, it is shown in \cite{K07} that the bijection \eqref{eq:RSK} is an isomorphism of $(\mathfrak{gl}_{m|n}, \mathfrak{gl}_r)$-bicrystals, where $\gl_{m|n}$ is a general linear Lie superalgebra associated to a superspace $\C^{m|n}$. 
This explains the $(\mathfrak{gl}_{m|n}, \mathfrak{gl}_r)$-duality \cite{CW01, H89} in terms of crystals (see \cite{K09} for the case of arbitrary $\mathcal{A}$).
}
\end{rem}


\section{Spinor model of symplectic type}\label{sec:spinor model}

\subsection{Spinor model of type $C$}\label{subsec:spinor model} 
Let us recall our main combinatorial object, which is introduced in \cite{K15,K16}.
Let 
\begin{equation*}
\mc{P}({\rm Sp})=\{\,(\la,\ell)\,|\,\ell\geq 1,\ \la\in \cP_\ell\,\}.
\end{equation*}
Hereafter $\mc{A}$ denotes a $\mathbb{Z}_2$-graded linearly ordered set unless otherwise specified.
For $a\in \Z_+$, let
\begin{equation*}
{\bf T}_\A(a)=\bigsqcup_{c\in\Z^+}SST_{\A}(\la(a,0,c)).
\end{equation*}
Note that, for $a \ne 0$, ${\bf T}_\A(a)$ is non-empty if and only if $SST_\A((1^a))\neq \emptyset$.
For $T\in {\bf T}_\A(a)$, we have ${\mc E}^a T\in SST_\A(\la(0,a,c))$ and  define 
\begin{equation*}\label{eq:LT RT}
{}^{\tt L}T=({\mc E}^a T)^{\tt L}, \quad   {}^{\tt R}T= ({\mc E}^a T)^{\tt R}.
\end{equation*}

\begin{ex} 
{\rm 

Suppose $\mathcal{A} = \mathbb{I}_{4|3}$ and $T\in SST_{\A}(\la(2,0,2))$ as below.

\[ \begin{tikzpicture}
\node [above right] at (0, 0) { 
\begin{ytableau} 
\tl 2 & \tl{$2'$} \\ 
\tl{$1'$} & \tl{$2'$} \\ 
\tl{$1'$} \\ 
\tl{$3'$}
\end{ytableau}};
\node [above right] at (2.5, 0) {
\begin{ytableau} 
\none & \tl 2 \\ 
\tl{$1'$} & \tl{$2'$} \\ 
\tl{$1'$} & \tl{$2'$} \\ 
\tl{$3'$}
\end{ytableau}};
\node [above right] at (5, 0) {
\begin{ytableau} 
\none & \tl 2 \\ 
\none & \tl{$2'$} \\ 
\tl{$1'$} & \tl{$2'$} \\ 
\tl{$1'$} & \tl{$3'$}  
\end{ytableau}};
\node [above right] at (7.5, 0) {
\begin{ytableau} 
\tl{$1'$} & \tl 2 \\ 
\tl{$1'$} & \tl{$2'$} \\ 
\none & \tl{$2'$} \\ 
\none & \tl{$3'$} 
\end{ytableau}};

\node at (-0.5, 3em) {$T=$}; 
\node at (1.87, 3em) {$\longrightarrow$}; 
\node at (4.4, 3em) {$\longrightarrow$};
\node [above] at (1.87, 3em) {$\mc E$}; 
\node [above] at (4.4, 3em) {$\mc E$};
\draw [dotted] (-0.3, 0.362em) -- (10, 0.362em);
\node [below right] at (-0.2, -0.1) 
{\small $\llfloor T^{\tt L}, T^{\tt R}\rrfloor_{(0,2)}$};
\node [below right] at (4.9, -0.1) 
{\small $\llfloor{}^{\tt L}T, {}^{\tt R}T\rrfloor$};
\node [below] at (8.5, -0.1) 
{\small $\llfloor {{}^{\tt L}T}, {{}^{\tt R}T} \rrfloor_{(2,0)}$};
\node [below] at (10.5, 0.4) {\small $L$};
\end{tikzpicture}
\] } \end{ex}

\begin{df}\label{def:spinor tableaux}{\rm \mbox{}
\begin{itemize}
\item[(1)]
For $a_1, a_2\in \Z_+$ with $a_2\leq a_1$ and $(T_2,T_1)\in {\bf T}_\A(a_2)\times {\bf T}_\A(a_1)$, we define  
\begin{equation*}
\text{$T_2\prec T_1$\quad if  $\llfloor {}^{\tt R}T_2, {T}^{\tt L}_1 \rrfloor$ and $\llfloor T_2^{\tt R},{}^{\tt L}T_1 \rrfloor_{(a_2,a_1)}$ are $\A$-semistandard.}
\end{equation*}

\item[(2)] For $(\la,\ell)\in \mc{P}({\rm Sp})$,
we define
\begin{equation*}
{\bf T}_{\A}(\la,\ell)=
\{\,{\bf T}=(T_\ell,\dots,T_1)\ |\ \text{$T_{\ell}\prec \cdots \prec T_{1}$}\,\}\subset {\bf T}_\A(\la_\ell)\times \cdots \times {\bf T}_\A(\la_1). 
\end{equation*}
\end{itemize}
}
\end{df}

\begin{ex} 
{\rm
Let $\mathcal{A} = \mathbb{I}_{4|3}$, and take $S \in \mathbf{T}_\mathcal{A}(1)$ and $T \in \mathbf{T}_\mathcal{A}(2)$ as follows.
\[ \begin{tikzpicture}
\node [above right] at (0, 0) 
{$\begin{ytableau}
\none \\
\tl 2 & \tl 2 \\ 
\tl 3 & \tl 4 \\
\tl{$1'$} 
\end{ytableau}$};
\node [above right] at (3, 0) 
{$\begin{ytableau} 
\none & \none \\ 
\tl 2 & \tl 2 \\ 
\tl 3 & \tl 4 \\ 
\none & \tl{$1'$}
\end{ytableau}$};
\node [above right] at (6.5, 0) 
{$\begin{ytableau} 
\tl 1 & \tl 1 \\ 
\tl 3 & \tl{$2'$} \\ 
\tl{$1'$} \\ 
\tl{$2'$} 
\end{ytableau}$};
\node [above right] at (9.5, 0) 
{$\begin{ytableau} 
\tl 1 & \tl 1 \\ 
\tl{$1'$} & \tl 3 \\ 
\none & \tl{$2'$} \\ 
\none & \tl{$2'$} 
\end{ytableau}$};
\draw [dotted] (-0.3, 0.362em) -- (11.2, 0.362em);
\node [below] at (0.9, -0.2) {$\llfloor S^\texttt{L}, S^\texttt{R} \rrfloor_{(0, 1)}$}; 
\node [below right] at (2.85, -0.2) {$\llfloor {^\texttt{L}S}, {^\texttt{R}S} \rrfloor_{(1, 0)}$}; 
\node [below] at (7.4, -0.2) {$\llfloor T^\texttt{L}, T^\texttt{R} \rrfloor_{(0, 2)}$}; 
\node [below right] at (9.35, -0.2) {$\llfloor {^\texttt{L}T}, {^\texttt{R}T} \rrfloor_{(2, 0)}$}; 
\end{tikzpicture} \]
Then $S \prec T$ since $\llfloor {^\texttt{R}S}, {T^\texttt{L}} \rrfloor$ and $\llfloor S^\texttt{R}, {^\texttt{L}T} \rrfloor_{(1, 2)}$ form semistandard tableaux.

\[ \begin{tikzpicture}
\node [above right] at (0, 0) {$\begin{ytableau} \none & \tl 1 \\ \tl 2 & \tl 3 \\ \tl 4 & \tl{$1'$} \\ \tl{$1'$} & \tl{$2'$} \end{ytableau}$};
\node [above right] at (3, 0) {$\begin{ytableau} \none & \tl 1 \\ \tl 2 & \tl{$1'$} \\ \tl 4 \\ \none \end{ytableau}$};

\draw [dotted] (-0.3, 0.362em) -- (4.6, 0.362em);
\node [below right] at (-0.1, -0.2) {$\llfloor {^\texttt{R}S}, {T^\texttt{L}} \rrfloor$};
\node [below right] at (2.9, -0.2) {$\llfloor {S^\texttt{R}}, {^\texttt{L}T} \rrfloor_{(1, 2)}$};
\end{tikzpicture} \]
}
\end{ex}

Put $$\mc{P}({\rm Sp})_\A=\{\,(\la,\ell)\in \mc{P}({\rm Sp})\,|\,{\bf T}_{\A}(\la,\ell)\neq \emptyset\,\}.$$ 
Then, for any partition $\lambda \ne \emptyset$, we have $SST_{\A}(\la')\neq \emptyset$ if and only if $(\la,\ell)\in \mc{P}({\rm Sp})_\A$.

Let ${\bf x}_\A=\{\,x_a\,|\,a\in \A\,\}$ be the set of formal commuting variables indexed by $\A$.
Let $\mu\in \cP$ be given such that $SST_\A(\mu)\neq \emptyset$.
For $T\in SST_\A(\mu)$, let
${\bf x}_\A^T =\prod_{a\in\A}x_a^{m_a}$, where $m_a$ is the number of occurrences of $a$ in $T$. Note that $s_\mu({\bf x}_\A)=\sum_{T\in SST_{\A}(\mu)}{\bf x}_\A^T$ is a super-analogue of Schur function corresponding to $\mu$.

Let $t$ be a variable commuting with $x_a$ ($a\in \A$).
For $(\la,\ell)\in \mc{P}({\rm Sp})_\A$, we define the character of ${\bf T}_{\A}(\la,\ell)$ to be 
\begin{equation}\label{eq:character of spinor}
S_{(\la,\ell)}({\bf x}_\A)=
t^\ell \sum_{(T_\ell,\dots,T_1)\in {\bf T}_\A(\la,\ell)}
{\bf x}_\A^{T_\ell}\,\cdots\, {\bf x}_\A^{T_1}.  
\end{equation}

\begin{rem}\label{rem:type Howe duality}
{\rm
 The character of ${\bf T}_\A(\la,\ell)$ has an important application in representation theory. 
Indeed, it is motivated by the $(\mf{g}_\A, {\rm Sp}_{2\ell})$-duality \eqref{eq:Howe duality type C} for a Lie (super)algebra $\g_\A$.

Let us first recall the decomposition \eqref{eq:Howe duality type C} for various choices of $\A$.
If $\A=[\ov{m}]$, then we have $(\mf{sp}_{2m}, {\rm Sp}_{2\ell})$-duality, where $V_{\mf{sp}_{2m}}(\la,\ell)$ is a finite-dimensional irreducible $\mf{sp}_{2m}$-module. If $\A=[n]'$, then we have $(\mf{so}_{2n}, {\rm Sp}_{2\ell})$-duality, where $V_{\mf{so}_{2n}}(\la,\ell)$ is an infinite-dimensional irreducible $\mf{so}_{2n}$-module. See \cite{H89, H95, HTW} for these dualities.
In general, when $\mathcal{A} = \mathbb{I}_{m|n}$, we have $(\mf{spo}_{2m|2n}, {\rm Sp}_{2\ell})$-duality \cite{CZ04}, which includes both of the above cases with $n=0$ and $m=0$, respectively. 
Here $\mf{spo}_{2m|2n}$ is the orthosymplectic Lie superalgebra whose Dynkin diagram is given by

  \[ \resizebox{0.8\textwidth}{!}{
  \begin{tikzpicture}
    \draw (-6.5, 0) circle [radius=0.15];
    \node at (-6.5, -0.5) {\small $0$};
    \draw (-5, 0) circle [radius=0.15];
    \node at (-5, -0.5) {\small$1$};
    \draw (-3.5, 0) circle [radius=0.15];
    \node at (-3.5, -0.5) {\small$2$};
    \node at (-2.5, 0) {\small$\cdots$};
    \draw (-1.5, 0) circle [radius=0.15];
    \node at (-1.5, -0.5) {\small${m-1}$};
    \draw (0, 0) circle [radius=0.15];
    \node at (0, -0.5) {\small$m$};    
    \draw (1.5, 0) circle [radius=0.15];
    \node at (1.5, -0.5) {\small$1'$};
    \node at (2.5, 0) {\small$\cdots$};
    \draw (3.5, 0) circle [radius=0.15];
    \node at (3.5, -0.5) {\small${(n-1)'}$};
    \draw (0.1, 0.1) -- (-0.1, -0.1);
    \draw (0.1, -0.1) -- (-0.1, 0.1);

    \draw (-6.5+0.2, 0.05) -- (-5.0-0.15, 0.05);
    \draw (-6.5+0.2, -0.05) -- (-5.0-0.15, -0.05);
    \draw (-6.5+0.15, 0) -- (-6.5+0.3, 0.15);
    \draw (-6.5+0.15, 0) -- (-6.5+0.3, -0.15);
    \draw (-5.0+0.15, 0) -- (-3.5-0.15, 0);
    \draw (-3.5+0.15, 0) -- (-2.5-0.4, 0);
    \draw (-2.5+0.4, 0) -- (-1.5-0.15, 0);
    \draw (-1.5+0.15, 0) -- (-0.0-0.15, 0);
    \draw (0.0+0.15, 0) -- (1.5-0.15, 0);
    \draw (1.5+0.15, 0) -- (2.5-0.4, 0);
    \draw (2.5+0.4, 0) -- (3.5-0.15, 0);
  \end{tikzpicture}
  } \]
The dualities when $\A$ is an infinite $\Z_2$-graded set can be found in \cite{LZ,W99}.

It is shown in \cite{K15} that the character of ${\bf T}_\A(\la,\ell)$ is equal 
to the character of $V_{\g_\A}(\la,\ell)$ when $\mathcal{A} = \mathbb{I}_{m|n}$. 
Indeed, this will also follow from comparing the character identities of \eqref{eq:Howe duality type C} and \eqref{eq:main result} for any $(\mf{g}_\A, {\rm Sp}_{2\ell})$-duality (see Theorem \ref{thm:Cauchy identity}). 
In this sense, we call ${\bf T}_\A(\la,\ell)$ a {\em spinor model of irreducible symplectic characters.}

 }
\end{rem}

\begin{rem}{\rm
The definition of ${\bf T}_\A(\la,\ell)$ in Definition \ref{def:spinor tableaux} is the same as in \cite{K15}. We remark that there is another definition of ${\bf T}_\A(\la,\ell)$ in \cite{K18-3}, which is slightly different from Definition \ref{def:spinor tableaux}, but which has the same character.
}
\end{rem}


\subsection{Schur expansion} 

Let $(\la,\ell)\in \mc{P}({\rm Sp})_\A$ be given. 
Consider an embedding of sets given by
\begin{equation*}
\xymatrixcolsep{4pc}\xymatrixrowsep{0.5pc}\xymatrix{
{\bf T}_\A(\la,\ell)  \ar@{->}[r]^{\iota}  & \  {\bf E}_\A^{2\ell} \\
{\bf T}=(T_\ell,\dots,T_1)  \ar@{|->}[r] & (T_\ell^{\tt L},T_\ell^{\tt R},\dots,T_1^{\tt L},T_1^{\tt R})}
\end{equation*}
and identify ${\bf T}_\A(\la,\ell)$ with its image under $\iota$.
By composing with $\kappa_\A$ in \eqref{eq:RSK}, we have an embedding
\begin{equation*}
\xymatrixcolsep{4pc}\xymatrixrowsep{0.5pc}\xymatrix{
{\bf T}_\A(\la,\ell) \ar@{->}^{\hskip -15mm \Phi_\A}[r] & \  \displaystyle{\bigsqcup_{\mu}}\, SST_\A(\mu)\times SST_{[2\ell]}(\mu') \\
{\bf T}  \ar@{|->}[r] & (P({\bf T}), Q({\bf T}))},
\end{equation*}
where the union is over $\mu\in \cP$ with $\ell(\mu')\leq 2\ell$.

Let us describe the image of $\Phi_\A$ explicitly. Since $SST_{[2\ell]}(\mu')$ is a regular $\mathfrak{gl}_{2\ell}$-crystal, 
there exists an action of the Weyl group of $\mf{gl}_{2\ell}$ on $SST_{[2\ell]}(\mu')$, 
which is isomorphic to $\mf{S}_{2\ell}$ generated by the simple reflection $r_i$ for $1\leq i\leq 2\ell-1$. 
For $Q\in SST_{[2\ell]}(\mu')$, let us define the weight of $Q$ to be the sequence $(m_1,\dots,m_{2\ell})$, 
where $m_i$ is the number of occurrences of $i$ in $Q$, and the $i$-signature of $Q$ to be the sequence 
$(\varepsilon_i(Q),\varphi_i(Q))$ ($1\leq i\leq 2\ell-1$).

For $\mu\in\mathscr{P}$ with $\ell(\mu')\leq 2\ell$, 
let ${K}_{\mu (\lambda, \ell)}$ be the set of $Q \in SST_{[2\ell]}(\mu')$ 
such that its weight $(m_1, \dots, m_{2\ell})$ satisfies the following conditions:
\begin{enumerate}
\item $m_{2k}-m_{2k-1} = \lambda_k$ for $1 \le k \le \ell$,
\item $m_{2k} \ge m_{2k+2}$ for $1 \le k \le \ell-1$,
\item the $(2k-1)$-signature of $Q$ is $(\lambda_k, 0)$ for $1 \le k \le \ell$,
\item the $(2k)$-signature of $r_{2k+1}Q$ is $(0, m_{2k}-m_{2k+2})$ for $1 \le k \le \ell-1$,
\item the $(2k)$-signature of $r_{2k-1}Q$ is $(\lambda_k-\lambda_{k+1}-p, m_{2k}-m_{2k+2}-p)$ with some $p \ge 0$ for $1 \le k \le \ell-1$.
\end{enumerate} 
Then we have

\begin{thm}{\cite[Theorem 6.12]{K15}}\label{thm:Ainsertion}
For $(\lambda, \ell) \in \mc{P}({\rm Sp})_\A$, $\Phi_\A$ induces a bijection
\begin{equation*} 
\xymatrixcolsep{4pc}\xymatrixrowsep{0.5pc}\xymatrix{
{\bf T}_\A(\la,\ell) \ar@{->}^{\hskip -10mm \Phi_\A}[r] & \  \displaystyle{\bigsqcup_{\mu}}\, SST_\A(\mu)\times {K}_{\mu (\lambda, \ell)}},
\end{equation*}
where the union is over $\mu\in \cP$ with $\ell(\mu')\leq 2\ell$.
\end{thm}

\begin{cor}
For $(\lambda, \ell) \in \mc{P}({\rm Sp})_\A$, we have
\begin{equation*}
S_{(\la,\ell)}({\bf x}_\A)  = t^\ell \sum_{\mu} c_{\mu(\la,\ell)}s_\mu({\bf x}_\A),
\end{equation*}
where $c_{\mu(\la,\ell)} = |{K}_{\mu(\la,\ell)}|$.
\end{cor}

\begin{ex}{\rm
Suppose that  $\mathcal{A} = \mathbb{I}_{4|3}$ and $(\la,\ell)=((3, 2, 1), 3)$.
\[ \begin{tikzpicture}
\node [above right] at (0, 0) {$\begin{ytableau}
\tl 2 & \tl 2 \\
\tl 3 & \tl 4 \\
\tl{$1'$}
\end{ytableau}$};
\node [above right] at (1.5, 0) {$\begin{ytableau}
\tl 1 & \tl 1 \\
\tl 3 & \tl{$2'$} \\
\tl{$1'$} \\
\tl{$2'$}
\end{ytableau}$};
\node [above right] at (3, 0) {$\begin{ytableau}
\tl 2 & \tl{$2'$} \\
\tl{$1'$} & \tl{$2'$} \\
\tl{$1'$} \\
\tl{$3'$} \\
\tl{$3'$}
\end{ytableau}$};

\node [above right] at (6.5, 0) {$\left(\ \raisebox{2em}{\begin{ytableau}
\tl 1 & \tl 1 & \tl 2 & \tl{$1'$} & \tl{$2'$} \\
\tl 2 & \tl 2 & \tl 3 & \tl{$2'$} & \tl{$3'$} \\
\tl 3 & \tl 4 & \tl{$1'$} & \tl{$2'$} \\
\tl{$1'$} & \tl{$2'$} & \tl{$3'$} \\
\tl{$1'$}
\end{ytableau}},\quad
\raisebox{2em}{\begin{ytableau}
\tl 1 & \tl 1 & \tl 2 & \tl 2 & \tl 2 \\
\tl 2 & \tl 2 & \tl 3 & \tl 4 \\
\tl 3 & \tl 4 & \tl 4 & \tl 6 \\
\tl 4 & \tl 5 & \tl 5 \\
\tl 6 & \tl 6
\end{ytableau}}\ \right)$};

\draw [dotted] (-0.2, 0.362em) -- (4.7, 0.362em);
\node [below] at (2.5, -0.2) {$\mathbf{T} \in \mathbf{T}_\mathcal{A}(\la,\ell)$};
\node [left] at (6.7, 3.6em) {$\Phi_\mathcal{A}(\mathbf{T}) =$};
\end{tikzpicture} \]
} \end{ex}


\section{Kashiwara-Nakashima tableaux and symplectic jeu de taquin}\label{sec:KN tableaux}
\subsection{KN tableaux of type $C_n$}\label{subsec:KN tableaux} 
Let us review the notion of Kashiwara-Nakashima tableaux (KN tableaux for short) of type $C_n$ \cite{KashNaka}.

Let $P=\bigoplus_{i=1}^n\Z\epsilon_i$, where $\{\,\epsilon_i\,|\,1\leq i\leq n\,\}$ is an orthonormal basis with respect to a symmetric bilinear form $(\,,\,)$. 
Suppose that $\g=\mf{sp}_{2n}$ of type $C_n$ with Dynkin diagram
\begin{center} 
\setlength{\unitlength}{0.19in}
\begin{picture}(20,3)
\put(5.6,2){\makebox(0,0)[c]{$\bigcirc$}}
\put(12.6,2){\makebox(0,0)[c]{$\bigcirc$}}
\put(10.4,2){\makebox(0,0)[c]{$\bigcirc$}}
\put(14.85,2){\makebox(0,0)[c]{$\bigcirc$}}
\put(6,2){\line(1,0){1.3}} \put(8.7,2){\line(1,0){1.3}} \put(10.82,2){\line(1,0){1.3}}
%
\put(13.7,2){\makebox(0,0)[c]{$\Longleftarrow$}}

\put(8,1.95){\makebox(0,0)[c]{$\cdots$}}
\put(5.6,1){\makebox(0,0)[c]{\tiny ${\alpha}_1$}}
\put(12.7,1){\makebox(0,0)[c]{\tiny ${\alpha}_{n-1}$}}
\put(10.4,1){\makebox(0,0)[c]{\tiny ${\alpha}_{n-2}$}}
\put(15,1){\makebox(0,0)[c]{\tiny ${\alpha}_n$}}

\end{picture}
\end{center}\vskip -3mm 
where $\alpha_i=\epsilon_i-\epsilon_{i+1}$ for $1\leq i\leq n-1$, and $\alpha_n=2\epsilon_n$.
The set of dominant integral weights is given by $P^+=\{\,\omega_\la\,|\,\la\in \cP_n\,\}$, where $\omega_\la=\la_1\epsilon_1+\cdots+\la_n\epsilon_n$. 

\begin{df}\label{def:KN-C}{\rm
For $\la\in \cP_n$, let ${\bf KN}_{\la}$ be the set of $T\in SST_{\mc{I}_n}(\la)$ satisfying  
\begin{itemize}
\item[(1)] if $T({i_1,j})=\ov{a}$ and $T({i_2,j})=a$ for some $a$ and $1\leq i_1< i_2\leq \la'_j$, then we have $i_1+(\la'_j-i_2+1)\leq a$,

\item[(2)] if either $T({p,j})=\ov{a}$, $T({q,j})=\ov{b}$, $T({r,j})=b$, $T({s,j+1})=a$ or 
 $T({p,j})=\ov{a}$, $T({q,j+1})=\ov{b}$, $T({r,j+1})=b$, $T({s,j+1})=a$ for some $1\leq a\leq b\leq n$, and $p\leq q<r\leq s$, then we have $(q-p)+(s-r)<b-a$,

\end{itemize}
where $T({i,j})$ denotes the entry of $T$ in the $i$th row from the bottom and the $j$th column from the right. We call ${\bf KN}_{\la}$ the set of {\em KN tableaux of type $C_n$ of shape $\la$.}}
\end{df}

The set ${\bf KN}_{(1)}$ has an $\mf{sp}_{2n}$-crystal structure such that ${\bf KN}_{(1)}\cong B(\epsilon_1)$, where
\begin{equation*}
  \xymatrixcolsep{2pc}\xymatrixrowsep{0pc}\xymatrix{
  \boxed{1}
  \, \ar@{->}[r]^{1} & \, \boxed{2} \, \ar@{->}[r]^{2} & \, \cdots \, \ar@{->}[r]^{n-1}  
  & \, \boxed{n}\,  \ar@{->}[r]^{n} & \, \boxed{\ov{n}}\, \ar@{->}[r]^{n-1} &   \cdots   \ar@{->}[r]^{2} & \boxed{\ov{2}}\, \ar@{->}[r]^{1} & \, \boxed{\ov{1}}}
\end{equation*}
with ${\rm wt}(\,
\resizebox{.02\hsize}{!}{\def\lr#1{\multicolumn{1}{|@{\hspace{.6ex}}c@{\hspace{.6ex}}|}{\raisebox{-.25ex}{$#1$}}}\raisebox{-.65ex}
{$\begin{array}[b]{c}
\cline{1-1} 
\lr{{i}} \\
\cline{1-1} 
\end{array}$}}
\,)=\epsilon_i$ and 
${\rm wt}(\,
\resizebox{.02\hsize}{!}{\def\lr#1{\multicolumn{1}{|@{\hspace{.6ex}}c@{\hspace{.6ex}}|}{\raisebox{-.25ex}{$#1$}}}\raisebox{-.65ex}
{$\begin{array}[b]{c}
\cline{1-1} 
\lr{\ov{i}} \\
\cline{1-1} 
\end{array}$}}
\,)=-\epsilon_i$ for $1\leq i\leq n$. Here $b\stackrel{i}{\rightarrow}b'$ means $\tf_ib=b'$.
For $\la\in \cP_n$ and $1 \le i \le n$, we define $\te_i$ and $\tf_i$ on ${\bf KN}_{\la}$ under the identification of $T\in {\bf KN}_{\la}$ with 
$\resizebox{.04\hsize}{!}{\def\lr#1{\multicolumn{1}{|@{\hspace{.6ex}}c@{\hspace{.6ex}}|}{\raisebox{-.1ex}{$#1$}}}\raisebox{-.9ex}
{$\begin{array}[b]{c}
\cline{1-1} 
\lr{w_1} \\
\cline{1-1} 
\end{array}$}}\otimes \cdots \otimes 
\resizebox{.04\hsize}{!}{\def\lr#1{\multicolumn{1}{|@{\hspace{.6ex}}c@{\hspace{.6ex}}|}{\raisebox{-.1ex}{$#1$}}}\raisebox{-.9ex}
{$\begin{array}[b]{c}
\cline{1-1} 
\lr{w_r} \\
\cline{1-1} 
\end{array}$}}\in ({\bf KN}_{(1)})^{\otimes r}$ when $w(T)=w_1\cdots w_r$. 
Then ${\bf KN}_{\la}$ is an $\mf{sp}_{2n}$-crystal with respect to $\te_i$ and $\tf_i$ for $1 \le i \le n$, and ${\bf KN}_{\la}\cong B(\omega_\la)$.

 
\subsection{Bijection between KN tableaux and spinor model} 
Suppose that $\A=[\ov{n}]$. Then we have $\mc{P}({\rm Sp})_n:=\mc{P}({\rm Sp})_{[\ov n]}=\{\,(\la,\ell)\in \mc{P}({\rm Sp}) \,|\,\la_1\leq n\,\}$. For simplicity, we put  
\begin{equation*}
\begin{split}
{\bf T}_n(a) & = {\bf T}_{[\ov{n}]}(a) \quad (0\leq a\leq n),\\
{\bf T}_n(\la,\ell) & ={\bf T}_{[\ov{n}]}(\la,\ell) \quad ((\la,\ell)\in \mc{P}({\rm Sp})_n).
\end{split}
\end{equation*}
For $(\la,\ell)\in \mc{P}({\rm Sp})_n$,  
put 
$$\rho_{n}(\la,\ell)=(n-\la_\ell,n-\la_{\ell-1},\dots,n-\la_1)',$$
which is the conjugate of the rectangular complement of $\la$ in $(n^\ell)$.

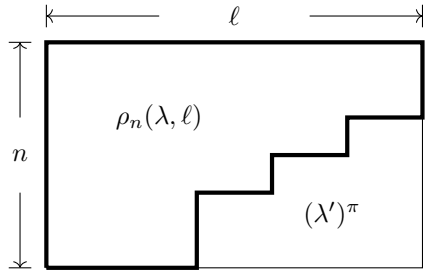
\[ \begin{tikzpicture}
\draw (0, 0) rectangle (5, 3);
\draw [ultra thick] (0, 0) -- (2, 0) -- (2, 1) -- (3, 1) -- (3, 1.5) -- (4, 1.5) -- (4, 2) -- (5, 2) -- (5, 3) -- (0, 3) -- (0, 0);
\node at (1.5, 2) {$\rho_n(\lambda, \ell)$};
\node at (3.8, 0.7) {$(\lambda')^\pi$};

\draw (-0.5, 3) -- (-0.2, 3); \draw (-0.5, 0) -- (-0.2, 0);
\draw [->] (-0.35, 2) -- (-0.35, 3); \draw [->] (-0.35, 1) -- (-0.35, 0); \node at (-0.35, 1.5) {$n$};
\draw (0, 3.2) -- (0, 3.5); \draw (5, 3.2) -- (5, 3.5);
\draw [->] (1.5, 3.35) -- (0, 3.35); \draw [->] (3.5, 3.35) -- (5, 3.35); \node at (2.5, 3.35) {$\ell$};
\end{tikzpicture} \]
 
For $U\in SST_{[\bar{n}]}((1^m))$, let $U^c$ be the tableau in $SST_{[n]}((1^{n-m}))$ such that $k$ appears in $U^c$ if and only if $\ov{k}$ does not appear in $U$ for each $k\in [n]$.
 
For $T\in {\bf T}_n(a)$, we define the following:
\begin{itemize}
\item ${T}^{\tt ad}$ : the tableau obtained by putting ${}^{\tt L}T$ below $({}^{\tt R}T)^c$,

\item ${T}^{\tt ad*}$ : the tableau obtained by putting ${T}^{\tt R}$ below $({T}^{\tt L})^c$.
\end{itemize}
Then the map $T \longmapsto {T}^{\tt ad}$ is a bijection from ${\bf T}_n(a)$ to ${\bf KN}_{(1^{n-a})}$ \cite[Lemma 3.11]{K18-2}. 

\begin{ex}{\rm Suppose that $n=5$ and $T \in \mathbf{T}_5(1)$ is given as follows. Then
\[ \begin{tikzpicture}
\node at (-1.5, 1.2) {\begin{ytableau} \bar{\tl 4} & \bar{\tl 2} \\ \bar{\tl 3} & \bar{\tl 1} \\ \bar{\tl 1} \end{ytableau}};
\node at (0.5, 1.2) {\begin{ytableau} \bar{\tl 3} & \bar{\tl 4} \\ \bar{\tl 1} & \bar{\tl 2} \\ \none & \bar{\tl 1} \end{ytableau}};

\node at (4, 0) {\begin{ytableau} \tl 2 \\ \tl 5 \end{ytableau}};
\node at (4, 2.5) {\begin{ytableau} \tl 3 \\ \tl 5 \end{ytableau}};

\node at (6.5, 0) {\begin{ytableau} \tl 2 \\ \tl 5 \\ \bar{\tl 2} \\ \bar{\tl 1} \end{ytableau}};
\node at (6.5, 2.5) {\begin{ytableau} \tl 3 \\ \tl 5 \\ \bar{\tl 3} \\ \bar{\tl 1} \end{ytableau}};

\draw [dotted] (-2.3, 0.5) -- (1.2, 0.5);

\node at (-2.5, 1.2) {$T=$};
\node [left] at (-1.4, 0.2) {$T^\texttt{L}$}; \node [right] at (-1.5, 0.2) {$T^\texttt{R}$};
\node [left] at (0.5, 0.2) {${^\texttt{L}T}$}; \node [right] at (0.5, 0.2) {${^\texttt{R}T}$};
\node at (3, 0) {$({T^\texttt{L}})^c=$}; \node at (5.5, 0) {$T^\texttt{ad*}=$};
\node at (3, 2.5) {$({^\texttt{R}T})^c = $}; \node at (5.5, 2.5) {$T^\texttt{ad} =$};
\end{tikzpicture}. \]
} \end{ex}

Let $C\in SST_{\mc{I}_n}((1^m))$ be given with $m\ge 1$. 
We call $C$ {\it admissible} (or an {\it admissible column}) if $C \in {\bf KN}_{(1^m)}$, equivalently $C=T^{\tt ad}$ for some $T\in {\bf T}_n(n-m)$, 
and {\it coadmissible} (or a {\it coadmissible column}) if $C=T^{\tt ad*}$ for some $T\in {\bf T}_n(n-m)$. We remark that the definition of an coadmissible column is equivalent to the one in \cite{Le02}.

Let $C\in SST_{\mc{I}_n}((1^m))$ be an admissible column with $C=T^{\tt ad}$ for $T\in {\bf T}_n(n-m)$. Then we define the following:
\begin{itemize}
\item  $lC$ : the tableau obtained by putting ${}^{\tt L}T$ below $(T^{\tt L})^c$, 

\item  $rC$ : the tableau obtained by putting ${T}^{\tt R}$ below $({}^{\tt R}T)^c$,

\item  $C^* = T^{\tt ad*}$. 
\end{itemize} 

For $V\in SST_{\mc{I}_n}((1^m))$, let $V_+$ and $V_-$ be the subtableau of $V$ 
consisting of its positive (unbarred) and negative (barred) letters, respectively.

\begin{lem}\label{lem:lCrC}
Under the above hypothesis, we have
\begin{equation*}
\begin{array}{ll}
(lC)_+= (T^{\tt L})^c = (T^{\tt ad*})_+, & (rC)_+= ({}^{\tt R}T)^c=(T^{\tt ad})_+,\\
(lC)_-= {}^{\tt L}T = (T^{\tt ad})_-,    & (rC)_-= {T}^{\tt R} = (T^{\tt ad*})_- .
\end{array}
\end{equation*}
\end{lem}

\begin{rem}{\rm
The notion of admissible and coadmissible columns has been introduced in \cite{De,Sh}.
For an admissible column $C$, one can check that $lC$ and $rC$ in this paper are equal to those in \cite[Definition 2.2.1]{Le02}.
}
\end{rem}

Let ${\rm spl}(C)=\llceil lC,rC \rrceil$ and call it the splitting form of $C$.
Note that $\mbox{spl}(C)$ is $\mc{I}_n$-semistandard.

Let $C_i\in SST_{\mc{I}_n}((1^{m_i}))$ ($i=1,2$) be an admissible column with $m_2\geq m_1$ and $C_i=T_i^{\tt ad}$ for $T_i\in {\bf T}_n(n-m_i)$. Following \cite{Le02}, define  
\begin{equation}\label{eq:ab-config}
\text{$C_2\prec C_1$ \quad if $\llceil rC_2,lC_1 \rrceil$ is $\mc{I}_n$-semistandard.}
\end{equation}
 
\begin{lem}\label{lem:equiv on admissibility}
Under the above hypothesis, we have $C_2\prec C_1$ if and only if $T_2\prec T_1$.
\end{lem} 
\pf It can be checked in a straightforward manner by using Lemma \ref{lem:lCrC} that Definition \ref{def:spinor tableaux}(1) is equivalent to \eqref{eq:ab-config}.
\qed\vskip 2mm

The following is another characterization of KN tableaux \cite[Theorem A.4]{Sh}.
\begin{prop} \label{prop:admissible}
Let $\mu\in \cP_n$ be given. 
Let $m_1,\dots,m_r$ be positive integers such that $\mu'=(m_r,\dots,m_1)$.
For $(C_r,\dots,C_1)\in SST_{\mc{I}_n}((1^{m_r})) \times \cdots \times SST_{\mc{I}_n}((1^{m_1}))$, we have
\begin{equation*}
\text{$\llceil C_r,\dots, C_1 \rrceil \in {\bf KN}_\mu$\quad if and only if\quad $C_{i+1}\prec C_i$ for $1\leq i\leq r-1$.}
\end{equation*}
\end{prop}

\begin{prop}\label{prop:bij spin to KN}
For $(\la,\ell)\in \mc{P}({\rm Sp})_n$, we have a bijection
\begin{equation}\label{eq:Psi_lambda:C}
\xymatrixcolsep{2pc}\xymatrixrowsep{0.5pc}\xymatrix{
 {\bf T}_n(\la,\ell)  \ar@{->}[r]  & \ {\bf KN}_{\rho_{n}(\la,\ell)} \\
  {\bf T}=(T_\ell,\dots,T_1)  \ar@{|->}[r] & {\bf T}^{\tt ad}:=\llceil T_\ell^{\tt ad},\dots, T_1^{\tt ad} \rrceil }.
\end{equation}
\end{prop}
\pf
Let ${\bf T}=(T_\ell,\dots,T_1)\in {\bf T}_{n}(\lambda, \ell)$ be given.
Put $C_i=T^{\tt ad}_i$ for $1\leq i\leq \ell$. Then 
$${\bf T}^{\tt ad}=\llceil C_\ell,\dots,C_1\rrceil$$ 
is a tableau of shape $\rho_{n}(\la,\ell)$. 
By Proposition \ref{prop:admissible}, ${\bf T}^{\tt ad}\in {\bf KN}_{\rho_{n}(\la,\ell)}$ 
if and only if $C_{i+1}\prec C_i$ for $1\leq i\leq \ell-1$. 
Hence by Lemma \ref{lem:equiv on admissibility}, we have ${\bf T}\in {\bf T}_n(\la,\ell)$ if and only if $T^{\tt ad} \in {\bf KN}_{\rho_{n}(\la,\ell)}$. The proof completes.
\qed

\begin{rem}{\rm
The bijectivity of the map in Proposition \ref{prop:bij spin to KN} is also proved in \cite[Theorem 3.14]{K18-2} by constructing an isomorphism of $\mathfrak{sp}_{2n}$-crystals between them. But we do not use the crystal structures on ${\bf T}_n(\la,\ell)$ and ${\bf KN}_{\rho_{n}(\la,\ell)}$ in this paper.
}
\end{rem}


\subsection{Admissibility and $\mathfrak{sl}_2$-crystal} 

In this subsection, we describe the relation between semistandard column tableaux with letters in $\mc{I}_n$ and admissible column tableaux using the crystal operators $\mathcal{E}$ and $\mathcal{F}$ \eqref{eq:E and F via jdt}.

Let
\begin{equation}\label{eq:single column tab set}
{\bf F}_n=\bigsqcup_{0 \le m \le 2n}SST_{\mc{I}_n}((1^m)),
\end{equation}
where we assume that $SST_{\mc{I}_n}((1^m))$ is the set of the empty tableau when $m=0$.

Let $C\in SST_{\mc{I}_n}((1^m))$ given. Recall that $(C_+)^c$ is the tableau in $SST_{[\ov{n}]}(1^{n-l})$ such that $\bar{k}$ appears in $(C_+)^c$ if and only if $k$ does not appear in $C_+$ for each $k\in [n]$, where $C_+\in SST_{[n]}((1^l))$.  

We define $\mc{E}C$ to be the unique tableau $C' \in SST_{\mc{I}_n}((1^{m-2}))$ such that 
\begin{equation*}
(C'_-,(C'_+)^c)=\mc{E}(C_-,(C_+)^c),
\end{equation*}
when $\mc{E}(C_-,(C_+)^c)\neq {\bf 0}$, and $\mc{E}C={\bf 0}$ otherwise (see \eqref{eq:E and F via jdt}). 
We define $\mc{F}C$ in a similar way.

\begin{lem}\label{lem:admissibility of a column}
Under the above hypothesis, 
\begin{itemize}
\item[(1)] ${\bf F}_n$ is a regular $\mf{sl}_2$-crystal with respect to $\mc{E}$ and $\mc{F}$,

\item[(2)] $C\in {\bf F}_n$ is admissible if and only if $\mc{E}C=0$,

\item[(3)] we have a bijection
\begin{equation} \label{eq:sl2 bijection}
\xymatrixcolsep{2pc}\xymatrixrowsep{0.5pc}\xymatrix{
{\bf F}_n \ar@{->}[r]  & \ 
\displaystyle{\bigsqcup_{0\leq a\leq n} 
{\bf KN}_{(1^{n-a})}\times \Z/(a+1)\Z } \\
C  \ar@{|->}[r] & (T,\varepsilon(C))},
\end{equation}
where $\varepsilon(C)=\max\{\,k\,|\,\mc{E}^kC\neq {\bf 0}\,\}$, 
$T=\mc{E}^{\rm max}C=\mc{E}^{\varepsilon(C)}C$, 
and $\mathbb{Z}/(a+1)\mathbb{Z}$ is understood as the set $\{ 0, 1, \dots, a \}$.

\end{itemize}
\end{lem}
\pf (1) It is clear that ${\bf F}_n$ is a regular $\mf{sl}_2$-crystal.
(2) Let $C\in SST_{\mc{I}_n}((1^m))$ be given with $\varepsilon=\varepsilon(C)$. We may identify $(C_-,(C_+)^c)$ as a skew tableau $T$ of shape $\la(\varepsilon,b,c)$ for some $b,c\in\Z_+$ with $\mf{r}_T=0$. 
Let $\ov{z}$ be a letter in $C$ or $T$. Then $\boxed{\ov{z}}$ can be moved to the right in $T$ by jeu de taquin or $\mc{E}$ if and only if $z$ is the smallest letter in $C$ which does not satisfy Definition \ref{def:KN-C}(1).
This implies that $C$ is admissible if and only if $\varepsilon=0$. 

(3) Let $C$ be given. By (2), $\mc{E}^{\rm max}C$ is admissible. If we let $C'=\mc{E}^{\varepsilon}C$, then 
$(C'_-,(C'_+)^c)$ is a skew tableau of shape $\la(0,a,c)$ for some $a,c\in\Z_+$,
and hence $C'\in {\bf KN}_{(1^{n-a})}$ by Proposition \ref{prop:bij spin to KN}. Since ${\bf F}_n$ is a regular $\mf{sl}_2$-crystal, we have $\varepsilon(C)\leq a$ and $\varphi(C)+\varepsilon(C)=a$.
The map sending $C$ to $(\mc{E}^{\rm max}C,\varepsilon(C))$ is reversible, and hence it gives the bijection.
\qed

\begin{ex}{\rm For $C \in \mathbf{F}_5$ below, we have
\[ C \ =\ \raisebox{3em}{$\begin{ytableau}
\tl 1 \\
\tl 2 \\
\tl 3 \\
\tl 4 \\
\bar{\tl 5} \\
\bar{\tl 4}
\end{ytableau}$}
\quad , \quad
(C_-, (C_+)^c) \ = \ \raisebox{0.5em}{$\begin{ytableau}
\bar{\tl 5} & \bar{\tl 5} \\
\bar{\tl 4}
\end{ytableau}$}
\quad ,  \quad
\mathcal{E}(C_-, (C_+)^c) \ = \ \raisebox{0.5em}{$\begin{ytableau}
\none & \bar{\tl 5} \\
\bar{\tl 5} & \bar{\tl 4}
\end{ytableau}$}. \]
Hence the image of $C$ under \eqref{eq:sl2 bijection} is
\[ (T, \varepsilon(C)) = \left(\ \raisebox{1.5em}{$\begin{ytableau}
\tl{$1$} \\
\tl{$2$} \\
\tl{$3$} \\
\bar{\tl 5}
\end{ytableau}$}\ , \
1 \right). \]
}\end{ex}
\vskip 2mm
 

\subsection{Symplectic insertion} 
Let us review the insertion algorithm for KN tableaux of type $C$ in \cite{Le02}.

Let $Pl(C_n)$ be the quotient of the free monoid generated by $\mc{I}_n$ subject to the relations:
\begin{enumerate} 
\item $yzx = yxz$ for $x \le y < z$ with $z \ne \bar{x}$,

\item $xzy = zxy$ for $x<y \le z$ with $z \ne \bar{x}$,
\item $y\overline{(x-1)}(x-1) = yx\bar{x}$ for $1 < x \le n$ and $x \le y \le \bar{x}$,

\item $x\bar{x}y = \overline{(x-1)}(x-1)y$ for $1 < x \le n$ and $x \le y \le \bar{x}$,

\item if $w=w(C)$ for a non-admissible column $C$ such that every proper subword is an admissible column word (the word of an admissible column), 
and $z$ is the smallest letter in $[n]$ such that the pair $(z, \bar{z})$ occurs in $w$ with $N(z)>z$, where $N(z)$ is the number of letters $x$ in $C$ such that $x\leq z$ or $x \geq \ov{z}$, then
\[ w = \widetilde{w}, \]
where $\widetilde{w}$ is the column word obtained by removing the pair $(z, \bar{z})$ in $w$.
\end{enumerate}

Denote by $\mc{W}_{\mc{I}_n}$ the set of words with letters in $\mc{I}_n$. For $w, w'\in \mc{W}_{\mc{I}_n}$, write $w\Cequiv w'$ if $w=w'$ in $Pl(C_n)$.

\begin{prop}[\cite{Le02}]\label{prop:plactic for C} 
For $w \in \mc{W}_{\mc{I}_n}$, we have a unique KN tableau $T$ such that $w \Cequiv w(T)$, which we denote by $P(w)$.
\end{prop}

For $w \in \mc{W}_{\mc{I}_n}$, we may obtain $P(w)$ by insertion algorithm.
For $x \in \mc{I}_n$ and $T\in {\bf KN}_\la$, 
let us first define $x\rightarrow T$, {\em the insertion of $x$ into $T$} as follows. 
Assume that $T=\llceil C_r,\dots, C_1 \rrceil$.

\begin{enumerate}
\item[{\em Case 1.}]
Suppose that $w(C_r)x$ is the reading word of an admissible column $C_r'$. 
Then
\[ (x \to T) := \llceil C_r', C_{r-1}, \dots, C_1\rrceil. \]

\item[{\em Case 2.}]
Suppose that $w(C_r)x$ is not the word of a column tableau and then there exists a letter $x'$ such that $w(C_r)x \Cequiv x'w(C_r')$ for some admissible column $C_r'$. 
Then we consider the insertion of $x'$ into $T' = \llceil C_{r-1}, \dots, C_1\rrceil$. If it belongs to {\em Case 1}, then we have $(x'\to T')$ and let 
\begin{equation}\label{eq:insertion case 2}
(x \to T) = \llceil C_r', (x' \to T') \rrceil.
\end{equation} 
Otherwise repeat the above step until we get to {\em Case 1}.

\item[{\em Case 3.}] Suppose that $w(C_r)x$ is the word of a non-admissible column whose proper subwords are admissible and suppose $\widetilde{w(C_r)x} = y_1 \cdots y_s$. Then
\begin{equation}\label{eq:insertion case 3}
(x \to T) := (y_s \to ( y_{s-1} \to (\cdots \to( y_1 \to T') \cdots))),
\end{equation} 
where $T' = \llceil C_{r-1}, \dots, C_1\rrceil$.
\end{enumerate}
We remark that {\em Case 3} does not occur during the insertion on the right-hand side of \eqref{eq:insertion case 2} and \eqref{eq:insertion case 3}.

Now, let $w  = x_1 \cdots x_r \in \mc{W}_{\mc{I}_n}$ be given.
Define
\begin{equation*}
{\tt P}(w)= (x_r \to ( x_{r-1} \to (\cdots \to( x_2 \to \boxed{x_1}) \cdots))),
\end{equation*}
and define ${\tt Q}(w)$ to be the sequence $(Q_1, \dots, Q_r)$ of partitions, 
where $Q_i$ is the shape of ${\tt P}(x_1 \cdots x_i)$ for $1\leq i\leq  r$. 
Note that ${\tt Q}(w)$ is an $n$-oscillating tableau of shape $\mu={\rm sh}({\tt P}(w))$, that is, a sequence of partitions in $\cP_n$ such that $Q_r=\mu$ and each pair $(Q_i,Q_{i+1})$ differs by one box for $1\leq i\leq r-1$, i.e., $Q_{i+1}/Q_i = \ytableausetup{smalltableaux} \ydiagram{1}$ or $Q_i/Q_{i+1} = \ydiagram{1} \ytableausetup{nosmalltableaux}$.

\begin{prop}[\cite{Le02}]\label{prop:RSK for KN}
We have a bijection 
\begin{equation*} 
\xymatrixcolsep{4pc}\xymatrixrowsep{0.5pc}\xymatrix{
\W_{\mc{I}_n} \ar@{->}[r] & \  \displaystyle{\bigsqcup_{\mu\in\cP_n}} {\bf KN}_\mu \times OT_n(\mu) \\
w  \ar@{|->}[r] & ({\tt P}(w), {\tt Q}(w))},
\end{equation*}
where $OT_n(\mu)$ is the set of $n$-oscillating tableaux of shape $\mu$. Moreover, the map is an isomorphism of $\mf{sp}_{2n}$-crystals, where the operators $\te_i$ and $\tf_i$ act on the first component on the right-hand side.
\end{prop}

\begin{ex} \label{ex:Lecouvey bijection}{\rm
For $w = 1\bar{5}\bar{4} 2\bar{5}\bar{4}\bar{3} 14\bar{4}\bar{2} \in \mathcal{W}_{\mathcal{I}_5}$, we have
\[ 
  \texttt{P}(w) = \raisebox{1.2em}{\ $\begin{ytableau}
  \tl{$1$} & \tl{$1$} & \bar{\tl 5} \\
  \tl{$5$} & \bar{\tl 5} & \bar{\tl 4} \\
  \bar{\tl 5} & \bar{\tl 3} \\
  \bar{\tl 4}
  \end{ytableau}$},
\]
\[ \ytableausetup{smalltableaux}
\texttt{Q}(w) = \left( \raisebox{1.5em}{
\ydiagram{1, 0, 0, 0, 0}\ , \ydiagram{1, 1, 0, 0, 0}\ , \ydiagram{1, 1, 1, 0, 0}\ ,
\ydiagram{2, 1, 1, 0, 0}\ , \ydiagram{2, 2, 1, 0, 0}\ , \ydiagram{2, 2, 1, 1, 0}\ , \ydiagram{2, 2, 1, 1, 1}\ ,
\ydiagram{3, 2, 1, 1, 1}\ , \ydiagram{3, 3, 1, 1, 1}\ , \ydiagram{3, 3, 2, 1, 1}\ , \ydiagram{3, 3, 2, 1, 0}\
} \right). \]
} \end{ex}


\subsection{Symplectic jeu de taquin}\label{subsec:sjdt for KN} 
Let us briefly recall the algorithm of {\em symplectic jeu de taquin} for KN tableaux of type $C$ introduced in \cite{Sh}. Our review is based on \cite{Le02}.

Let $T$ be a tableau of skew shape with letters in $\mc{I}_n$.
We call $T$ admissible if its columns are admissible columns, and its splitting form ${\rm spl}(T)$ is $\mc{I}_n$-semistandard, where ${\rm spl}(T)$ is the tableau obtained by replacing each column $T_i$ with ${\rm spl}(T_i)$. Let us denote by ${\bf KN}_{\zeta/\eta}$ the set of admissible tableaux of shape $\zeta/\eta$.

Let us say that $T$ is punctured if some boxes in $T$ are removed. 
We put $\bullet$ in each position of removed boxes and regard them as entries in $T$ 
to keep track of the empty box in the process of sliding. 
We also say that a punctured tableau $T$ is admissible if its columns are admissible and its splitting form is $\mc{I}_n$-semistandard when we ignore the punctures $\bullet$. 
Here the punctures $\bullet$ are also duplicated in the splitting form.

Now, the symplectic jeu de taquin can be described in the following steps.\vskip 2mm

\noindent {\em Step 1}. Suppose that $T$ is a punctured admissible tableau of skew shape with two columns, say $T=\llceil C_2,C_1 \rrceil^{(c_2, c_1)}$. Suppose that $\bullet$ is in the left column and ${\rm spl}(T)$ is as follows:
\begin{equation}\label{eq:punctured}
\ytableausetup{boxsize=1.5em} 
\begin{ytableau} 
\none[\vdots] & \none[\vdots] & \none[\vdots] & \none[\vdots] \\ 
\mbox{$\bullet$} & \mbox{$\bullet$} & b & b' \\ 
a & a' & \none[\raisebox{-0.4em}{\vdots}] & \none[\raisebox{-0.4em}{\vdots}] \\ 
\none[\raisebox{-0.4em}{\vdots}] & \none[\raisebox{-0.4em}{\vdots}]
\end{ytableau} 
\end{equation}
where the boxes filled with $a, a'$ and $b, b'$ may be empty.
Let $T'$ be the tableau obtained by applying the following process.
\begin{enumerate}
\item Suppose that $a' \le b$ or the domino
\raisebox{-0.45em}{
$
\ytableausetup{boxsize=1.3em} 
\begin{ytableau} b & b' 
\end{ytableau}
$}
is empty. Then $T'$ is such that ${\rm spl}(T')$ is given by switching
\raisebox{-0.45em}{
$
\ytableausetup{boxsize=1.3em} 
\begin{ytableau} 
\bullet & \bullet
\end{ytableau}
$} and
\raisebox{-0.45em}{
$
\ytableausetup{boxsize=1.3em} 
\begin{ytableau} 
a & a' 
\end{ytableau}
$}. 

\item Suppose that $a'>b$ or the domino
\raisebox{-0.45em}{
$
\ytableausetup{boxsize=1.3em} 
\begin{ytableau} a & a' 
\end{ytableau}
$} 
is empty.

\begin{enumerate} 

\item[(2-a)] If $b\in [n]$, then we have a column tableau $C_2'$ given by exchanging 
\raisebox{-0.45em}{$
\ytableausetup{boxsize=1.3em} 
\begin{ytableau} \bullet
\end{ytableau}
$}  
with 
\raisebox{-0.45em}{$
\ytableausetup{boxsize=1.3em} 
\begin{ytableau} b
\end{ytableau}
$} in $C_2$, and a unique punctured admissible column $C_1'$ such that $(C_1')^\ast$
is given by exchanging 
\raisebox{-0.45em}{$
\ytableausetup{boxsize=1.3em} 
\begin{ytableau} b
\end{ytableau}
$} 
with 
\raisebox{-0.45em}{$
\ytableausetup{boxsize=1.3em} 
\begin{ytableau} \bullet
\end{ytableau}
$} 
in $C_1^\ast$. 
Note that $C_2'$ may not be admissible.

\item[(2-b)] If $b\in [\ov{n}]$, then we have a unique admissible column $C'_2$ such that $(C_2')^\ast$
is given by exchanging 
\raisebox{-0.45em}{$
\ytableausetup{boxsize=1.3em} 
\begin{ytableau} \bullet
\end{ytableau}
$}  
with 
\raisebox{-0.45em}{$
\ytableausetup{boxsize=1.3em} 
\begin{ytableau} b
\end{ytableau}
$}  
in $(C_2)^\ast$, and a unique punctured admissible column $C_1'$ given by exchanging 
\raisebox{-0.45em}{$
\ytableausetup{boxsize=1.3em} 
\begin{ytableau} b
\end{ytableau}
$} 
with 
\raisebox{-0.45em}{$
\ytableausetup{boxsize=1.3em} 
\begin{ytableau} \bullet
\end{ytableau}
$} 
in $C_1$.
\end{enumerate}
Then we put  
$$T'=\llceil C_2',C'_1 \rrceil^{(c_2-1, c_1)}.$$
\end{enumerate}

\begin{ex}{\rm
$$ 
T_1 = 
\raisebox{2em}{
\ytableausetup{boxsize = 1.2em}
\begin{ytableau}
\bullet & \tl 3 \\
\tl 3 & \bar{\tl 5} \\ 
\bar{\tl 5} & \bar{\tl 3} \\ 
\bar{\tl 4} \\ 
\bar{\tl 2} 
\end{ytableau}} \quad \quad
{\rm spl}(T_1) = 
\raisebox{2em}{
\begin{ytableau} 
\bullet & \color{red} \bullet & \color{blue} \tl 2 & \tl 3 \\ 
\tl 1 & \color{blue} \tl 3 & \bar{\tl 5} & \color{red} \bar{\tl 5} \\ 
\color{blue} \bar{\tl 5} & \bar{\tl 5} & \bar{\tl 3} & \color{red} \bar{\tl 2} \\ 
\color{blue} \bar{\tl 4} & \bar{\tl 4} \\ 
\color{blue} \bar{\tl 2} & \bar{\tl 2} 
\end{ytableau}}\quad \quad
T_1' = 
\raisebox{2em}{
\begin{ytableau} 
\color{blue} \tl 2 & \color{red} \bullet \\ 
\color{blue} \tl 3 & \color{red} \bar{\tl 5} \\ 
\color{blue} \bar{\tl 5} & \color{red} \bar{\tl 2} \\ 
\color{blue} \bar{\tl 4} \\ 
\color{blue} \bar{\tl 2} 
\end{ytableau}} = 
\raisebox{2em}{
\begin{ytableau} 
\tl 2 & \bar{\tl 5} \\ 
\tl 3 & \bar{\tl 2} \\ 
\bar{\tl 5} \\ 
\bar{\tl 4} \\ 
\bar{\tl 2} 
\end{ytableau}}
$$ 

$$
T_2 = 
\raisebox{3em}{
\begin{ytableau}
\bullet & \tl 2 \\
\bullet & \tl 3 \\
\tl 3 & \bar{\tl 5} \\
\tl 5 & \bar{\tl 4} \\
\bar{\tl 3} & \bar{\tl 2} \\
\bar{\tl 1}
\end{ytableau}} \quad\quad 
{\rm spl}(T_2) = 
\raisebox{3em}{
\begin{ytableau}
\bullet & \bullet & \tl 1 & \color{blue} \tl 2 \\
\color{red} \tl 2 & \tl 3 & \tl 3 & \color{blue} \tl 3 \\
\color{red} \tl 5 & \tl 5 & \color{blue} \bar{\tl 5} & \bar{\tl 5} \\
\bullet & \color{blue} \bullet & \color{red} \bar{\tl 4} & \bar{\tl 4} \\
\bar{\tl 3} & \color{red} \bar{\tl 2} & \color{blue} \bar{\tl 2} & \bar{\tl 1} \\
\bar{\tl 1} & \color{red} \bar{\tl 1}
\end{ytableau} } \quad\quad 
T_2' = 
\raisebox{3em}{
\begin{ytableau} 
\bullet & \color{blue} \tl 2 \\
\color{red} \tl 3 & \color{blue} \tl 3 \\ 
\color{red} \tl 5 & \color{blue} \bar{\tl 5} \\ 
\color{red} \bar{\tl 4} & \color{blue} \bullet \\ 
\color{red} \bar{\tl 3} & \color{blue} \bar{\tl 2} \\ 
\color{red} \bar{\tl 1} \end{ytableau}} =
\raisebox{3em}{
\begin{ytableau} 
\bullet & \tl 2 \\
\tl 3 & \tl 3 \\ 
\tl 5 & \bar{\tl 5} \\ 
\bar{\tl 4} & \bar{\tl 2} \\ 
\bar{\tl 3} \\ 
\bar{\tl 1} \end{ytableau}}$$
}
\end{ex}

\noindent {\em Step 2}. Let $T$ be an admissible tableau of skew shape, say 
$$T=\llceil C_r,\dots, C_1 \rrceil^{(c_r,\dots,c_1)}.$$ Let $c$ be an inner corner at the 
$i$-th column (from the right). 
Let $C_i^\bullet$ denote the column $C_i$ with 
\raisebox{-0.45em}{$
\ytableausetup{boxsize=1.3em} 
\begin{ytableau} \bullet
\end{ytableau}
$} 
placed at the top. 

We apply {\em Step 1} to
$\llceil C_i^{\bullet},C_{i-1}\rrceil^{(c_{i}-1, c_{i-1})}$ to have 
$\llceil C'_i,C'_{i-1}\rrceil^{(c'_{i},c'_{i-1})}$ where 
\raisebox{-0.45em}{$
\ytableausetup{boxsize=1.3em} 
\begin{ytableau} \bullet
\end{ytableau}
$} has moved to $C'_{i-1}$. Now, we repeat {\em Step 1} until 
\raisebox{-0.45em}{$
\ytableausetup{boxsize=1.3em} 
\begin{ytableau} \bullet
\end{ytableau}
$}
is placed at the bottom of the $j$-th column to have
$$T'=\llceil \dots, C_{i+1}, C'_i,\dots, C'_j,C_{j-1},\dots \rrceil^{(\dots,c'_i,\dots,c'_j,\dots)},$$
for some $j \le i$, $C'_i,\dots, C'_j$ and $c'_i,\dots,c'_j$, where we ignore 
\raisebox{-0.45em}{$
\ytableausetup{boxsize=1.3em} 
\begin{ytableau} \bullet
\end{ytableau}
$} in the $j$-th column. 
It is shown in \cite{Sh} that if $T'$ is not admissible, then $C'_i$ is the only column which is not admissible. 
Let $C''_i$ be the unique admissible column such that $w(C''_i)\Cequiv w(C'_i)$ with respect to $\mc{I}_{n}$. Put
$$T''=\llceil \dots, C_{i+1}, C''_i,\dots, C'_j,C_{j-1},\dots \rrceil^{(\dots,c_i,\dots,c'_j,\dots)}.$$
In this case, we have $c'_i=c_i-1$.
Now we define
\begin{equation}\label{eq:SJDT for KN}
{{\tt jdt}}_{KN}(T, c) = 
\begin{cases} 
T' & \mbox{if $T'$ is admissible}, \\ 
T'' & \mbox{if $T'$ has a non-admissible column.} 
\end{cases}
\end{equation}
If ${\rm sh}(T)=\zeta/\eta$, then
\begin{equation*}
{\rm sh}({{\tt jdt}}_{KN}(T, c))=
\begin{cases}
\alpha/\eta  & \text{for some $\alpha \subsetneq \zeta$ if $T'$ is non-admissible},\\
\alpha/\beta  & \text{for some $\alpha \subsetneq \zeta, \beta \subsetneq \eta$ if $T'$ is admissible}.
\end{cases}
\end{equation*}
Moreover, by \cite[Theorem 6.3.8]{Le02}, we have 
$$w({{\tt jdt}}_{KN}(T, c))\Cequiv w(T).$$
Hence by applying the ${{\tt jdt}}_{KN}(\,\cdot\,,c)$ successively to the inner corners, we obtain a unique KN tableau ${\tt P}(w(T))$ by induction on $|\zeta|=\sum_i\zeta_i$.


\section{Sliding algorithm for spinor model}\label{sec:sliding for spinor}

\subsection{Spinor model of a skew shape} 

\begin{df}
{\rm 
Let
${\bf T}=(T_\ell,\dots,T_1)\in {\bf T}_\A(a_\ell)\times\cdots\times{\bf T}_\A(a_1)$ be given for some $a_1,\dots,a_\ell\in \Z_+$. 
Let $\la/\mu$ be a skew diagram with $\la,\mu\in\cP_\ell$. We say that
\begin{itemize}
\item[(1)] ${\bf T}$ is of shape $\la/\mu$ if  
\begin{equation*}
a_i=\la_i-\mu_i,\quad {\rm ht}(T^{\tt L}_{i+1})+\mu_{i+1}\leq {\rm ht}(T^{\tt L}_i)+\mu_i\quad (1\leq i\leq \ell),
\end{equation*} 
where we put ${\rm ht}(C) = c$ for a tableau $C$ of shape $(1^c)$,

\item[(2)] ${\bf T}$ is $\A$-admissible of shape $\la/\mu$ if ${\bf T}$ is of shape $\la/\mu$ and for $1\leq i\leq \ell-1$
\begin{equation*}
\text{
$\llfloor {}^{\tt R}T_{i+1}, {T}^{\tt L}_i \rrfloor_{(\mu_{i+1},\mu_i)}$ 
and 
$\llfloor T_{i+1}^{\tt R},{}^{\tt L}T_i \rrfloor_{(\la_{i+1},\la_i)}$ 
are $\A$-semistandard.}
\end{equation*}
We denote by ${\bf T}_\A(\la/\mu,\ell)$ the set of $\A$-admissible tableaux of shape $\la/\mu$.
\end{itemize}}\vskip 2mm
\end{df}

When ${\bf T}$ is of shape $\la/\mu$, let us often identify ${\bf T}$ with a tableau in $\mathbb{P}_L$ given by
\begin{equation*}
\llfloor {\bf T} \rrfloor_{(\mu_\ell,\dots,\mu_1)}= \llfloor T_\ell,\dots,T_1 \rrfloor_{(\mu_\ell,\dots,\mu_1)} := 
\llfloor\, T^{\tt L}_\ell,T^{\tt R}_\ell,\dots,T^{\tt L}_1,T^{\tt R}_1 \,\rrfloor_{(\mu_\ell,\lambda_\ell,\dots,\mu_1,\lambda_1)}.
\end{equation*}
Note that ${\bf T}$ does not correspond to a KN tableau of skew shape unless $\A = [\overline{n}]$. 
So we may not apply the algorithm in Section \ref{subsec:sjdt for KN} to have ${\bf T}'\in {\bf T}_\A(\nu,\ell)$ for some $(\nu,\ell)\in \mc{P}({\rm Sp})_\A$. 
To overcome this problem, we introduce the notion of {\it $n$-conjugate} of ${\bf T}$.

\begin{df}{\rm
Let ${\bf T}=(T_\ell,\dots,T_1)\in {\bf T}_\A(a_\ell)\times\cdots\times{\bf T}_\A(a_1)$ be given. 
Let $(P({\bf T}), Q({\bf T}))=\kappa_\A({\bf T})=\kappa_\mathcal{A}(T_\ell^{\tt L}, T_\ell^{\tt R}, \dots, T_1^{\tt L}, T_1^{\tt R})$ with $\nu={\rm sh}(P({\bf T}))$, where $\kappa_\A$ is in \eqref{eq:RSK}.

For $n\geq \ell(\nu)$, we define the $n$-conjugate of ${\bf T}$ to be the unique 
$\ov{\bf T}=(\ov T_\ell,\dots,\ov T_1)\in {\bf T}_{[\ov{n}]}(a_\ell)\times\cdots\times{\bf T}_{[\ov{n}]}(a_1)$ such that
\begin{equation*}
\kappa_{[\ov{n}]}\left(\ov{\bf T}\right) = (P(\ov{\bf T}), Q(\ov{\bf T})) =(H_{\nu}, Q({\bf T})),
\end{equation*}
where $H_\nu\in SST_{[\ov{n}]}(\nu)$ is the highest weight element, that is, the $i$-th row from the top is filled with $\ov{n-i+1}$ for $1\leq i\leq n$. 
}
\end{df}

Note that if we replace $n$ with $m> n$, then the corresponding $\ov{\bf T}$ is given by replacing $\ov{a}$ with $\ov{a+m-n}$ for all $a$.

\begin{ex} \label{ex:conjugate}{\rm
Suppose $\mathcal{A} = \mathbb{I}_{4|3}$ and take $\mathbf{T} = (T_3, T_2, T_1) \in \mathbf{T}_\mathcal{A}(1) \times \mathbf{T}_\mathcal{A}(1) \times \mathbf{T}_\mathcal{A}(2)$ below. 
Then we have
\[
\mathbf{T} = \left(\
 \raisebox{1.5em}{$\ytableausetup{boxsize=1.2em} \begin{ytableau}
\none \\
\tl 2 & \tl 4 \\
\tl 3 & \tl{$2'$} \\
\tl{$1'$}
\end{ytableau}$ }, \
\raisebox{1.5em}{$\begin{ytableau}
\tl 1 & \tl 1 \\
\tl 2 & \tl 3 \\
\tl{$1'$} & \tl{$2'$} \\
\tl{$3'$}
\end{ytableau}$ }, \
\raisebox{1.5em}{$\begin{ytableau}
\tl 2 & \tl{$2'$} \\
\tl{$1'$} & \tl{$2'$} \\
\tl{$1'$} \\
\tl{$3'$}
\end{ytableau}$}\ \right)
\ \ , \ \
\kappa_\mathcal{A}(T) =
\left(\ \raisebox{2em}{\begin{ytableau}
\tl 1 & \tl 1 & \tl 2 & \tl{$1'$} & \tl{$2'$} \\
\tl 2 & \tl 2 & \tl 3 & \tl{$2'$} & \tl{$3'$} \\
\tl 3 & \tl 4 & \tl{$1'$} & \tl{$2'$} \\
\tl{$1'$} & \tl{$2'$} & \tl{$3'$} \\
\tl{$1'$}
\end{ytableau}},\quad
\raisebox{2em}{\begin{ytableau}
\tl 1 & \tl 1 & \tl 2 & \tl 2 & \tl 4 \\
\tl 2 & \tl 2 & \tl 3 & \tl 5 \\
\tl 3 & \tl 3 & \tl 4 & \tl 6 \\
\tl 4 & \tl 4 & \tl 5 \\
\tl 6 & \tl 6
\end{ytableau}}\ \right).
\]
Hence the $5$-conjugate $\ov{\mathbf{T}}$, the inverse image of $(H_\nu, Q(\mathbf{T}))$ under $\kappa_{[\ov{5}]}$ with $\nu = (5, 5, 4, 3, 1)$, is
\[ \overline{\mathbf{T}} = \left( \raisebox{1.5em}{
\begin{ytableau}
\none \\
\bar{\tl 5} & \bar{\tl 3} \\
\bar{\tl 4} & \bar{\tl 2} \\
\bar{\tl 2}
\end{ytableau}}\ , \
\raisebox{1.5em}{
\begin{ytableau}
\bar{\tl 5} & \bar{\tl 5} \\
\bar{\tl 4} & \bar{\tl 4} \\
\bar{\tl 3} & \bar{\tl 3} \\
\bar{\tl 1}
\end{ytableau}}\ , \
\raisebox{1.5em}{
\begin{ytableau}
\bar{\tl 5} & \bar{\tl 5} \\
\bar{\tl 4} & \bar{\tl 4} \\
\bar{\tl 3} \\
\bar{\tl 2}
\end{ytableau} }
\right). \]
}\end{ex}

The following two lemmas play an important role in this paper.
\begin{lem}\label{lem:admin between skew tableaux-1}
Under the above hypothesis, ${\bf T}$ is $\A$-admissible of shape $\la/\mu$ if and only if its $n$-conjugate $\ov{\bf T}$ is $[\ov{n}]$-admissible of shape $\la/\mu$.
\end{lem}
\pf It follows directly from \cite[Lemma 6.2]{K15} and $Q({\bf T})=Q(\ov{\bf T})$.
\qed

\begin{lem}\label{lem:admin between skew tableaux-2}
Under the above hypothesis, $\ov{\bf T}$ is $[\ov{n}]$-admissible of shape $\la/\mu$ if and only if 
$$\llceil\, \ov{T}_\ell^{\tt ad},\dots,\ov{T}_1^{\tt ad}\,\rrceil^{\rho_{\mu_1}(\mu,\ell)}$$ 
is admissible of shape $\rho_{n+\mu_1}(\la,\ell)/\rho_{\mu_1}(\mu,\ell)$.
\end{lem}

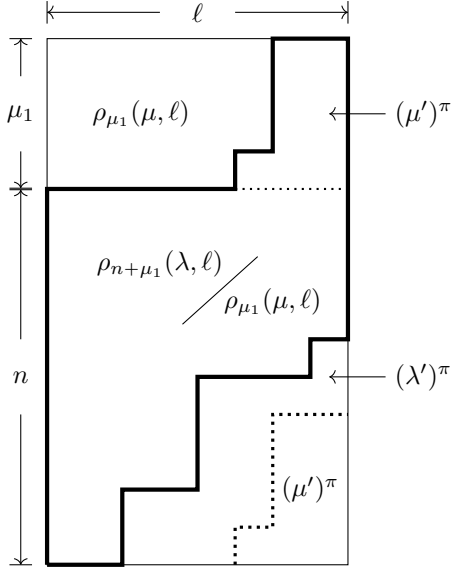
\[
\begin{tikzpicture}
\draw (0, 0) rectangle (4, 7);
\draw [ultra thick] (0, 0) -- (1, 0) -- (1, 1) -- (2, 1) -- (2, 2.5) -- (3.5, 2.5) -- (3.5, 3) -- (4, 3) -- (4, 7) -- (3, 7) -- (3, 5.5) -- (2.5, 5.5) -- (2.5, 5) -- (0, 5) -- (0, 0);
\draw [very thick, dotted] (2.5, 0) -- (2.5, 0.5) -- (3, 0.5) -- (3, 2) -- (4, 2);
\draw [thick, dotted] (2.5, 5) -- (4, 5);
\node at (3.5, 1) {$(\mu')^\pi$};
\node at (1.5, 4) {$\rho_{n+\mu_1}(\lambda, \ell)$}; \node at (3, 3.5) {$\rho_{\mu_1}(\mu, \ell)$};
\draw (1.8, 3.2) -- (2.8, 4.1);
\node [right] at (4.5, 2.5) {$(\lambda')^\pi$}; \draw [->] (4.5, 2.5) -- (3.75, 2.5);
\node [right] at (4.5, 6) {$(\mu')^\pi$}; \draw [->] (4.5, 6) -- (3.75, 6);
\node at (1.25, 6) {$\rho_{\mu_1}(\mu, \ell)$};

\draw (-0.5, 7) -- (-0.2, 7); \draw (-0.5, 4.99) -- (-0.2, 4.99); \draw (-0.5, 5.01) -- (-0.2, 5.01); \draw (-0.5, 0) -- (-0.2, 0);
\draw [->] (-0.35, 3) -- (-0.35, 4.99); \draw [->] (-0.35, 2) -- (-0.35, 0); \node at (-0.35, 2.5) {$n$};
\draw [->] (-0.35, 6.5) -- (-0.35, 7); \draw [->] (-0.35, 5.5) -- (-0.35, 5.01); \node at (-0.35, 6) {$\mu_1$};
\draw (0, 7.2) -- (0, 7.5); \draw (4, 7.2) -- (4, 7.5);
\draw [->] (1.5, 7.35) -- (0, 7.35); \draw [->] (2.5, 7.35) -- (4, 7.35); \node at (2, 7.35) {$\ell$};
\end{tikzpicture}
\]

\pf Suppose that $\ov{\bf T}$ is $[\ov{n}]$-admissible of shape $\la/\mu$. 
Let $m$ be a sufficiently large positive integer.
Let $X=\{\,u_m<\cdots<u_1<1<\cdots<n \,\}$ be a linearly ordered set  of degree $0$, where $u_1,\dots,u_m$ are formal symbols, 
and let $\ov{X}=\{\,\ov{n}<\cdots<\ov{1}< \ov{u}_1<\cdots<\ov{u}_m\,\}$ of degree $0$ similarly.
Since $\ov{\bf T}$ is $[\ov{n}]$-admissible of shape $\la/\mu$, there exists 
${\bf S}=(S_\ell,\dots,S_1)\in {\bf T}_{\ov{X}}(\la,\ell)$, 
where $\ov{\bf T}$ can be obtained from ${\bf S}$ by ignoring the letters in $\{\,\ov{u}_1,\dots,\ov{u}_m\,\}$.

Let $\mc{I}_{m,n}=\{\,u_m<\cdots<u_1<1<\cdots<n<\ov{n}<\cdots<\ov{1}< \ov{u}_1<\cdots<\ov{u}_m\,\}$. 
We may define ${\bf S}^{\tt ad}=(S^{\tt ad}_\ell,\dots,S^{\tt ad}_1)$ as in Proposition \ref{prop:bij spin to KN}, which is a KN tableaux of shape $\rho_{m+n}(\la, \ell)$
with letters in $\mc{I}_{m,n}$.
Then we can check that 
$\llceil T_\ell^{\tt ad},\dots,T_1^{\tt ad}\rrceil^{\rho_{\mu_1}(\mu,\ell)}$
can be obtained from 
$\llceil S_\ell^{\tt ad},\dots,S_1^{\tt ad}\rrceil$
by ignoring the letters in $\{\,u_m, \ldots, u_1\,\}$.
This implies that $\llceil T_\ell^{\tt ad},\dots,T_1^{\tt ad}\rrceil^{\rho_{\mu_1}(\mu,\ell)}$ is $\mc{I}_n$-admissible of shape $\rho_{n+\mu_1}(\la,\ell)/\rho_{\mu_1}(\mu,\ell)$.
The proof of the converse is similar.
\qed
\vskip 2mm

The following is an analogue of Proposition \ref{prop:bij spin to KN} for skew shapes.
\begin{cor}\label{cor:bij spin to KN of skew shape}
We have a bijection    
\begin{equation}\label{eq:skewPsi_lambda:C}
\xymatrixcolsep{2pc}\xymatrixrowsep{0.5pc}\xymatrix{
 {\bf T}_n(\la/\mu,\ell)  \ar@{->}[r]  & \ {\bf KN}_{\rho_{n+\mu_1}(\la,\ell)/\rho_{\mu_1}(\mu,\ell)} \\
  {\bf T}=\llfloor T_\ell,\dots,T_1 \rrfloor_{(\mu_\ell, \dots, \mu_1)}  \ar@{|->}[r] & {\bf T}^{\tt ad}:=\llceil T_\ell^{\tt ad},\dots, T_1^{\tt ad} \rrceil^{\rho_{\mu_1}(\mu,\ell)}}.
\end{equation}
\end{cor}

\begin{ex}{\rm
For $\mathbf{T} = \llfloor T_3, T_2, T_1 \rrfloor_{(0, 1, 2)} \in \mathbf{T}_5(\lambda/\mu, 3)$ with $\lambda = (4, 2, 1)$ and $\mu = (2, 1, 0)$ given below, we have $\mathbf{T}^\texttt{ad} \in \mathbf{KN}_{(6, 5, 3)'/(2, 1, 0)'}$ as follows.
\[ \begin{tikzpicture}
\node [above right] at (0, 0) {$\begin{ytableau}
\bar{\tl 4} & \bar{\tl 2} \\
\bar{\tl 3} & \bar{\tl 1} \\
\bar{\tl 1}
\end{ytableau}$};
\node [above right] at (1.5, 1.2em) {$\begin{ytableau}
\bar{\tl 5} & \bar{\tl 5} \\
\bar{\tl 4} & \bar{\tl 4} \\
\bar{\tl 2} & \bar{\tl 2} \\
\bar{\tl 1}
\end{ytableau}$};
\node [above right] at (3, 2.4em) {$\begin{ytableau}
\bar{\tl 5} & \bar{\tl 5} \\
\bar{\tl 4} & \bar{\tl 2} \\
\bar{\tl 3} \\
\bar{\tl 1}
\end{ytableau}$};

\node [above right] at (7, 0) {$\begin{ytableau}
\none & \none & \tl 3 \\
\none & \tl 3 & \bar{\tl 5} \\
\tl 3 & \bar{\tl 5} & \bar{\tl 3} \\
\tl 5 & \bar{\tl 4} \\
\bar{\tl 3} & \bar{\tl 2} \\
\bar{\tl 1}
\end{ytableau}$};

\draw [dotted] (-0.2, 0.362em) -- (4.7, 0.362em);
\draw [dotted] (6.8, 7.8em) -- (9, 7.8em);
\node [below] at (2.5, -0.2) {$\mathbf{T} =\llfloor T_3, T_2, T_1 \rrfloor_{(0, 1, 2)}$};
\node [below] at (8.2, -0.2) {$\mathbf{T}^\texttt{ad} = \llceil T_3^\texttt{ad}, T_2^\texttt{ad}, T_1^\texttt{ad} \rrceil^{(2, 1, 0)}$};
\end{tikzpicture} \]
} \end{ex}


\subsection{Jeu de taquin for spinor model of a skew shape} 
Now, let us introduce an analogue of jeu de taquin for ${\bf T}\in {\bf T}_\A(\la/\mu,\ell)$. 
\subsubsection{}
We first consider the case when $\ell=2$.
Suppose that ${\bf T}=(T_2,T_1)\in {\bf T}_\A(a_2)\times {\bf T}_\A(a_1)$ is given for some $a_1, a_2\in \Z_+$. Let 
\begin{equation*}
d(T_1,T_2)=\min\left\{\,d\,\Big|\,d\in\Z_+,\ \text{$\llfloor T_2,T_1 \rrfloor_{(0,d)}$ is $\A$-admissible (of a skew shape)}\,\right\}.
\end{equation*}
Note that we have $T_2\prec T_1$ if and only if $d(T_1,T_2)=0$.
Let us assume that ${\bf T}=(T_2^{\tt L},T_2^{\tt R},T_1^{\tt L},T_1^{\tt R})\in {\bf E}^4_\A$. Recall that
\begin{equation}\label{eq:comb R}
\mc{E}_3^{a_2}{\bf T}=({}^{\tt L}T_2,{}^{\tt R}T_2,T_1^{\tt L},T_1^{\tt R}),\quad
\mc{E}_1^{a_1}{\bf T}=(T_2^{\tt L},T_2^{\tt R},{}^{\tt L}T_1,{}^{\tt R}T_1).
\end{equation}

Suppose that $d=d(T_1,T_2)>0$. Let ${\bf T}'$ be given by applying a sequence of jeu de taquin's as follows: 

\begin{itemize}
\item[{\em Case 1}.] Suppose that $\llfloor {}^{\tt R}T_2, T_1^{\tt L} \rrfloor_{(0,d-1)}$ is not $\A$-semistandard. Then we put
\begin{equation}\label{eq:sliding-1}
{\bf T}'=(U_4,U_3,U_2,U_1)=
\begin{cases}
\mc{F}^{a_2-1}_3\mc{E}_2\mc{E}_3^{a_2}{\bf T} & \text{if $\varepsilon_3(\mc{E}_2\mc{E}_3^{a_2}{\bf T})=0$},\\
\mc{F}^{a_2}_3\mc{E}_2\mc{E}_3^{a_2}{\bf T} & \text{if $\varepsilon_3(\mc{E}_2\mc{E}_3^{a_2}{\bf T})=1$}.\\
\end{cases}
\end{equation}

\item[{\em Case 2}.] Suppose that $\llfloor {}^{\tt R}T_2, T_1^{\tt L} \rrfloor_{(0,d-1)}$ is $\A$-semistandard, but $\llfloor T_{2}^{\tt R},{}^{\tt L}T_1 \rrfloor_{(0,d-1)}$ is not.
Then we put
\begin{equation}\label{eq:sliding-2}
{\bf T}'=(U_4,U_3,U_2,U_1)=
\mc{F}^{a_1+1}_1\mc{F}_2\mc{E}_1^{a_1}{\bf T}.
\end{equation}
\end{itemize}

\begin{lem}
Under the above hypothesis, ${\bf T}'$ is well-defined.
\end{lem}
\pf First, consider {\em Case 1}. 
Since $\llfloor {}^{\tt R}T_2, T_1^{\tt L} \rrfloor_{(0,d-1)}$ is not $\A$-semistandard,
we have $\varepsilon_2(\mc{E}_3^{a_2}{\bf T})>0$ and hence $\mc{E}_2\mc{E}_3^{a_2}{\bf T}\neq {\bf 0}$. Let $$(V_4,V_3,V_2,V_1)=\mc{E}_2\mc{E}_3^{a_2}{\bf T}.$$
Since $\mathcal{E}_2 \llfloor {^{\tt L}T_2}, {^{\tt R}T_2}, {T_1^{\tt L}}, {T_1^{\tt R}} \rrfloor_{(0, 0, d, d+a_1)} = \llfloor V_4, V_3, V_2, V_1 \rrfloor_{(0, 1, d-1, d+a_1)}$, it is not difficult to see that $\varepsilon_3(\mc{E}_2\mc{E}_3^{a_2}{\bf T})=0, 1$. Hence ${\bf T}'$ is well-defined.
Similarly, we can check the well-definedness of ${\bf T}'$ in {\em Case 2}.
\qed

\begin{ex}\label{ex:sliding}
{\rm Suppose that $\mathcal{A} = \mathbb{I}_{4|3}$.

(1) The following is an example of {\em Case 1} with $\varepsilon_3(\mathcal{E}_2\mathcal{E}_3^{a_2}\mathbf{T}) = 0$.
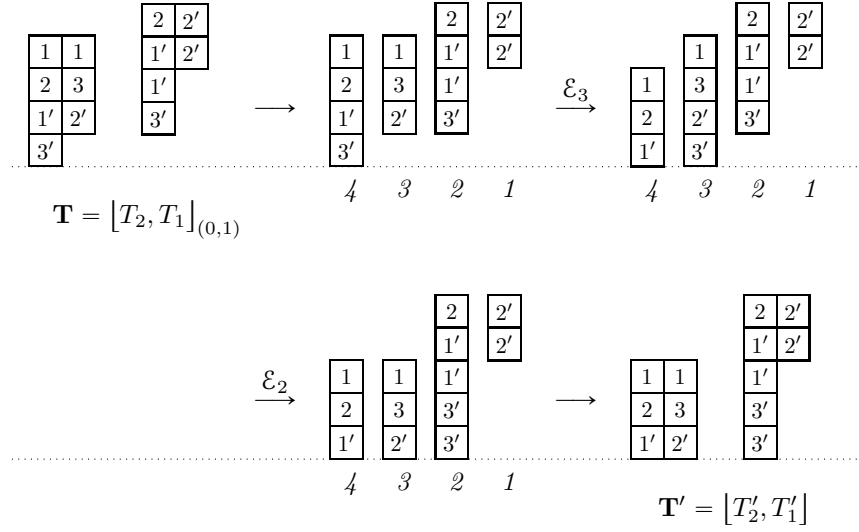
\[ \begin{tikzpicture}
\node [above right] at (0, 0) {$\begin{ytableau}
\tl 1 & \tl 1 \\
\tl 2 & \tl 3 \\
\tl{$1'$} & \tl{$2'$} \\
\tl{$3'$}
\end{ytableau}$};
\node [above right] at (1.5, 1.2em) {$\begin{ytableau}
\tl 2 & \tl{$2'$} \\
\tl{$1'$} & \tl{$2'$} \\
\tl{$1'$} \\
\tl{$3'$}
\end{ytableau}$};

\node [above right] at (4, 0) {$\begin{ytableau}
\tl 1 \\
\tl 2 \\
\tl{$1'$} \\
\tl{$3'$}
\end{ytableau}$};
\node [above right] at (4.7, 0) {$\begin{ytableau}
\tl 1 \\
\tl 3 \\
\tl{$2'$} \\
\none
\end{ytableau}$};
\node [above right] at (5.4, 0) {$\begin{ytableau}
\tl 2 \\
\tl{$1'$} \\
\tl{$1'$} \\
\tl{$3'$} \\
\none
\end{ytableau}$};
\node [above right] at (6.1, 0) {$\begin{ytableau}
\tl{$2'$} \\
\tl{$2'$} \\
\none \\
\none \\
\none
\end{ytableau}$};

\node [above right] at (8, 0) {$\begin{ytableau}
\tl 1 \\
\tl 2 \\
\tl{$1'$}
\end{ytableau}$};
\node [above right] at (8.7, 0) {$\begin{ytableau}
\tl 1 \\
\tl 3 \\
\tl{$2'$} \\
\tl{$3'$}
\end{ytableau}$};
\node [above right] at (9.4, 0) {$\begin{ytableau}
\tl 2 \\
\tl{$1'$} \\
\tl{$1'$} \\
\tl{$3'$} \\
\none
\end{ytableau}$};
\node [above right] at (10.1, 0) {$\begin{ytableau}
\tl{$2'$} \\
\tl{$2'$} \\
\none \\
\none \\
\none
\end{ytableau}$};

\draw [dotted] (-0.1, 0.362em) -- (11.2, 0.362em);
\node at (3.4, 2.5em) {$\longrightarrow$}; \node at (7.4, 2.5em) {$\longrightarrow$};
\node [above] at (7.4, 2.5em) {$\mathcal{E}_3$};
\node [below] at (1.7, -0.2) {$\mathbf{T} = \llfloor T_2, T_1 \rrfloor_{(0, 1)}$};
\node [below] at (4.4, 0.2em) {\em 4}; \node [below] at (5.1, 0.2em) {\em 3}; \node [below] at (5.8, 0.2em) {\em 2}; \node [below] at (6.5, 0.2em) {\em 1};
\node [below] at (8.4, 0.2em) {\em 4}; \node [below] at (9.1, 0.2em) {\em 3}; \node [below] at (9.8, 0.2em) {\em 2}; \node [below] at (10.5, 0.2em) {\em 1};
\end{tikzpicture} \]
\[ \begin{tikzpicture}
\node [above right] at (4, 0) {$\begin{ytableau}
\tl 1 \\
\tl 2 \\
\tl{$1'$}
\end{ytableau}$};
\node [above right] at (4.7, 0) {$\begin{ytableau}
\tl 1 \\
\tl 3 \\
\tl{$2'$}
\end{ytableau}$};
\node [above right] at (5.4, 0) {$\begin{ytableau}
\tl 2 \\
\tl{$1'$} \\
\tl{$1'$} \\
\tl{$3'$} \\
\tl{$3'$}
\end{ytableau}$};
\node [above right] at (6.1, 0) {$\begin{ytableau}
\tl{$2'$} \\
\tl{$2'$} \\
\none \\
\none \\
\none
\end{ytableau}$};

\node [above right] at (8, 0) {$\begin{ytableau}
\tl 1 & \tl 1 \\
\tl 2 & \tl 3 \\
\tl{$1'$} & \tl{$2'$}
\end{ytableau}$};
\node [above right] at (9.5, 0) {$\begin{ytableau}
\tl 2 & \tl{$2'$} \\
\tl{$1'$} & \tl{$2'$} \\
\tl{$1'$} \\
\tl{$3'$} \\
\tl{$3'$}
\end{ytableau}$};

\draw [dotted] (-0.1, 0.362em) -- (11.2, 0.362em);
\node at (3.4, 2.5em) {$\longrightarrow$}; \node at (7.4, 2.5em) {$\longrightarrow$};
\node [above] at (3.4, 2.5em) {$\mathcal{E}_2$};
\node [below] at (9.5, -0.2) {$\mathbf{T}' = \llfloor T_2', T_1' \rrfloor$};
\node [below] at (4.4, 0.2em) {\em 4}; \node [below] at (5.1, 0.2em) {\em 3}; \node [below] at (5.8, 0.2em) {\em 2}; \node [below] at (6.5, 0.2em) {\em 1};
\end{tikzpicture} \]

(2) The following is an example of {\em Case 2}.
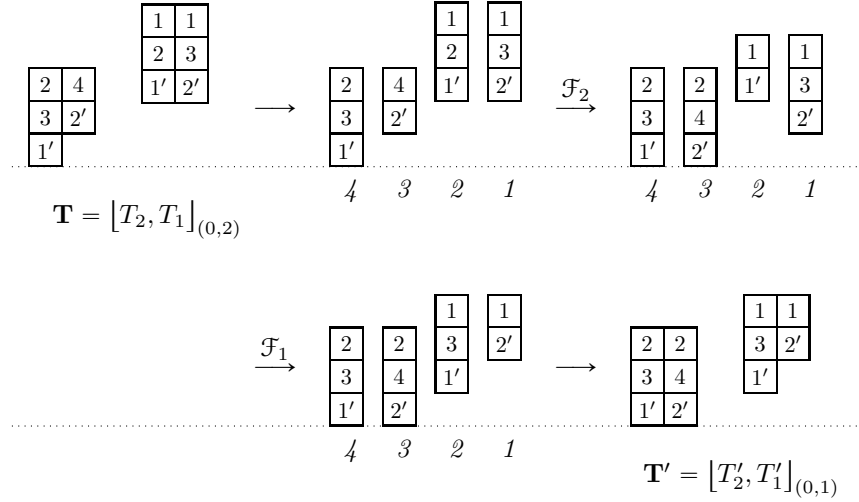
\[ \begin{tikzpicture}
\node [above right] at (0, 0) {$\begin{ytableau}
\tl 2 & \tl 4 \\
\tl 3 & \tl{$2'$} \\
\tl{$1'$}
\end{ytableau}$};
\node [above right] at (1.5, 2.4em) {$\begin{ytableau}
\tl 1 & \tl 1 \\
\tl 2 & \tl 3 \\
\tl{$1'$} & \tl{$2'$}
\end{ytableau}$};

\node [above right] at (4, 0) {$\begin{ytableau}
\tl 2 \\
\tl 3 \\
\tl{$1'$}
\end{ytableau}$};
\node [above right] at (4.7, 0) {$\begin{ytableau}
\tl 4 \\
\tl{$2'$} \\
\none
\end{ytableau}$};
\node [above right] at (5.4, 0) {$\begin{ytableau}
\tl 1 \\
\tl 2 \\
\tl{$1'$} \\
\none \\
\none
\end{ytableau}$};
\node [above right] at (6.1, 0) {$\begin{ytableau}
\tl 1 \\
\tl 3 \\
\tl{$2'$} \\
\none \\
\none
\end{ytableau}$};

\node [above right] at (8, 0) {$\begin{ytableau}
\tl 2 \\
\tl 3 \\
\tl{$1'$}
\end{ytableau}$};
\node [above right] at (8.7, 0) {$\begin{ytableau}
\tl 2 \\
\tl 4 \\
\tl{$2'$}
\end{ytableau}$};
\node [above right] at (9.4, 0) {$\begin{ytableau}
\tl 1 \\
\tl{$1'$} \\
\none \\
\none
\end{ytableau}$};
\node [above right] at (10.1, 0) {$\begin{ytableau}
\tl 1 \\
\tl 3 \\
\tl{$2'$} \\
\none
\end{ytableau}$};

\draw [dotted] (-0.1, 0.362em) -- (11.2, 0.362em);
\node at (3.4, 2.5em) {$\longrightarrow$}; \node at (7.4, 2.5em) {$\longrightarrow$};
\node [above] at (7.4, 2.5em) {$\mathcal{F}_2$};
\node [below] at (1.7, -0.2) {$\mathbf{T} = \llfloor T_2, T_1 \rrfloor_{(0, 2)}$};
\node [below] at (4.4, 0.2em) {\em 4}; \node [below] at (5.1, 0.2em) {\em 3}; \node [below] at (5.8, 0.2em) {\em 2}; \node [below] at (6.5, 0.2em) {\em 1};
\node [below] at (8.4, 0.2em) {\em 4}; \node [below] at (9.1, 0.2em) {\em 3}; \node [below] at (9.8, 0.2em) {\em 2}; \node [below] at (10.5, 0.2em) {\em 1};
\end{tikzpicture} \]
\[ \begin{tikzpicture}
\node [above right] at (4, 0) {$\begin{ytableau}
\tl 2 \\
\tl 3 \\
\tl{$1'$}
\end{ytableau}$};
\node [above right] at (4.7, 0) {$\begin{ytableau}
\tl 2 \\
\tl 4 \\
\tl{$2'$}
\end{ytableau}$};
\node [above right] at (5.4, 0) {$\begin{ytableau}
\tl 1 \\
\tl 3 \\
\tl{$1'$} \\
\none
\end{ytableau}$};
\node [above right] at (6.1, 0) {$\begin{ytableau}
\tl 1 \\
\tl{$2'$} \\
\none \\
\none
\end{ytableau}$};

\node [above right] at (8, 0) {$\begin{ytableau}
\tl 2 & \tl 2 \\
\tl 3 & \tl 4 \\
\tl{$1'$} & \tl{$2'$}
\end{ytableau}$};
\node [above right] at (9.5, 0) {$\begin{ytableau}
\tl 1 & \tl 1 \\
\tl 3 & \tl{$2'$} \\
\tl{$1'$} \\
\none
\end{ytableau}$};

\draw [dotted] (-0.1, 0.362em) -- (11.2, 0.362em);
\node at (3.4, 2.5em) {$\longrightarrow$}; \node at (7.4, 2.5em) {$\longrightarrow$};
\node [above] at (3.4, 2.5em) {$\mathcal{F}_1$};
\node [below] at (9.6, -0.2) {$\mathbf{T}' = \llfloor T_2', T_1' \rrfloor_{(0, 1)}$};
\node [below] at (4.4, 0.2em) {\em 4}; \node [below] at (5.1, 0.2em) {\em 3}; \node [below] at (5.8, 0.2em) {\em 2}; \node [below] at (6.5, 0.2em) {\em 1};
\end{tikzpicture} \]

(3) The following is an example of {\em Case 1} with $\varepsilon_3(\mathcal{E}_2\mathcal{E}_3^{a_2}\mathbf{T}) = 1$.
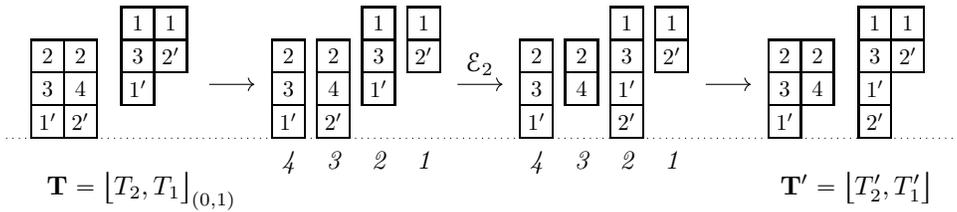
\[ \begin{tikzpicture}
\node [above right] at (0, 0) {$\begin{ytableau}
\tl 2 & \tl 2 \\
\tl 3 & \tl 4 \\
\tl{$1'$} & \tl{$2'$}
\end{ytableau}$};
\node [above right] at (1.2, 0) {$\begin{ytableau}
\tl 1 & \tl 1 \\
\tl 3 & \tl{$2'$} \\
\tl{$1'$} \\
\none
\end{ytableau}$};

\node [above right] at (3.2, 0) {$\begin{ytableau}
\tl 2 \\
\tl 3 \\
\tl{$1'$}
\end{ytableau}$};
\node [above right] at (3.8, 0) {$\begin{ytableau}
\tl 2 \\
\tl 4 \\
\tl{$2'$}
\end{ytableau}$};
\node [above right] at (4.4, 0) {$\begin{ytableau}
\tl 1 \\
\tl 3 \\
\tl{$1'$} \\
\none
\end{ytableau}$};
\node [above right] at (5.0, 0) {$\begin{ytableau}
\tl 1 \\
\tl{$2'$} \\
\none \\
\none
\end{ytableau}$};

\node [above right] at (6.5, 0) {$\begin{ytableau}
\tl 2 \\
\tl 3 \\
\tl{$1'$}
\end{ytableau}$};
\node [above right] at (7.1, 0) {$\begin{ytableau}
\tl 2 \\
\tl 4 \\
\none
\end{ytableau}$};
\node [above right] at (7.7, 0) {$\begin{ytableau}
\tl 1 \\
\tl 3 \\
\tl{$1'$} \\
\tl{$2'$}
\end{ytableau}$};
\node [above right] at (8.3, 0) {$\begin{ytableau}
\tl 1 \\
\tl{$2'$} \\
\none \\
\none
\end{ytableau}$};

\node [above right] at (9.8, 0) {$\begin{ytableau}
\tl 2 & \tl 2 \\
\tl 3 & \tl 4 \\
\tl{$1'$}
\end{ytableau}$};
\node [above right] at (11.0, 0) {$\begin{ytableau}
\tl 1 & \tl 1 \\
\tl 3 & \tl{$2'$} \\
\tl{$1'$} \\
\tl{$2'$}
\end{ytableau}$};

\draw [dotted] (-0.2, 0.362em) -- (12.5, 0.362em); 
\draw [->] (2.5, 2.4em) -- (3.1, 2.4em);
\draw [->] (5.8, 2.4em) -- (6.4, 2.4em);
\draw [->] (9.1, 2.4em) -- (9.7, 2.4em);
\node [above] at (6.1, 2.4em) {$\mathcal{E}_2$};
\node [below] at (3.55, 0.2em) {\em 4}; \node [below] at (4.15, 0.2em) {\em 3}; \node [below] at (4.75, 0.2em) {\em 2}; \node [below] at (5.35, 0.2em) {\em 1};
\node [below] at (6.85, 0.2em) {\em 4}; \node [below] at (7.45, 0.2em) {\em 3}; \node [below] at (8.05, 0.2em) {\em 2}; \node [below] at (8.65, 0.2em) {\em 1};
\node [below] at (1.6, -0.2) {$\mathbf{T} =\llfloor T_2, T_1 \rrfloor_{(0, 1)}$};
\node [below] at (11.1, -0.2) {$\mathbf{T}' =\llfloor T_2', T_1' \rrfloor$};
\end{tikzpicture} \]
} \end{ex}

\begin{prop}\label{prop:jdt for spin col 2 well-defined}
Under the above hypothesis, there exists a unique pair $\mathbf{T}' = (T'_2,T'_1)$ such that 
\begin{gather*}
{\bf T}'=(T_2',T_1')\in {\bf T}_\A(a_2+2\varepsilon-1)\times  {\bf T}_\A(a_1+1),\\
d(T'_1,T'_2) \leq d(T_1,T_2)-1,
\end{gather*}
where $\varepsilon=\varepsilon_3(\mc{E}_2\mc{E}_3^{a_2}{\bf T})$ in Case 1, and $\varepsilon=0$ otherwise.
\end{prop}
\pf Suppose first that $\A=[\ov{n}]$.
Let $d=d(T_1,T_2)$.
Since $\llfloor T_2,T_1 \rrfloor_{(0,d)}$ is $[\ov{n}]$-admissible, we have 
$\llceil T^{\tt ad}_2, T^{\tt ad}_1 \rrceil^{(d,0)}$ 
is admissible by Lemma \ref{lem:admin between skew tableaux-2}. 

Let $T=\llceil C_2,C_1 \rrceil^{(c_2,c_1)}=\llceil T^{\tt ad}_2, T^{\tt ad}_1 \rrceil^{(d,0)}$.
Let $c$ be the inner corner of $T$.
We claim that the algorithm  in Section \ref{subsec:sjdt for KN} to have $T'={\tt jdt}_{KN}(T,c)$ corresponds to either \eqref{eq:sliding-1} or \eqref{eq:sliding-2}.

First, we put 
\raisebox{-0.45em}{$
\ytableausetup{boxsize=1.3em} 
\begin{ytableau} \bullet
\end{ytableau}
$}  
at the top of $C_2$ and apply {\em Step 1.}(1) as far as possible to have \eqref{eq:punctured}. 
Since $d=d(T_1,T_2)>0$, we should apply {\em Step 1.}(2) to \eqref{eq:punctured} with $a'>b$. 

Suppose that $b\in [\ov{n}]$. Then it is straightforward to see from Lemma \ref{lem:lCrC} and \eqref{eq:comb R} that 
applying {\em Step 1}.(2-b) and sliding 
\raisebox{-0.45em}{$
\ytableausetup{boxsize=1.3em} 
\begin{ytableau} \bullet
\end{ytableau}
$}
to the bottom of the column corresponds to \eqref{eq:sliding-2}. 
This implies that $(U_4,U_3)=((T'_2)^{\tt L}, (T'_2)^{\tt R})$ for some $T'_2\in {\bf T}_n(a_2-1)$ and $(U_2,U_1)=((T'_1)^{\tt L}, (T'_1)^{\tt R})$ for some $T'_1\in {\bf T}_n(a_1+1)$. 
Furthermore, since $T'=\llceil C_2',C'_1 \rrceil^{(c_2-1,c_1)}$, 
it follows from Lemma \ref{lem:admin between skew tableaux-2} that
$\llfloor T'_2,T'_1 \rrfloor_{(0,d-1)}$ is $\A$-admissible and
\begin{equation}\label{eq:aux}
({\bf T}')^{\tt ad}=T'=\llceil C_2',C'_1 \rrceil^{(c_2-1,c_1)},
\end{equation}
where $({\bf T}')^{\tt ad}$ is given in Corollary \ref{cor:bij spin to KN of skew shape}.
This implies that $d(T'_1,T'_2) \leq d-1$.

Next, suppose that $b\in [{n}]$. It is easy to see that 
exchanging 
\raisebox{-0.45em}{$
\ytableausetup{boxsize=1.3em} 
\begin{ytableau} \bullet
\end{ytableau}
$}  
with 
\raisebox{-0.45em}{$
\ytableausetup{boxsize=1.3em} 
\begin{ytableau} b
\end{ytableau}
$} in {\em Step 1.}(2-a) to have $C_1'$
corresponds to $\mc{E}_2\mc{E}_3^{a_2}{\bf T}$. 
If $C_2'$ is admissible, which is equivalent to $\varepsilon_3(\mc{E}_2\mc{E}_3^{a_2}{\bf T})=0$, 
then the process to have $C_2'$ corresponds to applying $\mc{F}_3^{a_2-1}$ to $\mc{E}_2\mc{E}_3^{a_2}{\bf T}$. 
As in the previous case, we have $(U_4, U_3) = ((T_2')^{\tt L}, (T_2')^{\tt R})$ and $(U_2, U_1) = ((T_1')^{\tt L}, (T_1')^{\tt R})$ 
for some $T_2' \in \mathbf{T}_n(a_2-1)$ and $T_1' \in \mathbf{T}_n(a_1+1)$,
and $\llfloor T'_2,T'_1 \rrfloor_{(0,d-1)}$ is $\A$-admissible with \eqref{eq:aux}, which implies $d(T'_1,T'_2) \leq d(T_1,T_2)-1$.
On the other hand, if $C_2'$ is not admissible, which is equivalent to $\varepsilon_3(\mc{E}_2\mc{E}_3^{a_2}{\bf T})=1$, 
then it is not difficult to see that the process to have $C_2''$ corresponds to applying $\mc{F}_3^{a_2}$ to $\mc{E}_2\mc{E}_3^{a_2}{\bf T}$.
Hence, we have $(U_4,U_3)=((T'_2)^{\tt L}, (T'_2)^{\tt R})$ for some $T'_2\in {\bf T}_n(a_2+1)$ and $(U_2,U_1)=((T'_1)^{\tt L}, (T'_1)^{\tt R})$ for some $T'_1\in {\bf T}_n(a_1+1)$. 
In this case, we have $T'=\llceil C''_2, C'_1 \rrceil^{(c_2,c_1)}$, and hence $\llfloor T'_2,T'_1 \rrfloor_{(0,d)}$ is $\A$-admissible.

Moreover, it can be shown that $\llceil C_2'', C_1' \rrceil^{(c_2-1, c_1)}$ is also admissible, hence $\llfloor T_2', T_1' \rrfloor_{(0, d-1)}$ is $\mathcal{A}$-admissible with \eqref{eq:aux}, which implies $d(T_1', T'_2) \le d(T_1, T_2)-1$.

Now, suppose that $\A$ is arbitrary. Let $\ov{\bf T}$ be the $n$-conjugate of ${\bf T}$ for a sufficiently large $n$. Let $\mc{X}$ be the composite of operators $\mc{E}_i$ and $\mc{F}_i$ in \eqref{eq:sliding-1} or \eqref{eq:sliding-2}. Then we have by definition of $\ov{\bf T}$ and Lemma \ref{lem:sl crystal}
\begin{equation}\label{eq:recording is the same}
Q(\mc{X}{\bf T})=\mc{X}Q({\bf T})=\mc{X}Q(\ov{\bf T})=Q(\mc{X}\ov{\bf T}).
\end{equation}
By the previous arguments for the case of $\A=[\ov{n}]$, there exists $\ov{\bf T}'=(\ov{T}'_2,\ov{T}'_1)$ such that
$$\mc{X}\ov{\bf T}=\ov{\bf T}'\in {\bf T}_n(a_2+2\varepsilon-1)\times{\bf T}_n(a_1+1).$$
By Lemma \cite[Lemma 6.2]{K15} and \eqref{eq:recording is the same}, we have ${\bf T}'=({T}'_2,{T}'_1)$ such that
$$\mc{X}{\bf T}={\bf T}'\in {\bf T}_{\A}(a_2+2\varepsilon-1)\times{\bf T}_{\A}(a_1+1),$$
and the $n$-conjugate of $T'_i$ is $\ov{T}'_i$ for $i=1,2$.
Finally, it follows from the arguments for $\A=[\ov{n}]$ and Lemma \ref{lem:admin between skew tableaux-1} that $d(T'_1,T'_2)\leq d(T_1,T_2)-1$.
\qed\vskip 2mm

\begin{df}\label{def:sjdt spin 2 col}
{\rm
For ${\bf T}=(T_2,T_1)\in {\bf T}_\A(a_2)\times {\bf T}_\A(a_1)$ with $d(T_1,T_2)>0$, we define
\begin{equation*}
{\tt jdt}_{spin}({\bf T})={\bf T}',
\end{equation*}
where ${\bf T}'$ is given in \eqref{eq:sliding-1} and \eqref{eq:sliding-2}.
}
\end{df}

The following corollaries follow from the proof of Proposition \ref{prop:jdt for spin col 2 well-defined}.

\begin{cor}\label{cor:sjdt spin=KN}
Under the above hypothesis, if $\A=[\ov{n}]$, then we have
\begin{equation*} 
\left({\tt jdt}_{spin}({\bf T})\right)^{\tt ad}
={\tt jdt}_{KN}\left({\bf T}^{\tt ad},c\right),
\end{equation*}
where $c$ is the inner corner of $\mathbf{T}^{\tt ad}$.
\end{cor}

\begin{cor}\label{cor:compatibility for jdt}
Under the above hypothesis, we have
\begin{align*}
\ov{{\tt jdt}_{spin}({\bf T})} 
&={\tt jdt}_{spin}(\ov{\bf T}),\\
\left({\tt jdt}_{spin}(\ov{\bf T})\right)^{\tt ad}
&={\tt jdt}_{KN}\left(\ov{\bf T}^{\tt ad},c\right),
\end{align*}
where $\ov{\, \cdot \, }$ denotes the $n$-conjugate for a sufficiently large $n$, $(\,\cdot \, )^{\tt ad}$ is given in \eqref{eq:Psi_lambda:C} or \eqref{eq:skewPsi_lambda:C}, and $c$ is the inner corner of $\ov{\bf T}^{\tt ad}$.
\end{cor}


\subsubsection{}
Now we consider a general case.
Let $\la/\mu$ be a skew diagram with $\la,\mu\in\cP_\ell$ and 
let ${\bf T}=(T_\ell,\dots,T_1)\in {\bf T}_\A(\la/\mu,\ell)$ be given.
Let $c$ be an inner corner of $\la/\mu$ in the $i$-th row from the top.

Let us define an analogue of jeu de taquin sliding on ${\bf T}$ with respect to $c$.
We first take a sufficiently large $n$ and the $n$-conjugate $\ov{\bf T}$ of ${\bf T}$.  
Then we consider 
\begin{equation}\label{eq:sjdt at c}
{\tt jdt}_{KN}\left(\ov{\bf T}^{\tt ad}, b \right),
\end{equation} 
where $b$ is the inner corner of $\overline{\bf T}^\texttt{ad}$ in the $(i+1)$-th column from the right and $(\,\cdot \, )^{\tt ad}$ is as in \eqref{eq:skewPsi_lambda:C}. 
Recall that \eqref{eq:sjdt at c} is obtained by applying a sequence of ${\tt jdt}_{KN}$ to two neighboring components in $\ov{\bf T}^{\tt ad}$ (see {\em Step 2} in Section \ref{subsec:sjdt for KN}).
By Definition \ref{def:sjdt spin 2 col} and Corollary \ref{cor:compatibility for jdt}, there exists a composite of operators $\mc{E}_i$ and $\mc{F}_i$, say $\mc{X}$, such that
$\left( \mc{X}\,\ov{\bf T}\right)^{\tt ad} = {\tt jdt}_{KN}\left(\ov{\bf T}^{\tt ad},b\right)$.
\begin{df}{\rm 
Under the above hypothesis, we define
\begin{equation}\label{eq:sjdt spin general}
{\tt jdt}_{spin}({\bf T},c)=\mc{X}\,{\bf T}.
\end{equation}}
\end{df}
Note that ${\tt jdt}_{spin}({\bf T},c)$ is independent of the choice of $\mc{X}$ since
\begin{equation*}
\left(\ov{{\tt jdt}_{spin}({\bf T},c)}\right)^{\tt ad} 
= \left(\ov{\mc{X}\,{\bf T}}\right)^{\tt ad}
= \left(\mc{X}\,\ov{{\bf T}}\right)^{\tt ad}
= {\tt jdt}_{KN}\left(\ov{\bf T}^{\tt ad}, b \right),
\end{equation*}
and $\ov{\, \cdot \, }$ and $(\,\cdot \, )^{\tt ad}$ are injective on the connected component of $\mathbf{T}$.

\begin{thm}\label{thm:P tableau for spinor}
Let $\la/\mu$ be a skew diagram with $\la,\mu\in\cP_\ell$ and 
let ${\bf T}=(T_\ell,\dots,T_1)\in {\bf T}_\A(\la/\mu,\ell)$ be given.
There exists a unique ${\tt P}({\bf T})\in {\bf T}_\A(\nu,\ell)$ for some $(\nu,\ell)\in \mc{P}({\rm Sp})_\A$, which can be obtained from ${\bf T}$ by applying ${\tt jdt}_{spin}(\,\cdot\,,c)$ finitely many times with respect to inner corners. In particular, if $\A=[\ov{n}]$, then we have
\begin{equation*}
{\tt P}({\bf T})^{\tt ad}={\tt P}\left(w\left({\bf T}^{\tt ad}\right)\right).
\end{equation*}

\end{thm}
\pf Let us first prove the existence of ${\tt P}({\bf T})$.
Let $\ov{\bf T}$ be the $n$-conjugate of ${\bf T}$ for a sufficiently large $n$.
Let ${\bf U}=\ov{\bf T}$ and ${\bf V}=\ov{\bf T}^{\tt ad}$.
By Section \ref{subsec:sjdt for KN}, there exists a sequence ${\bf V}={\bf V}_0,\dots,{\bf V}_r$ such that 
\begin{equation}\label{eq:seq for V}
{\bf V}_{i+1}={\tt jdt}_{KN}({\bf V}_i,b_i)\quad (1\leq i\leq r-1),
\end{equation} 
for some inner corner $b_i$ in ${\rm sh}({\bf V}_i)$, and ${\bf V}_r \in {\bf KN}_\delta$ for some $\delta\in \cP_n$ with $\delta_1\leq \ell$.
By Corollary \ref{cor:compatibility for jdt}, there exists a sequence ${\bf U}={\bf U}_0,\dots,{\bf U}_r$ such that ${\bf U}_i^{\tt ad}={\bf V}_i$ and 
\begin{equation}\label{eq:seq for U}
{\bf U}_{i+1}={\tt jdt}_{spin}({\bf U}_i,c_i)\quad (1\leq i\leq r-1)
\end{equation}
for some inner corners $c_i$ in $\mbox{sh}(\mathbf{U}_i)$.
Again by Corollary \ref{cor:compatibility for jdt}, there exists a sequence ${\bf T}={\bf T}_0,\dots,{\bf T}_r$ such that $\ov{{\bf T}}_i={\bf U}_i$ (the $n$-conjugate of ${\bf T}_i$) and 
\begin{equation}\label{eq:seq for T}
{\bf T}_{i+1}={\tt jdt}_{spin}({\bf T}_i,c_i)\quad (1\leq i\leq r-1).
\end{equation}
Since ${\bf V}_r \in {\bf KN}_\delta$, we conclude from Lemmas \ref{lem:admin between skew tableaux-1} and \ref{lem:admin between skew tableaux-2} that ${\bf T}_r\in {\bf T}_\A(\nu,\ell)$ with $\rho_{n}(\nu,\ell)=\delta$. We put ${\tt P}({\bf T})={\bf T}_r$.

Now let us prove the uniqueness. 
Suppose that there exists a sequence 
${\bf T}={\bf T}'_0,\dots,{\bf T}'_s$ such that 
${\bf T}'_{i+1}={\tt jdt}_{spin}({\bf T}'_i,c'_i)$ for some $c'_i \ (1\leq i\leq s-1)$ and ${\bf T}'_s\in {\bf T}_\A(\xi,\ell)$ for some $(\xi,\ell)\in \mc{P}({\rm Sp})_\A$. 
We claim that ${\bf T}_r={\bf T}'_s$.

If we put ${\bf U}'_i=\ov{{\bf T}'_i}$ and ${\bf V}'_i=({\bf U}'_i)^{\tt ad}$ for $1\leq i\leq s$, then they also satisfy \eqref{eq:seq for V} and \eqref{eq:seq for U} by Corollary \ref{cor:compatibility for jdt}. 
We have ${\bf V}'_s\in {\bf KN}_{\rho_{n}(\xi,\ell)}$ by Lemmas \ref{lem:admin between skew tableaux-1} and \ref{lem:admin between skew tableaux-2}. 
By Proposition \ref{prop:plactic for C}, we have
$${\bf V}'_s={\tt P}(w({\bf V}'_0))={\tt P}(w({\bf V}_0))={\bf V}_r ,$$
and hence $\xi=\nu$. 
On the other hand, let $\mc{X}$ and $\mc{X}'$ be composites of $\mc{E}_i$ and $\mc{F}_i$ $(1\leq i\leq 2\ell-1)$ such that 
\begin{equation*}
\mc{X}{\bf T}={\bf T}_r,\quad \mc{X}'{\bf T}={\bf T}'_s, 
\end{equation*}
which implies $\mc{X}{\bf U}={\bf U}_r$ and $\mc{X}'{\bf U}={\bf U}'_s$, respectively.
Since ${\bf V}_r={\bf V}'_s$, we have $\mc{X}{\bf U}={\bf U}_r={\bf U}'_s=\mc{X}'{\bf U}$, and 
\begin{equation*}
\begin{split}
Q(\mc{X}{\bf T})&=\mc{X}Q({\bf T})=\mc{X}Q({\bf U})=Q(\mc{X}{\bf U})\\
&=Q(\mc{X}'{\bf U})=\mc{X}'Q({\bf U})=\mc{X}Q({\bf T})=Q(\mc{X}'{\bf T}).
\end{split}
\end{equation*}
Recall that $Q({\bf T})=Q({\bf U})$ by definition of the $n$-conjugate.
Since $P(\mc{X}{\bf T})=P({\bf T})=P(\mc{X}'{\bf T})$, we have 
$\Phi_\A({\bf T}_r)=\Phi_\A({\bf T}'_s)$ and hence ${\bf T}_r={\bf T}'_s$ by Theorem \ref{thm:Ainsertion}. 

Finally, if $\A=[\ov{n}]$, then it follows directly from Corollary \ref{cor:sjdt spin=KN} that ${\tt P}({\bf T})^{\tt ad}={\tt P}\left(w\left({\bf T}^{\tt ad}\right)\right)$.  This completes the proof.
\qed

\begin{ex} \label{ex:insertion tableau}{\rm
Let ${\bf T}=\llfloor T_3, T_2, T_1 \rrfloor_{(0, 1, 2)}$ be the tableau in Example \ref{ex:conjugate}. 
Then ${\tt P}({\bf T})$ can be obtained as follows (for detailed computation, see Example \ref{ex:sliding})
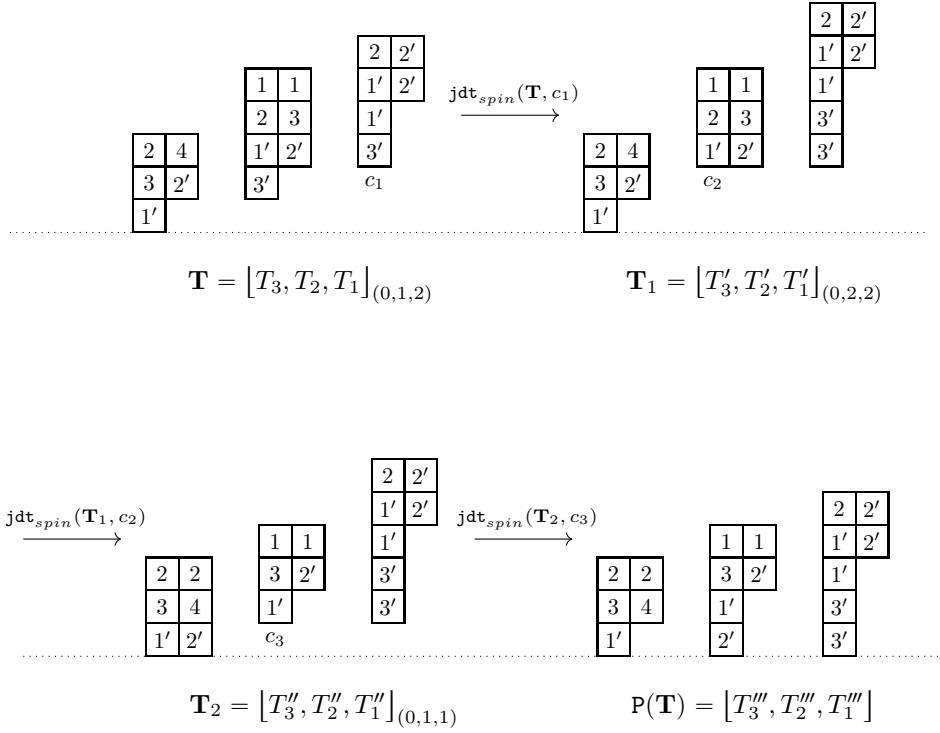
\[ \begin{tikzpicture} \ytableausetup{boxsize=1.2em}
\node [above right] at (0, 0) {$\begin{ytableau}
\tl 2 & \tl 4 \\
\tl 3 & \tl{$2'$} \\
\tl{$1'$}
\end{ytableau}$};
\node [above right] at (1.5, 0) {$\begin{ytableau}
\tl 1 & \tl 1 \\
\tl 2 & \tl 3 \\
\tl{$1'$} & \tl{$2'$} \\
\tl{$3'$} \\
\none
\end{ytableau}$};
\node [above right] at (3, 0) {$\begin{ytableau}
\tl 2 & \tl{$2'$} \\
\tl{$1'$} & \tl{$2'$} \\
\tl{$1'$} \\
\tl{$3'$} \\
\none[\tl{$c_1$}] \\
\none
\end{ytableau}$};

\node [above right] at (6, 0) {$\begin{ytableau}
\tl 2 & \tl 4 \\
\tl 3 & \tl{$2'$} \\
\tl{$1'$}
\end{ytableau}$};
\node [above right] at (7.5, 0) {$\begin{ytableau}
\tl 1 & \tl 1 \\
\tl 2 & \tl 3 \\
\tl{$1'$} & \tl{$2'$} \\
\none[\tl{$c_2$}] \\
\none
\end{ytableau}$};
\node [above right] at (9, 0) {$\begin{ytableau}
\tl 2 & \tl{$2'$} \\
\tl{$1'$} & \tl{$2'$} \\
\tl{$1'$} \\
\tl{$3'$} \\
\tl{$3'$} \\
\none \\
\none
\end{ytableau}$};

\draw [dotted] (-1.5, 0.362em) -- (10.7, 0.362em); 
\draw [->] (4.5, 4.8em) -- (5.8, 4.8em);
\node [above] at (5.2, 4.8em) {\scriptsize $\texttt{jdt}_{spin}(\mathbf{T}, c_1)$};
\node [below] at (2.5, -0.2) {$\mathbf{T} =\llfloor T_3, T_2, T_1 \rrfloor_{(0, 1, 2)}$};
\node [below] at (8.4, -0.2) {$\mathbf{T}_1 = \llfloor T_3', T_2', T_1' \rrfloor_{(0, 2, 2)}$};
\end{tikzpicture} \]

\[ \begin{tikzpicture}
\node [above right] at (0, 0) {$\begin{ytableau}
\tl 2 & \tl 2 \\
\tl 3 & \tl 4 \\
\tl{$1'$} & \tl{$2'$}
\end{ytableau}$};
\node [above right] at (1.5, 0) {$\begin{ytableau}
\tl 1 & \tl 1 \\
\tl 3 & \tl{$2'$} \\
\tl{$1'$} \\
\none[\tl{$c_3$}]
\end{ytableau}$};
\node [above right] at (3, 0) {$\begin{ytableau}
\tl 2 & \tl{$2'$} \\
\tl{$1'$} & \tl{$2'$} \\
\tl{$1'$} \\
\tl{$3'$} \\
\tl{$3'$} \\
\none
\end{ytableau}$};

\node [above right] at (6, 0) {$\begin{ytableau}
\tl 2 & \tl 2 \\
\tl 3 & \tl 4 \\
\tl{$1'$}
\end{ytableau}$};
\node [above right] at (7.5, 0) {$\begin{ytableau}
\tl 1 & \tl 1 \\
\tl 3 & \tl{$2'$} \\
\tl{$1'$} \\
\tl{$2'$}
\end{ytableau}$};
\node [above right] at (9, 0) {$\begin{ytableau}
\tl 2 & \tl{$2'$} \\
\tl{$1'$} & \tl{$2'$} \\
\tl{$1'$} \\
\tl{$3'$} \\
\tl{$3'$}
\end{ytableau}$};

\draw [dotted] (-1.5, 0.362em) -- (10.7, 0.362em); 
\draw [->] (-1.5, 4.8em) -- (-0.2, 4.8em);
\node [above] at (-0.8, 4.8em) {\scriptsize $\texttt{jdt}_{spin}(\mathbf{T}_1, c_2)$};
\draw [->] (4.5, 4.8em) -- (5.8, 4.8em);
\node [above] at (5.2, 4.8em) {\scriptsize $\texttt{jdt}_{spin}(\mathbf{T}_2, c_3)$};
\node [below] at (2.5, -0.2) {$\mathbf{T}_2 = \llfloor T_3'', T_2'', T_1'' \rrfloor_{(0, 1, 1)}$};
\node [below] at (8.2, -0.2) {$\texttt{P}(\mathbf{T})=\llfloor T_3''', T_2''', T_1''' \rrfloor$};
\end{tikzpicture} \]
where
\begin{align*}
\texttt{jdt}_{spin}(\mathbf{T}, c_1) &= \mathcal{E}_2 \mathcal{E}_3 \mathbf{T}, \\
\texttt{jdt}_{spin}(\mathbf{T}_1, c_2) &= \mathcal{F}_3 \mathcal{F}_4 \mathbf{T}_1, \\
\texttt{jdt}_{spin}(\mathbf{T}_2, c_3) &= \mathcal{E}_4 \mathbf{T}_2.
\end{align*}
The corresponding jeu de taquin for the  $5$-conjugate ${\bf U}$ and ${\bf V}={\bf U}^{\tt ad}$ is given as follows.
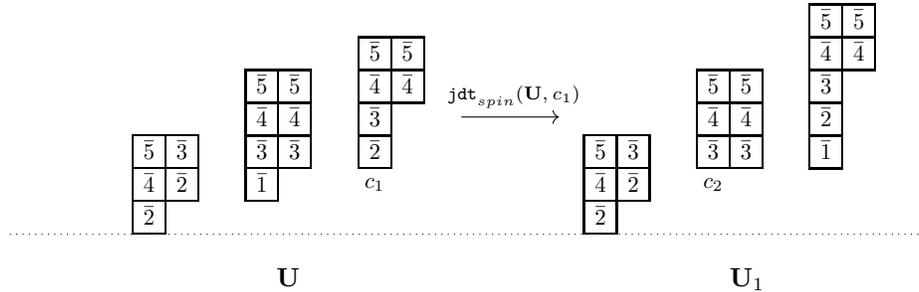
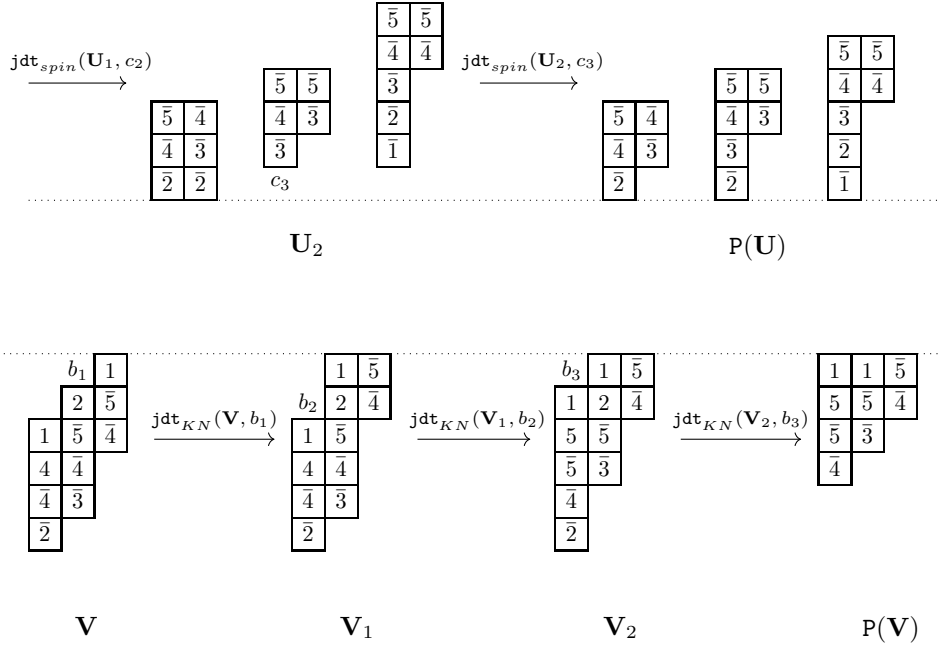
\[ \begin{tikzpicture} \ytableausetup{boxsize=1.2em}
\node [above right] at (0, 0) {$\begin{ytableau}
\bar{\tl 5} & \bar{\tl 3} \\
\bar{\tl 4} & \bar{\tl 2} \\
\bar{\tl 2}
\end{ytableau}$};
\node [above right] at (1.5, 0) {$\begin{ytableau}
\bar{\tl 5} & \bar{\tl 5} \\
\bar{\tl 4} & \bar{\tl 4} \\
\bar{\tl 3} & \bar{\tl 3} \\
\bar{\tl 1} \\
\none
\end{ytableau}$};
\node [above right] at (3, 0) {$\begin{ytableau}
\bar{\tl 5} & \bar{\tl 5} \\
\bar{\tl 4} & \bar{\tl 4} \\
\bar{\tl 3} \\
\bar{\tl 2} \\
\none[\tl{$c_1$}] \\
\none
\end{ytableau}$};

\node [above right] at (6, 0) {$\begin{ytableau}
\bar{\tl 5} & \bar{\tl 3} \\
\bar{\tl 4} & \bar{\tl 2} \\
\bar{\tl 2}
\end{ytableau}$};
\node [above right] at (7.5, 0) {$\begin{ytableau}
\bar{\tl 5} & \bar{\tl 5} \\
\bar{\tl 4} & \bar{\tl 4} \\
\bar{\tl 3} & \bar{\tl 3} \\
\none[\tl{$c_2$}] \\
\none
\end{ytableau}$};
\node [above right] at (9, 0) {$\begin{ytableau}
\bar{\tl 5} & \bar{\tl 5} \\
\bar{\tl 4} & \bar{\tl 4} \\
\bar{\tl 3} \\
\bar{\tl 2} \\
\bar{\tl 1} \\
\none \\
\none
\end{ytableau}$};

\draw [dotted] (-1.5, 0.362em) -- (10.7, 0.362em); 
\draw [->] (4.5, 4.8em) -- (5.8, 4.8em);
\node [above] at (5.2, 4.8em) {\scriptsize $\texttt{jdt}_{spin}(\mathbf{U}, c_1)$};
\node [below] at (2.2, -0.2) {$\mathbf{U}$};
\node [below] at (8.3, -0.2) {$\mathbf{U}_1$};
\end{tikzpicture} \]

\[ \begin{tikzpicture}
\node [above right] at (0, 0) {$\begin{ytableau}
\bar{\tl 5} & \bar{\tl 4} \\
\bar{\tl 4} & \bar{\tl 3} \\
\bar{\tl 2} & \bar{\tl 2}
\end{ytableau}$};
\node [above right] at (1.5, 0) {$\begin{ytableau}
\bar{\tl 5} & \bar{\tl 5} \\
\bar{\tl 4} & \bar{\tl 3} \\
\bar{\tl 3} \\
\none[\tl{$c_3$}]
\end{ytableau}$};
\node [above right] at (3, 0) {$\begin{ytableau}
\bar{\tl 5} & \bar{\tl 5} \\
\bar{\tl 4} & \bar{\tl 4} \\
\bar{\tl 3} \\
\bar{\tl 2} \\
\bar{\tl 1} \\
\none
\end{ytableau}$};

\node [above right] at (6, 0) {$\begin{ytableau}
\bar{\tl 5} & \bar{\tl 4} \\
\bar{\tl 4} & \bar{\tl 3} \\
\bar{\tl 2}
\end{ytableau}$};
\node [above right] at (7.5, 0) {$\begin{ytableau}
\bar{\tl 5} & \bar{\tl 5} \\
\bar{\tl 4} & \bar{\tl 3} \\
\bar{\tl 3} \\
\bar{\tl 2}
\end{ytableau}$};
\node [above right] at (9, 0) {$\begin{ytableau}
\bar{\tl 5} & \bar{\tl 5} \\
\bar{\tl 4} & \bar{\tl 4} \\
\bar{\tl 3} \\
\bar{\tl 2} \\
\bar{\tl 1}
\end{ytableau}$};

\draw [dotted] (-1.5, 0.362em) -- (10.7, 0.362em); 
\draw [->] (-1.5, 4.8em) -- (-0.2, 4.8em);
\node [above] at (-0.8, 4.8em) {\scriptsize $\texttt{jdt}_{spin}(\mathbf{U}_1, c_2)$};
\draw [->] (4.5, 4.8em) -- (5.8, 4.8em);
\node [above] at (5.2, 4.8em) {\scriptsize $\texttt{jdt}_{spin}(\mathbf{U}_2, c_3)$};
\node [below] at (2.2, -0.2) {$\mathbf{U}_2$};
\node [below] at (8.2, -0.2) {$\texttt{P}(\mathbf{U})$};
\end{tikzpicture} \]

\[ \begin{tikzpicture}
\node [below right] at (0, 0) {$\begin{ytableau}
\none & \none[\tl{$b_1$}] & \tl{$1$} \\
\none & \tl{$2$} & \bar{\tl 5} \\
\tl{$1$} & \bar{\tl 5} & \bar{\tl 4} \\
\tl{$4$} & \bar{\tl 4} \\
\bar{\tl 4} & \bar{\tl 3} \\
\bar{\tl 2}
\end{ytableau}$};

\node [below right] at (3.5, 0) {$\begin{ytableau}
\none & \tl{$1$} & \bar{\tl 5} \\
\none[\tl{$b_2$}] & \tl{$2$} & \bar{\tl 4} \\
\tl{$1$} & \bar{\tl 5} \\
\tl{$4$} & \bar{\tl 4} \\
\bar{\tl 4} & \bar{\tl 3} \\
\bar{\tl 2}
\end{ytableau}$};

\node [below right] at (7, 0) {$\begin{ytableau}
\none[\tl{$b_3$}] & \tl{$1$} & \bar{\tl 5} \\
\tl{$1$} & \tl{$2$} & \bar{\tl 4} \\
\tl{$5$} & \bar{\tl 5} \\
\bar{\tl 5} & \bar{\tl 3} \\
\bar{\tl 4} \\
\bar{\tl 2}
\end{ytableau}$};

\node [below right] at (10.5, 0) {$\begin{ytableau}
\tl{$1$} & \tl{$1$} & \bar{\tl 5} \\
\tl{$5$} & \bar{\tl 5} & \bar{\tl 4} \\
\bar{\tl 5} & \bar{\tl 3} \\
\bar{\tl 4}
\end{ytableau}$};

\draw [dotted] (-0.2, -0.362em) -- (12.3, -0.362em); 
\draw [->] (1.8, -3.6em) -- (3.4, -3.6em);
\node [above] at (2.6, -3.6em) {\scriptsize $\texttt{jdt}_{KN}(\mathbf{V}, b_1)$};
\draw [->] (5.3, -3.6em) -- (6.8, -3.6em);
\node [above] at (6.1, -3.6em) {\scriptsize $\texttt{jdt}_{KN}(\mathbf{V}_1, b_2)$};
\draw [->] (8.8, -3.6em) -- (10.4, -3.6em);
\node [above] at (9.6, -3.6em) {\scriptsize $\texttt{jdt}_{KN}(\mathbf{V}_2, b_3)$};
\node [below] at (0.9, -3.5) {$\mathbf{V}$};
\node [below] at (4.5, -3.5) {$\mathbf{V}_1$};
\node [below] at (8.0, -3.5) {$\mathbf{V}_2$};
\node [below] at (11.6, -3.5) {$\texttt{P}(\mathbf{V})$};
\end{tikzpicture} \]
} \end{ex}


\section{Recording tableaux for spinor model} \label{sec:recording for spinor}

\subsection{Oscillating tableaux of shape $(\la,\ell)$} 

Recall that an {\it oscillating tableau} is a sequence of partitions $Q=(Q_1,\dots,Q_s)$ for some $s\geq 1$ such that 
each pair $(Q_i,Q_{i+1})$ differs by one box for $1\leq i\leq s-1$.
We say that an oscillating tableau $Q=(Q_1, \dots, Q_s)$ is {\it vertical} if 
$Q_1\subsetneq  \cdots \subsetneq Q_r \supsetneq \cdots \supsetneq Q_s$ for some $1 \le r \le s$ and $Q_r/Q_1$ and $Q_r/Q_s$ is a skew diagram of vertical strip. We denote by $|Q|=s$ the length of $Q$.

Let $(\la,\ell)\in \cP({\rm Sp})$ be given.
For $n\geq \la_1$, we define ${\bf O}(\la,\ell;n)$ to be the set of a sequence 
of oscillating tableaux $Q=(Q^{(1)}:\cdots :Q^{(\ell)})$ such that
\begin{itemize}
\item[(1)] $Q$ is itself an oscillating tableau,

\item[(2)] $Q^{(i)}=(Q_{i,1},\dots,Q_{i,s_i})$ is a vertical oscillating tableau for $1\leq i\leq \ell$,

\item[(3)] $\ell(Q_{i, j}) \le n$ for $1 \le i \le \ell$ and $1 \le j \le s_i$,

\item[(4)] $Q_{1,1}=\raisebox{0.2em}{\boxed{\mbox{}}}$ and $Q_{\ell,s_\ell}=\rho_n(\la,\ell)$.
\end{itemize} 

Let us consider a stable limit of $Q$. 
More precisely, let $\sigma(Q)=(\widehat{Q}^{(1)}:\cdots:\widehat{Q}^{(\ell)})$, 
where $\widehat{Q}^{(i)}$ is a vertical oscillating tableau with $|\widehat Q^{(i)}|=s_i+1$ given as follows;
\begin{equation*}
\widehat{Q}^{(i)} =((i) \cup Q_{i-1, s_{i-1}},(i)\cup Q_{i,1},\dots,(i)\cup Q_{i,s_i}) \quad (1\leq i\leq \ell).
\end{equation*}
Here we note that $Q_{i,k}$ is partition with $\ell(Q_{i,k}')\leq i$ for $1\leq k\leq s_i$ since $Q^{(j)}$ is vertical for $1\leq j\leq i-1$, and denote by $(i)\cup Q_{i,k}$ the partition obtained by adding $i$ to $Q_{i,k}$ as its first part. 
Hence
$\sigma(Q)=(\widehat{Q}^{(1)}:\cdots:\widehat{Q}^{(\ell)})\in {\bf O}(\la,\ell;n+1)$.
Indeed, $\sigma : {\bf O}(\la,\ell;n) \longrightarrow {\bf O}(\la,\ell;n+1)$ is injective for $n\geq \la_1$, 
and induces an equivalence relation on $\bigsqcup_{n\geq \la_1}{\bf O}(\la,\ell;n)\times \{n\}$, where 
\begin{equation}\label{eq:equiv under sigma}
\text{ $(Q',m)\sim (Q,n)$ if and only if $\sigma^{m-n}(Q)=Q'$,}
\end{equation}
for $Q'\in {\bf O}(\la,\ell;m)$ and $Q \in {\bf O}(\la,\ell;n)$ with $m\geq n$.
For example, if
\[ (Q^{(1)} : Q^{(2)}) = \left( \raisebox{1.5em}{$\ytableausetup{smalltableaux}\
\ydiagram{1, 0, 0, 0, 0}\ , \
\ydiagram{1, 1, 0, 0, 0}\ , \
\ydiagram{1, 1, 1, 0, 0}\ : \
\ydiagram{2, 1, 1, 0, 0}\ , \
\ydiagram{2, 2, 1, 0, 0}\ , \
\ydiagram{2, 2, 1, 1, 0}\ , \
\ydiagram{2, 2, 1, 1, 1}
$}\ \right), \]
then
\[ (\widehat{Q}^{(1)}, \widehat{Q}^{(2)}) = \left( \raisebox{2em}{$\
\ydiagram{1, 0, 0, 0, 0, 0}\ , \
\ydiagram{1, 1, 0, 0, 0, 0}\ , \
\ydiagram{1, 1, 1, 0, 0, 0}\ , \
\ydiagram{1, 1, 1, 1, 0, 0}\ : \
\ydiagram{2, 1, 1, 1, 0, 0}\ , \
\ydiagram{2, 2, 1, 1, 0, 0}\ , \
\ydiagram{2, 2, 2, 1, 0, 0}\ , \
\ydiagram{2, 2, 2, 1, 1, 0}\ , \
\ydiagram{2, 2, 2, 1, 1, 1}
$}\ \right). \]

We define 
\begin{equation*}
{\bf O}(\la,\ell) =\{\,[Q,n]\,|\, Q\in {\bf O}(\la,\ell;n)\ (n\geq \la_1) \,\},
\end{equation*}
where $[Q,n]$ is the equivalence class of $Q\in {\bf O}(\la,\ell;n)$ with respect to \eqref{eq:equiv under sigma}.
We call $[Q,n] \in \mathbf{O}(\lambda, \ell)$ an oscillating tableau of shape $(\la,\ell)$. 

\begin{rem}\label{rem:osc=tab}
{\rm
Let $Q=(Q^{(1)}:\cdots:Q^{(\ell)})\in {\bf O}(\la,\ell;n)$ with $|Q^{(i)}|=s_i$.
Each $Q^{(i)}$ can be identified with a tableau $U^{(i)}$ in $SST_{\mc{I}_n}((1^{s_i}))$ for $1\leq i\leq \ell$, where $a$ (resp. $\ov{a}$) occurs in $U^{(i)}$ if and only if a box is added (resp. removed) in the $a$-th row in $Q^{(i)}$. Hence we may view $Q^{(i)}\in {\bf F}_n$ \eqref{eq:single column tab set}, and apply $\mc{E}$ and $\mc{F}$. Indeed, $\mc{E}Q^{(i)}$ (resp. $\mc{F}Q^{(i)}$) corresponds to removing (resp. adding) two components in $Q^{(i)}$ if it is not ${\bf 0}$.}
\end{rem} 

Next, let us define the weights of the elements in ${\bf O}(\la,\ell)$.
Let ${\bf a}=(a_1,\dots,a_\ell)\in \Z_+^\ell$ be given. We assume that $n$ is sufficiently large so that $n-a_i>0$ for $1\leq i\leq \ell$.
Define ${\bf O}(\la,\ell;n)_{\bf a}$ to be the set of 
$Q=(Q^{(1)}:\cdots:Q^{(\ell)})\in {\bf O}(\la,\ell;n)$  such that
\begin{equation}\label{eq:condition for a}
|Q^{(i)}|=(n-a_i)+2\varepsilon(Q^{(i)})\quad (1\le i\le \ell),
\end{equation}
where $Q^{(i)}$ is viewed an element in the regular $\mf{sl}_2$-crystal ${\bf F}_n$ by Remark \ref{rem:osc=tab}.

\begin{lem}\label{lem:a weight compatible under sigma}
Under the above hypothesis,
\begin{itemize}
\item[(1)] $\varphi(Q^{(i)})+\varepsilon(Q^{(i)})=a_i$ for $1\leq i\leq \ell$,

\item[(2)] $\varepsilon(Q^{(i)})=\varepsilon(\widehat{Q}^{(i)})$ for $1\leq i\leq \ell$, 
where $\sigma(Q)=(\widehat{Q}^{(1)}:\cdots:\widehat{Q}^{(\ell)})$,

\item[(3)] $\sigma\left( {\bf O}(\la,\ell;n)_{\bf a} \right)\subset {\bf O}(\la,\ell;n+1)_{\bf a}$.
\end{itemize}
\end{lem}
\pf Consider $Q^{(i)}$. If $\mc{E}Q^{(i)}={\bf 0}$, then the tableau $U^{(i)}\in SST_{\mc{J}_n}((1^{n-a_i}))$ corresponding to $Q^{(i)}$ is admissible  by Lemma \ref{lem:admissibility of a column}(2). So, in general, each $Q^{(i)}$ belongs to a regular $\mf{sl}_2$-crystal with the highest weight $a_i$ by \eqref{eq:condition for a}, and $\varphi(Q^{(i)})+\varepsilon(Q^{(i)})=a_i$ for $1\leq i\leq \ell$. This proves (1).

Let $\widehat{U}^{(i)}$ correspond to $\widehat{Q}^{(i)}$. Then it is obtained from $U^{(i)}$ by replacing $k$ (resp. $\ov{k}$) with $k+1$ (resp. $\ov{k+1}$), and adding the box $\boxed{1}$ at the top.
Then both $(\widehat{U}^{(i)}_+)^c$ and $\widehat{U}^{(i)}_-$ are obtained by replacing $\bar{k}$ with $\overline{k+1}$ and so we obtain (2).
 
Finally, we have
$|\widehat{Q}^{(i)}|=|Q^{(i)}|+1= (n+1-a_i)+2\varepsilon(Q^{(i)})
=(n+1-a_i)+2\varepsilon(\widehat{Q}^{(i)})$, which implies (3).
\qed\vskip 2mm

Hence by Lemma \ref{lem:a weight compatible under sigma}, we have the following weight decomposition
\begin{equation*}
{\bf O}(\la,\ell) = \bigsqcup_{{\bf a}\in \Z_+^\ell} {\bf O}(\la,\ell)_{{\bf a}}, 
\end{equation*}
where ${\bf O}(\la,\ell)_{{\bf a}}$ is the set of equivalence classes $[Q,n]$ of $Q\in {\bf O}(\la,\ell;n)_{{\bf a}}$. 
We call $[Q,n]\in {\bf O}(\la,\ell)_{{\bf a}}$ an oscillating tableau of shape $(\la,\ell)$ with weight ${\bf a}$.

\begin{ex} \label{ex:oscillating tableau}{\rm
Consider an oscillating tableau
$Q = (Q^{(1)} : Q^{(2)} : Q^{(3)}) \in \mathbf{O}(\lambda, \ell; n)$
with $\lambda = (3, 2, 1), \ell = 3$, and $n = 5$ as follows.
\[ Q = \left( \raisebox{2em}{\ytableausetup{smalltableaux}
\ydiagram{1, 0, 0, 0, 0}\ , \
\ydiagram{1, 1, 0, 0, 0}\ , \
\ydiagram{1, 1, 1, 0, 0}\ , \
\ydiagram{1, 1, 1, 1, 0}\ , \
\ydiagram{1, 1, 1, 1, 1}\ , \
\ydiagram{1, 1, 1, 1, 0}\ , \
\ydiagram{1, 1, 1, 0, 0}\ : \
\ydiagram{2, 1, 1, 0, 0}\ , \
\ydiagram{2, 2, 1, 0, 0}\ , \
\ydiagram{2, 2, 1, 1, 0}\ , \
\ydiagram{2, 2, 1, 1, 1}\ : } \right. \]
\[\hskip 4cm \left. \raisebox{2em}{
\ydiagram{3, 2, 1, 1, 1}\ , \
\ydiagram{3, 3, 1, 1, 1}\ , \
\ydiagram{3, 3, 2, 1, 1}\ , \
\ydiagram{3, 3, 2, 2, 1}\ , \
\ydiagram{3, 3, 2, 2, 0}\ , \
\ydiagram{3, 3, 2, 1, 0}
}\ \right) \]
If we consider each $Q^{(i)}$ as an element in $\mathbf{F}_n$, then we see that
$\varepsilon(Q^{(1)}) = 2, \varepsilon(Q^{(2)}) = 0$, and $\varepsilon(Q^{(3)}) = 1$
and hence the weight of $[Q,n]$ is $\mathbf{a} = (2, 1, 1)$.
} \end{ex}


\subsection{Admissible oscillating tableaux} 

Let $(\la,\ell)\in \mc{P}({\rm Sp})$ be given. 
For $n\geq \la_1$,
let
\begin{equation*}
{\bf O}_{\circ}(\la,\ell;n)=\left\{\,Q\,\Big|\,Q=(Q^{(1)}:\cdots:Q^{(\ell)})\in {\bf O}(\la,\ell;n),\
\varepsilon(Q^{(i)})=0\ (1\leq i\leq\ell)\,\right\}.
\end{equation*}

For ${\bf a}\in \Z_+^\ell$ and a sufficiently large $n$, let 
${\bf O}_{\circ}(\la,\ell;n)_{\bf a}={\bf O}_{\circ}(\la,\ell;n)\cap {\bf O}(\la,\ell;n)_{\bf a}$. By Lemma \ref{lem:a weight compatible under sigma}, we have
$\sigma\left( {\bf O}_\circ(\la,\ell;n)_{\bf a} \right)\subset {\bf O}_\circ(\la,\ell;n+1)_{\bf a}$. Hence we have the following weight decomposition
\begin{equation*}
{\bf O}_\circ(\la,\ell):= \bigsqcup_{{\bf a}\in \Z_+^\ell} {\bf O}_\circ(\la,\ell)_{{\bf a}}, 
\end{equation*}
where ${\bf O}_\circ(\la,\ell)_{{\bf a}}$ is the set of equivalence classes $[Q,n]$ of $Q\in {\bf O}_\circ(\la,\ell;n)_{{\bf a}}$. 
We call $[Q, n] \in {\bf O}_\circ(\la,\ell)$ an {\it admissible oscillating tableau of shape $(\la,\ell)$}.
  
\begin{prop}
For $(\la,\ell)\in \mc{P}({\rm Sp})$ and ${\bf a}=(a_1,\dots,a_\ell)\in \Z_+^\ell$, we have a bijection
\begin{equation}\label{eq:bijection 3}
\xymatrixcolsep{2pc}\xymatrixrowsep{0.5pc}\xymatrix{
{\bf O}(\la,\ell)_{{\bf a}} \ar@{->}[r]  & \ {\bf O}_\circ(\la,\ell)_{{\bf a}}\times \Z/({\bf a}+{\bf 1})\Z \\ [Q,n] \ar@{|->}[r]  & \left([Q_\circ,n], {\varepsilon}(Q)\right)
},
\end{equation}
where $\Z/({\bf a}+{\bf 1})\Z=\Z/(a_1+1)\Z\times\cdots\times\Z/(a_\ell+1)\Z$, and
\begin{equation*}
\begin{split}
Q_\circ &=(\mc{E}^{\rm max}Q^{(1)}:\cdots:\mc{E}^{\rm max}Q^{(\ell)}),\quad
{\varepsilon}(Q)=\left(\varepsilon(Q^{(1)}),\dots,\varepsilon(Q^{(\ell)})\right)
\end{split}
\end{equation*}
for $Q = (Q^{(1)}:\cdots:Q^{(\ell)})\in {\bf O}(\la,\ell;n)_{\bf a}$.
\end{prop}
\pf The map is a well-defined bijection by Lemmas \ref{lem:admissibility of a column}(3) and \ref{lem:a weight compatible under sigma}(2). \qed

\begin{ex} \label{ex:osc to adm-osc}{\rm
Let $[Q, 5] \in \mathbf{O}(\lambda, \ell)$ be given in Example \ref{ex:oscillating tableau}. Then
the image of $[Q, 5]$ under \eqref{eq:bijection 3} is
\[ [(\mathcal{E}^{\rm max} Q^{(1)} : \mathcal{E}^{\rm max} Q^{(2)} : \mathcal{E}^{\rm max} Q^{(3)}), (2, 0, 1)], \]
where
\begin{align*}
\mathcal{E}^{\rm max} Q^{(1)} &= \left( \raisebox{1em}{\
\ydiagram{1, 0, 0}\ , \
\ydiagram{1, 1, 0}\ , \
\ydiagram{1, 1, 1} \ } \right), \\
\mathcal{E}^{\rm max} Q^{(2)} &= \left( \raisebox{2em}{\
\ydiagram{2, 1, 1, 0, 0}\ , \
\ydiagram{2, 2, 1, 0, 0}\ , \
\ydiagram{2, 2, 1, 1, 0}\ , \
\ydiagram{2, 2, 1, 1, 1} \ } \right), \\
\mathcal{E}^{\rm max} Q^{(3)} &= \left( \raisebox{1.5em}{\
\ydiagram{3, 2, 1, 1, 1}\ , \
\ydiagram{3, 3, 1, 1, 1}\ , \
\ydiagram{3, 3, 2, 1, 1}\ , \
\ydiagram{3, 3, 2, 1, 0} \ } \right).
\end{align*}

} \end{ex}


\subsection{King tableaux} 
Let us recall another combinatorial model for irreducible symplectic characters.
For $\ell\geq 2$, let 
$$\mc{J}_\ell=\{\,1<\ov{1}<2<\ov{2}<\cdots<\ell<\ov{\ell}\,\},$$
where we assume that all the entries are of degree $0$.

For $(\la,\ell)\in \mc{P}({\rm Sp})$, let ${\bf K}(\la,\ell)$ be the set of $T\in SST_{\mc{J}_\ell}(\la)$ such that all entries in the $i$th row are larger than or equal to $i$. It is known as the set of {\em King tableaux of shape $\la$}.
For $n\geq \la_1$, let ${\bf K}(\la,\ell;n)$ denote the set ${\bf K}(\la,\ell)$, 
where the columns of the tableaux are enumerated by $n, n-1,\dots, 1$ from the left.

Let $K\in {\bf K}(\la,\ell;n)$ be given. We define a sequence of vertical oscillating tableaux $Q(K;n)=(Q^{(1)}:\cdots:Q^{(\ell)})$ as follows: 
for $1 \le i \le \ell$ and $1 \le j \le n$,
\begin{enumerate}
\item the letter $i$ is contained in the $j$th column of $K$ if and only if there is no step in $Q^{(i)}$ such that a box is added in the $j$th row,

\item the letter $\bar{i}$ is contained in the $j$th column of $K$ if and only if there is a step in $Q^{(i)}$ such that a box is deleted in the $j$th row.
\end{enumerate}

\begin{thm}{\cite[Theorem 2.7]{Lee}} \label{thm:King bijection}
For $\lambda \subseteq (n^\ell)$, we have a bijection
\begin{equation*} 
\xymatrixcolsep{2pc}\xymatrixrowsep{0.5pc}\xymatrix{
{\bf K}(\la,\ell;n) \ar@{->}[r]  & \ {\bf O}(\la,\ell;n)\\
K \ar@{|->}[r]  & Q(K;n)}.
\end{equation*}
\end{thm}

\begin{ex} \label{ex:King to ost}{\rm 
Let $\lambda = (3, 2, 1) \subseteq (5^3)$ with $n=5$ and $\ell=3$, and

\[ 
K = 
\raisebox{1.0em}{\ytableausetup{boxsize=1.2em}
\begin{ytableau} 
\bar{\tl 1} & \bar{\tl 1} & \tl{$2$} \\
\tl{$3$} & \bar{\tl 3} \\
\bar{\tl 3}
\end{ytableau}}
\ \in\ {\bf K}(\lambda, \ell;n). 
\]
Then
\begin{align*}
Q(K; n) &= \left( \raisebox{2em}{\ytableausetup{smalltableaux}
\ydiagram{1, 0, 0, 0, 0}\ , \
\ydiagram{1, 1, 0, 0, 0}\ , \
\ydiagram{1, 1, 1, 0, 0}\ , \
\ydiagram{1, 1, 1, 1, 0}\ , \
\ydiagram{1, 1, 1, 1, 1}\ , \
\ydiagram{1, 1, 1, 1, 0}\ , \
\ydiagram{1, 1, 1, 0, 0}\ : \
\ydiagram{2, 1, 1, 0, 0}\ , \
\ydiagram{2, 2, 1, 0, 0}\ , \
\ydiagram{2, 2, 1, 1, 0}\ , \
\ydiagram{2, 2, 1, 1, 1}\ : } \right. \\
 & \hskip 4cm \left. \raisebox{2em}{
\ydiagram{3, 2, 1, 1, 1}\ , \
\ydiagram{3, 3, 1, 1, 1}\ , \
\ydiagram{3, 3, 2, 1, 1}\ , \
\ydiagram{3, 3, 2, 2, 1}\ , \
\ydiagram{3, 3, 2, 2, 0}\ , \
\ydiagram{3, 3, 2, 1, 0}
}\ \right). \ytableausetup{boxsize=1.2em}
\end{align*}
}
\end{ex}

\begin{cor}
For $\lambda \subseteq (n^\ell)$, we have a bijection
\begin{equation}\label{eq:bijection 4} 
\xymatrixcolsep{2pc}\xymatrixrowsep{0.5pc}\xymatrix{
{\bf K}(\la,\ell) \ar@{->}[r]  & \ {\bf O}(\la,\ell)\\
K \ar@{|->}[r]  & [Q(K;n),n]}.
\end{equation}
\end{cor}
\pf It follows from the definition of $Q(K;n)$ and $\sigma$ that 
$\sigma(Q(K;n))=Q(K;n+1)$ for $K\in {\bf K}(\la,\ell;n)$. Hence the map is a well-defined bijection.
\qed


\section{Symplectic RSK correspondence}\label{sec:RSK for spinor}

\subsection{Pieri rule for spinor model} 

Let ${\bf a}=(a_1,\dots,a_\ell)\in \Z_+^\ell$ be given.
Let ${\bf T}=(T_\ell,\dots,T_1)\in {\bf T}_\A(a_\ell)\times\cdots\times{\bf T}_\A(a_1)$ be given.
We may regard ${\bf T}=(T_\ell,\dots,T_1)\in {\bf T}_\A(\zeta/\eta,\ell)$ for some skew diagram $\zeta/\eta$ with $\zeta, \eta\in \cP_\ell$.

By Theorem \ref{thm:P tableau for spinor}, 
there exists a unique ${\tt P}({\bf T})\in {\bf T}_\A(\la,\ell)$ for some $(\la,\ell)\in \mc{P}({\rm Sp})_\A$, which is obtained by applying ${\tt jdt}_{spin}(\,\cdot\,,c)$ finitely many times with respect to inner corners $c$.

Let us define a recording tableau for ${\tt P}({\bf T})$.
We choose first a sufficiently large $n$.
Let $\ov{\bf T}=(\ov{T}_\ell,\dots,\ov{T}_1)$ be the $n$-conjugate of ${\bf T}$. 
By Corollary \ref{cor:bij spin to KN of skew shape}, $\ov{\bf T}^{\tt ad}\in {\bf KN}_{\alpha/\beta}$ for some skew diagram $\alpha/\beta$.
Let
\begin{equation*}
Q({\bf T};n)={\tt Q}\left(\ov{\bf T}^{\tt ad}\right),
\end{equation*}
where the right-hand side means ${\tt Q}(w)$ for $w=w(\ov{\bf T}^{\tt ad})$ given in Proposition \ref{prop:RSK for KN}.

\begin{lem}\label{lem:recording step 1}
Under the above hypothesis, we have
\begin{itemize}
\item[(1)] $Q({\bf T};n)\in {\bf O}_\circ(\la,\ell;n)_{\bf a}$,

\item[(2)] $\sigma(Q({\bf T};n))=Q({\bf T};n+1)$.

\end{itemize}
\end{lem}
\pf (1) Note that 
$$\ov{\bf T}^{\tt ad}=\llceil\, \ov{T}_\ell^{\tt ad},\dots, \ov{T}_1^{\tt ad} \,\rrceil^{\beta'},$$ 
where $\ov{T}_i^{\tt ad}\in SST_{\mc{I}_n}((1^{n-a_i}))$. 
Let $w^{(i)}=w_{i,1}\cdots w_{i,n-a_i}=w(\ov{T}_i^{\tt ad})$.
For $1\leq i\leq \ell$ and $1\leq k\leq n-a_i$, put 
\begin{equation}\label{eq:construction of Q}
Q_{i,k}={\rm sh}P\left(w^{(1)}\cdots w^{(i-1)}w_{i,1}\cdots w_{i,k}\right),
\end{equation}
where we assume that $w^{(0)}$ is the empty word.
Then by Lemma \ref{lem:admissibility of a column} and Proposition \ref{prop:RSK for KN}, we have
\begin{itemize}
\item[$\bullet$] $Q^{(i)}=(Q_{i,1},\dots, Q_{i,n-a_i})$ is a vertical oscillating tableau with $\varepsilon(Q^{(i)})=0$,

\item[$\bullet$] $Q({\bf T};n)={\tt Q}(w^{(1)}\cdots w^{(\ell)})=(Q^{(1)}:\cdots:Q^{(\ell)})$,
\end{itemize}
which implies that $Q({\bf T};n)\in {\bf O}_\circ(\la,\ell;n)_{\bf a}$.

(2) It can be checked in a straightforward manner. So we leave it to the reader.
\qed\vskip 2mm

Now, we define
\begin{equation}
Q_\circ({\bf T}) = [Q({\bf T};n),n]\in {\bf O}_\circ(\la,\ell)_{\bf a},
\end{equation}
which is well-defined by Lemma \ref{lem:recording step 1}.
The following is one of the main results in this paper.

\begin{thm}
For ${\bf a}=(a_1,\dots,a_\ell)\in \Z_+^\ell$, we have a bijection
\begin{equation}\label{eq:bijection 2}
\xymatrixcolsep{2pc}\xymatrixrowsep{0.5pc}\xymatrix{
{\bf T}_\A(a_\ell)\times\cdots\times{\bf T}_\A(a_1) \ar@{->}[r]  & 
\ \displaystyle{\bigsqcup_{(\la,\ell)\in \mc{P}({\rm Sp})_\A}{\bf T}_\A(\la,\ell)\times {\bf O}_\circ(\la,\ell)_{\bf a}}\\
{\bf T} \ar@{|->}[r]  & ({\tt P}({\bf T}), Q_\circ({\bf T}) )}.
\end{equation}
\end{thm}
\pf Let ${\bf T}$ be given. Choose a sufficiently large $n$ and let $\ov{\bf T}$ be the $n$-conjugate of ${\bf T}$. Put ${\bf U}=\ov{\bf T}$ and ${\bf V}={\bf U}^{\tt ad}$.
By Corollary \ref{cor:compatibility for jdt}, we have the following commuting diagram:
\begin{equation}\label{eq:compatibility of jdt} 
\xymatrixcolsep{4pc}\xymatrixrowsep{2pc}\xymatrix{
{\bf T} \ar@{|->}[r]^{\mc{Y}'}\ar@{|->}[d]_{\ov{\mbox{\ }}}  & {\tt P}({\bf T})\ar@{|->}[d]^{\ov{\mbox{\ }}} \\
{\bf U} \ar@{|->}[r]^{\mc{Y}'}\ar@{|->}[d]_{\mbox{}^{\tt ad}}  & {\tt P}({\bf U})\ar@{|->}[d]^{\mbox{}^{\tt ad}} \\
{\bf V} \ar@{|->}[r]^{\mc{Y}}    & {\tt P}({\bf V})}
\end{equation}
where $\mc{Y}$ is a sequence of ${\tt jdt}_{KN}$'s \eqref{eq:SJDT for KN}, and $\mc{Y}'$ is the corresponding sequence of ${\tt jdt}_{spin}$'s \eqref{eq:sjdt spin general}. Recall $Q_\circ(\mathbf{T}) = [\texttt{Q}(\mathbf{V}), n]$.

Let us first prove the injectivity of the map.
Let ${\bf T}'$ be given such that 
$({\tt P}({\bf T}),Q_\circ({\bf T}))=({\tt P}({\bf T}'),Q_\circ({\bf T}'))$.
Let ${\bf U}'$ be its $n$-conjugate and ${\bf V'}={\bf U'}^{\tt ad}$. 
By definition, this implies that ${\tt P}({\bf V})={\tt P}({\bf V}')$ and ${\tt Q}({\bf V})={\tt Q}({\bf V}')$. 
We claim that
\begin{equation}\label{eq:Pieri for KN} 
\xymatrixcolsep{2pc}\xymatrixrowsep{0.5pc}\xymatrix{
{\bf KN}_{(n-a_\ell)}\times\cdots\times {\bf KN}_{(n-a_1)} \ar@{->}[r]  & 
\ \displaystyle{\bigsqcup_{(\la,\ell)\in \mc{P}({\rm Sp})_n}{\bf KN}_{\rho_n(\la, \ell)}\times {\bf O}_\circ(\la,\ell;n)_{\bf a}}\\
{\bf V} \ar@{|->}[r]  & ({\tt P}({\bf V}),{\tt Q}({\bf V}))}
\end{equation}
is a bijection.
In fact, the map is a morphism of $\mf{sp}_{2n}$-crystals by Proposition \ref{prop:RSK for KN}, which in particular implies that ${\tt Q}({\bf V})$ is constant on its connected component. Using the combinatorial rule of tensor product decomposition \cite{Na}, we see that ${\tt Q}({\bf V})$ uniquely determines a highest weight element in the connected component and hence the map is injective. On the other hand, for a given pair $(H_{\rho_n(\la, \ell)},Q)$ on the right-hand side of \eqref{eq:Pieri for KN}, one can construct directly ${\bf T}$ such that $\te_i{\bf T}={\bf 0}$ for $1\leq i\leq n$ and ${\tt Q}({\bf T})=Q$ again by \cite{Na}.   
This implies the surjectivity of \eqref{eq:Pieri for KN}. 
Hence, we have ${\bf V}={\bf V}'$ and ${\bf T}={\bf T}'$ by \eqref{eq:compatibility of jdt}.

The surjectivity of the map follows from \eqref{eq:compatibility of jdt} and the bijection \eqref{eq:Pieri for KN}.
\qed

\begin{ex} \label{ex:RSK}{\rm
Let $\mathbf{T} = (T_3, T_2, T_1)$ be given in Example \ref{ex:conjugate}. By Example \ref{ex:insertion tableau}, we get

\[ \texttt{P}(\mathbf{T}) = \left( \ \raisebox{2em}{$\begin{ytableau}
\none \\
\none \\
\tl 2 & \tl 2 \\
\tl 3 & \tl 4 \\
\tl{$1'$}
\end{ytableau}$}\ , \
\raisebox{2em}{$\begin{ytableau}
\none \\
\tl 1 & \tl 1 \\
\tl 3 & \tl{$2'$} \\
\tl{$1'$} \\
\tl{$2'$}
\end{ytableau}$}\ , \
\raisebox{2em}{$\begin{ytableau}
\tl 2 & \tl{$2'$} \\
\tl{$1'$} & \tl{$2'$} \\
\tl{$1'$} \\
\tl{$3'$} \\
\tl{$3'$}
\end{ytableau}$ }\ \right).
\]
Since\[ \overline{\mathbf{T}}^{\tt ad} = \raisebox{2em}{$\begin{ytableau}
\none & \none & \tl{$1$} \\
\none & \tl{$2$} & \bar{\tl 5} \\
\tl{$1$} & \bar{\tl 5} & \bar{\tl 4} \\
\tl{$4$} & \bar{\tl 4} \\
\bar{\tl 4} & \bar{\tl 3} \\
\bar{\tl 2}
\end{ytableau}$}\ ,
\]
we have by Example \ref{ex:Lecouvey bijection}

\[ \ytableausetup{smalltableaux} Q_\circ(\mathbf{T}) = \left[ \left(\ \raisebox{1.5em}{$
\ydiagram{1, 0, 0, 0, 0}\ , \ydiagram{1, 1, 0, 0, 0}\ , \ydiagram{1, 1, 1, 0, 0}\ :
\ydiagram{2, 1, 1, 0, 0}\ , \ydiagram{2, 2, 1, 0, 0}\ , \ydiagram{2, 2, 1, 1, 0}\ , \ydiagram{2, 2, 1, 1, 1}\ :
\ydiagram{3, 2, 1, 1, 1}\ , \ydiagram{3, 3, 1, 1, 1}\ , \ydiagram{3, 3, 2, 1, 1}\ , \ydiagram{3, 3, 2, 1, 0} $}
\ \right)\ , 5 \right]\ . \]
} \end{ex}


\subsection{RSK correspondence} 

The goal of this subsection is to establish an analogue of RSK correspondence for spinor model. 
Let 
\[ \mathbf{F}_\mathcal{A}^\ell = \mathbf{E}_\mathcal{A}^{2\ell} \quad (\ell \ge 1). \]
 
\begin{lem}
We have a bijection
\begin{equation*}
\xymatrixcolsep{2pc}\xymatrixrowsep{0.5pc}\xymatrix{
{\bf F}_\A^1 \ar@{->}[r]  & 
\ \displaystyle{\bigsqcup_{a}{\bf T}_{\A}(a)\times \Z/(a+1)\Z}\\
{\bf T} \ar@{|->}[r]  & (\mc{F}^{\rm max}{\bf T}, \varphi({\bf T}))},
\end{equation*}
where the union is over $a\geq 0$ such that ${\bf T}_\A(a)\neq \emptyset$.
\end{lem} 
\pf Let ${\bf T} \in \mathbf{F}_\mathcal{A}^1$ be given. 
Then we have 
$\mc{F}^{\rm max}{\bf T}=\mc{F}^{\varphi}{\bf T}\in SST_{\A}(\la(a,0,c))$
for some $a, c\in \Z_+$ by Lemma \ref{lem:regularity of E}, where $\varphi=\varphi({\bf T})$. It is clear that the map is injective. The connected component of ${\bf T}$ is a regular $\mf{sl}_2$-crystal with highest weight $a$, and the map is also surjective.
\qed\vskip 2mm
 
For ${\bf a}=(a_1,\dots,a_\ell)\in \Z_+^\ell$, put
$$
{\bf T}_\A({\bf a})={\bf T}_\A(a_\ell)\times\cdots\times{\bf T}_\A(a_1).
$$ 
 
\begin{cor}
We have a bijection
\begin{equation}\label{eq:bijection 1}
\xymatrixcolsep{4pc}\xymatrixrowsep{0.5pc}\xymatrix{
{\bf F}_\A^\ell \ar@{->}[r]  & 
\ \displaystyle{\bigsqcup_{\bf a}{\bf T}_\A({\bf a})
\times \Z/({\bf a}+{\bf 1})\Z}\\
{\bf T} \ar@{|->}[r] & (\mc{F}^{\rm max}{\bf T}, \varphi({\bf T}))},
\end{equation}
where  
\begin{equation*}
\begin{split}
\mc{F}^{\rm max}{\bf T}&=(\mc{F}^{\rm max}(U_{2\ell},U_{2\ell-1}),\dots,\mc{F}^{\rm max}(U_{2},U_{1})),\\
\varphi({\bf T})&=(\varphi(U_{2\ell},U_{2\ell-1}),\dots,\varphi(U_{2},U_{1})),
\end{split}
\end{equation*}
for ${\bf T}=(U_{2\ell},\dots,U_1)\in {\bf T}_\A({\bf a})$ and the union is over $\mathbf{a} \in \mathbb{Z}_+^\ell$ such that $\mathbf{T}_\mathcal{A}(\mathbf{a}) \ne \emptyset$.
\end{cor}

\begin{ex} \label{ex:nonadmissible RSK}{\rm
Suppose that $\mathcal{A} = \mathbb{I}_{4|3}$. If 
\[ \ytableausetup{boxsize = 1.2em}
\mathbf{T} = \left( \ \raisebox{1.5em}{$\begin{ytableau}
\none \\
\none \\
\tl 3 \\
\tl{$1'$}
\end{ytableau}$}\ , \
\raisebox{1.5em}{$\begin{ytableau}
\none \\
\tl 2 \\
\tl 4 \\
\tl{$2'$}
\end{ytableau}$}\ , \
\raisebox{1.5em}{$\begin{ytableau}
\tl 1 \\
\tl 2 \\
\tl{$1'$} \\
\tl{$3'$}
\end{ytableau}$}\ , \
\raisebox{1.5em}{$\begin{ytableau}
\none \\
\tl 1 \\
\tl 3 \\
\tl{$2'$}
\end{ytableau}$}\ , \
\raisebox{1.5em}{$\begin{ytableau}
\none \\
\none \\
\tl{$1'$} \\
\tl{$1'$}
\end{ytableau}$}\ , \
\raisebox{1.5em}{$\begin{ytableau}
\tl 2 \\
\tl{$2'$} \\
\tl{$2'$} \\
\tl{$3'$}
\end{ytableau}$}\ \right)
\in \mathbf{F}_\mathcal{A}^3,
\]
then we have
\[ \mathcal{F}^{\tt max} \mathbf{T} = \left( \
\raisebox{1.5em}{$\begin{ytableau}
\tl 2 & \tl 4 \\
\tl 3 & \tl{$2'$} \\
\tl{$1'$}
\end{ytableau}$\ , \
$\begin{ytableau}
\tl 1 & \tl 1 \\
\tl 2 & \tl 3 \\
\tl{$1'$} & \tl{$2'$} \\
\tl{$3'$}
\end{ytableau}$\ , \
$\begin{ytableau}
\tl 2 & \tl{$2'$} \\
\tl{$1'$} & \tl{$2'$} \\
\tl{$1'$} \\
\tl{$3'$}
\end{ytableau}$}\ \right)\ ,
\quad \varphi(\mathbf{T}) = (1, 0, 2). \]
} \end{ex}
 
Now we are ready to state our main result in this paper.
Consider the composition of the following sequence of bijections.
\begin{equation*}
\xymatrixcolsep{3pc}\xymatrixrowsep{0pc}\xymatrix{
{\bf F}_\A^\ell
\ar@{->}[r]^{\hskip -4.5cm\eqref{eq:bijection 1}} & \ 
\quad \displaystyle{\bigsqcup_{{\bf a}\in \Z_+^\ell}{\bf T}_\A({\bf a})
\times \Z/({\bf a}+{\bf 1})\Z}\hskip 4.5cm \\
\mbox{}\quad \quad \ar@{->}[r]^{\hskip -4.5cm\eqref{eq:bijection 2}} &  \quad
\displaystyle{\bigsqcup_{{\bf a}\in \Z_+^\ell}\bigsqcup_{(\la,\ell)\in \mc{P}({\rm Sp})_\A}{\bf T}_\A(\la,\ell)\times {\bf O}_\circ(\la,\ell)_{\bf a}\times \Z/({\bf a}+{\bf 1})\Z}\quad \\
\mbox{}\quad \quad \ar@{->}[r]^{\hskip -4.5cm\eqref{eq:bijection 3}} &  \ 
\displaystyle{\bigsqcup_{(\la,\ell)\in \mc{P}({\rm Sp})_\A}{\bf T}_\A(\la,\ell)\times {\bf O}(\la,\ell)}\hskip 4.3cm\\
\mbox{}\quad \quad \ar@{->}[r]^{\hskip -4.5cm\eqref{eq:bijection 4}} &  \ 
\displaystyle{\bigsqcup_{(\la,\ell)\in \mc{P}({\rm Sp})_\A}{\bf T}_\A(\la,\ell)\times {\bf K}(\la,\ell)}\hskip 4.3cm
}
\end{equation*}

Let $(\texttt{P}(\mathbf{T}), Q(\mathbf{T}))$ denote the image of $\mathbf{T} \in \mathbf{F}_\mathcal{A}^\ell$ under the composition of \eqref{eq:bijection 1}, \eqref{eq:bijection 2}, and \eqref{eq:bijection 3}. 
Let $(\texttt{P}(\mathbf{T}), \texttt{Q}(\mathbf{T}))$ denote the image of $(\texttt{P}(\mathbf{T}), Q(\mathbf{T}))$ under  \eqref{eq:bijection 4}.
Hence we obtain the following correspondence, which is the main result in this paper.

\begin{thm}\label{thm:RSK}
For $\ell\geq 1$, we have a bijection
\begin{equation*}
\xymatrixcolsep{4pc}\xymatrixrowsep{0.5pc}\xymatrix{
{\bf F}_\A^\ell
\ar@{->}[r] & \ 
\displaystyle{\bigsqcup_{(\la,\ell)\in \mc{P}({\rm Sp})_\A}{\bf T}_\A(\la,\ell)\times {\bf K}(\la,\ell)} \\
{\rm \mathbf{T}}
\ar@{->}[r] & \
({\tt P}(\mathbf{T}), {\tt Q}(\mathbf{T}))
}.
\end{equation*}
\end{thm}

\begin{ex}{\rm
Let $\mathbf{T} \in {\bf F}_\mathcal{A}^3$ be the one given in Example \ref{ex:nonadmissible RSK}. 
Combining Examples \ref{ex:nonadmissible RSK}, \ref{ex:RSK}, \ref{ex:insertion tableau}, \ref{ex:osc to adm-osc}, and \ref{ex:oscillating tableau}, we have $(\texttt{P}(\mathbf{T}), Q(\mathbf{T})) \in \mathbf{T}_\mathcal{A}(\lambda, 3) \times \mathbf{O}(\lambda, 3)$ for $\lambda = (3, 2, 1)$, where

\begin{align*}
\texttt{P}(\mathbf{T}) &=
\raisebox{3em}{\ \ytableausetup{boxsize = 1.2em} $\begin{ytableau}
\none \\
\none \\
\tl 2 & \tl 2 \\
\tl 3 & \tl 4 \\
\tl{$1'$} \\
\none
\end{ytableau}$}\ \
\raisebox{3em}{$\begin{ytableau}
\none \\
\tl 1 & \tl 1 \\
\tl 3 & \tl{$2'$} \\
\tl{$1'$} \\
\tl{$2'$}
\end{ytableau}$}\ \
\raisebox{3em}{$\begin{ytableau}
\tl 2 & \tl{$2'$} \\
\tl{$1'$} & \tl{$2'$} \\
\tl{$1'$} \\
\tl{$3'$} \\
\tl{$3'$}
\end{ytableau}$} \\
Q(\mathbf{T}) &= \left( \raisebox{2em}{\ytableausetup{smalltableaux}
\ydiagram{1, 0, 0, 0, 0}\ , \
\ydiagram{1, 1, 0, 0, 0}\ , \
\ydiagram{1, 1, 1, 0, 0}\ , \
\ydiagram{1, 1, 1, 1, 0}\ , \
\ydiagram{1, 1, 1, 1, 1}\ , \
\ydiagram{1, 1, 1, 1, 0}\ , \
\ydiagram{1, 1, 1, 0, 0}\ : \
\ydiagram{2, 1, 1, 0, 0}\ , \
\ydiagram{2, 2, 1, 0, 0}\ , \
\ydiagram{2, 2, 1, 1, 0}\ , \
\ydiagram{2, 2, 1, 1, 1}\ : } \right. \\
 & \hskip 3cm \left. \raisebox{2em}{
\ydiagram{3, 2, 1, 1, 1}\ , \
\ydiagram{3, 3, 1, 1, 1}\ , \
\ydiagram{3, 3, 2, 1, 1}\ , \
\ydiagram{3, 3, 2, 2, 1}\ , \
\ydiagram{3, 3, 2, 2, 0}\ , \
\ydiagram{3, 3, 2, 1, 0}
}\ \right).
\end{align*}

The oscillating tableau $Q(\mathbf{T})$ corresponds to a King tableau $K$ in Example \ref{ex:King to ost} under \eqref{eq:bijection 4}.
Hence ${\tt Q}({\bf T})\in {\bf K}(\la,3)$, where
\[ 
\texttt{Q}(\mathbf{T}) = 
\ytableausetup{boxsize=1.2em} \raisebox{1em}
{$\begin{ytableau}
\bar{\tl 1} & \bar{\tl 1} & \tl{$2$} \\
\tl{$3$} & \bar{\tl 3} \\
\bar{\tl 3}
\end{ytableau}$\ .} 
\]
} \end{ex}

\begin{rem}{\rm
When $\A=[\ov{n}]$, the right-hand side of the bijection in Theorem \ref{thm:RSK} has an $(\mf{sp}_{2n},\mf{sp}_{2\ell})$-bicrystal structure. On the other hand, ${\bf F}_\A^\ell$ is an $\mf{sp}_{2n}$-crystal by \eqref{eq:bijection 1}, and the bijection is an isomorphism of $\mf{sp}_{2n}$-crystals. However, we do not know how to define an $\mf{sp}_{2\ell}$-crystal structure on ${\bf F}_\A^\ell$ directly so that the bijection is an isomorphism of $(\mf{sp}_{2n},\mf{sp}_{2\ell})$-bicrystals.
} 
\end{rem}


\subsection{Cauchy type identity} 

Let ${\bf z}={\bf z}_\ell=\{\,z_1,\dots,z_\ell\,\}$ be formal commuting variables, 
which commute with ${\bf x}={\bf x}_\A=\{\,x_a\,|\,a\in \A\,\}$ (cf.~ Section \ref{subsec:spinor model}). 
For $\mathcal{A} = [n]$, write $\mathbf{x}_n = \mathbf{x}_{[n]} = \{\, x_1, \dots, x_n \,\}$.

Let $(\la,\ell)\in {\mc P}({\rm Sp})$ be given. 
For $K\in {\bf K}(\la,\ell)$, let ${\bf z}^K =\prod_{i\in [\ell]}z_i^{m_i-m_{\ov{i}}}$, where $m_i$ (resp. $m_{\ov{i}}$) is the number of occurrences of $i$ (resp. $\ov{i}$) in $K$. 
Then put
\begin{equation*}
sp_{\la}({\bf z})= \sum_{K\in {\bf K}(\la,\ell)}{\bf z}^K.
\end{equation*}
It is well-known that $sp_{\la}({\bf z})$ is the character of the irreducible highest weight module of ${\rm Sp}_{2\ell}$ with highest weight corresponding to $\la$.

Let ${\bf U}=(U_{2\ell},\dots,U_1)\in {\bf F}_\A^\ell$ be given with $u_i = {\rm ht}(U_i)$. 
Let ${\bf x}^{{\bf U}}=\prod_{i=1}^{2\ell}{\bf x}^{U_i}$ and
${\bf z}^{{\bf U}}=\prod_{i\in [\ell]} z_i^{u_{2i}-u_{2i-1}}$.
Then we have
\begin{equation*}
{\rm ch} {\bf F}_\A^\ell
:= \sum_{\bf U} {\bf x}^{{\bf U}}{\bf z}^{{\bf U}}
=
\prod_{j=1}^\ell \frac{\prod_{a \in \mathcal{A}_0} (1+x_az_j)(1+x_az_j^{-1})}{ \prod_{a \in \mathcal{A}_1} (1-x_az_j)(1-x_az_j^{-1})}.
\end{equation*}

\begin{thm}\label{thm:Cauchy identity}
We have the following identity
\begin{equation*}
t^{\ell}\prod_{j=1}^\ell \frac{\prod_{a \in \mathcal{A}_0} (1+x_az_j)(1+x_az_j^{-1})}{ \prod_{a \in \mathcal{A}_1} (1-x_az_j)(1-x_az_j^{-1})} = 
\sum_{(\lambda, \ell) \in \mc{P}({\rm Sp})_\A}
S_{(\la,\ell)}({\bf x}_\A)sp_\la({\bf z}).
\end{equation*}
\end{thm}
\pf Let ${\bf U}=(U_{2\ell},\dots,U_1)\in {\bf F}_\A^\ell$ be given. 
If ${\bf U}$ is mapped to $({\bf T},K)$ by Theorem \ref{thm:RSK}, then it suffices to show that ${\bf x}^{{\bf U}}{\bf z}^{{\bf U}}={\bf x}^{\bf T}{\bf z}^K$.
Since ${\bf x}^{{\bf U}}={\bf x}^{\bf T}$ is clear, it remains to show that ${\bf z}^{{\bf U}}={\bf z}^{K}$.

Suppose that ${\bf U}$ is mapped to $({\bf T},\upvarphi)$ by \eqref{eq:bijection 1} where 
${\bf T}=(T_\ell,\dots,T_1)\in {\bf T}_\A({\bf a})$ for some ${\bf a}\in \Z_+^\ell$ and 
$\upvarphi=(\varphi_\ell,\dots,\varphi_1)=(\varphi(U_{2\ell},U_{2\ell-1}),\dots,\varphi(U_{2},U_{1}))\in \Z/({\bf a}+{\bf 1})\Z$. 

Choose a sufficiently large $n$, and take the $n$-conjugate $\ov{\bf T}=(\ov{T}_\ell,\dots,\ov{T}_1)$ of ${\bf T}$ 
and $\ov{\bf T}^{\tt ad}=(\ov{T}^{\tt ad}_\ell,\dots,\ov{T}^{\tt ad}_1)$. 
By \eqref{eq:bijection 2}, ${\bf T}$ is mapped to $({\tt P}({\bf T}), Q_0({\bf T}))$. Here ${\tt P}({\bf T})\in {\bf T}_\A(\la,\ell)$ for some $(\la,\ell)\in {\mc P}({\rm Sp})$ and $Q_0({\bf T})=[Q({\bf T};n),n]\in {\bf O}_\circ(\la,\ell)_{\bf a}$ with 
\begin{equation*}
Q({\bf T};n)=(Q^{(1)}:\cdots:Q^{(\ell)}),
\end{equation*}
as given in \eqref{eq:construction of Q}. Note that $|Q^{(i)}|=n-a_i$ for $1\leq i\leq \ell$.
By \eqref{eq:bijection 3}, $(Q({\bf T};n),\upvarphi)$ is mapped to $[Q',n]$ 
where $Q'=(Q^{'(1)}:\cdots:Q^{'(\ell)}) \in {\bf O}(\la,\ell;n)_{\bf a}$. 
Then we have $|Q^{'(i)}|=|Q^{(i)}|+2\varphi_i$ for $1\leq i\leq \ell$.

Let $u_j = {\rm ht}(U_j)$ for $1\leq j\leq 2\ell$, and $t^\pm_i = {\rm ht}(\ov{T}_i^{\tt ad})_\pm$ for $1\leq i\leq \ell$.
By considering the $\mathfrak{sl}_2$-weight of $(U_{2i}, U_{2i-1})$, we have
\begin{equation*}
u_{2i}-u_{2i-1} + 2\varphi_i = (n-t^+_i) - t^-_i = n - (t^+_i + t^-_i) =a_i.
\end{equation*}
On the other hand, it is straightforward to see from the bijection in Theorem \ref{thm:King bijection} that
$$
a_i -2\varphi_i = n-|Q^{'(i)}| = m_i-m_{\ov{i}}
$$
for $1 \le i \le \ell$. Hence $m_i-m_{\ov{i}}=u_{2i}-u_{2i-1}$. This proves ${\bf z}^{{\bf U}}={\bf z}^{K}$.
\qed\vskip 2mm

Let us end this section with well-known identities which can be recovered from Theorem \ref{thm:Cauchy identity} under special choices of $\mathcal{A}$.

First, assume that $\mathcal{A} = [\bar{n}]$. Let
$P = \oplus_{i=1}^n \mathbb{Z}\epsilon_i$ be the weight lattice for $\mf{sp}_{2n}$ in Section \ref{subsec:KN tableaux},
and let $\mathbb{Z}[P]$ be its group ring with a $\mathbb{Z}$-basis
$\{\, e^\mu\, |\, \mu \in P \,\}$. 
Note that
$\varpi_n = \epsilon_1 + \cdots + \epsilon_n$,
the $n$-th fundamental weight.
For $0 \le a \le n$ and $T \in \mathbf{T}_n(a)$, define
\[ {\rm wt}(T) = \varpi_n - \sum_{i=1}^n m_i \epsilon_i, \]
where $m_i$ is the number of occurrences of $\bar{i}$ in $T$.
For ${\bf a}=(a_1,\dots,a_\ell)\in \Z_+^\ell$ and $\mathbf{T} = (T_\ell, \dots, T_1) \in \mathbf{T}_n(\mathbf{a})$, define
$ {\rm wt}(\mathbf{T}) = \sum_{i=1}^\ell {\rm wt}(T_i) $.
For a tableau $K$ with letters in $\mathcal{I}_n$, let
\[ \mbox{wt}(K) = \sum_{i=1}^n (m_i - m_{\ov i}) \epsilon_i, \]
where $m_a$ is the number of occurrences of $a$ in $K$ for $a\in \mc{I}_n$.
It is easy to check that $\mbox{wt}(\mathbf{T}) = \mbox{wt}(\mathbf{T}^{\tt ad})$ for any $\mathbf{T} \in \mathbf{T}_n(\mathbf{a})$, and hence \eqref{eq:Psi_lambda:C} and \eqref{eq:skewPsi_lambda:C} are weight-preserving bijections.

By identifying $x_{\ov i}=x_i^{-1} = e^{-\epsilon_i}\in \Z[P]$ for $i \in [n]$ and $t = e^{\varpi_n} = x_1 \cdots x_n$ in \eqref{eq:character of spinor}, we have
\[ S_{(\lambda, \ell)}(\mathbf{x}_\mathcal{A})
= \sum_{\mathbf{T} \in \mathbf{T}_n(\lambda, \ell)} t^\ell \mathbf{x}_\mathcal{A}^\mathbf{T}
= \sum_{\mathbf{T} \in \mathbf{T}_n(\lambda, \ell)} e^{{\rm wt}(\mathbf{T})}
= \sum_{\mathbf{T}^{\tt ad} \in \mathbf{KN}_{\rho_n(\lambda, \ell)}} e^{{\rm wt}(\mathbf{T}^{\tt ad})}
= sp_{\rho_n(\lambda, \ell)} (\mathbf{x}_n).
\]
The following identity follows immediately from Theorem \ref{thm:Cauchy identity} and the identity $x_i+x_i^{-1}+z_j+z_j^{-1} = x_i(1+x_i^{-1}z_j)(1+x_i^{-1}z_j^{-1})$.

\begin{cor}[\cite{King75}]\label{cor:dual spinors}
For $n, \ell\geq 1$, we have
\[ \prod_{i=1}^n \prod_{j=1}^\ell (x_i+x_i^{-1}+z_j+z_j^{-1}) = 
\sum_{\lambda \subseteq (n^\ell)} sp_{\rho_n(\la, \ell)}(\mathbf{x}_n) sp_\lambda(\mathbf{z}). \]
\end{cor}

Next, assume that $\mathcal{A} = [n]'$. 
For $\ell \ge n$, there exists a bijection in \cite[Theorem 6.5]{K18-3}
\begin{equation*}
\xymatrixcolsep{4pc}\xymatrixrowsep{0pc}\xymatrix{
\mathbf{T}_\mathcal{A} (\lambda, \ell)
\ar@{->}[r] & \ 
\displaystyle{\bigsqcup_{\beta:\mbox{\scriptsize even}} SST_\mathcal{A}(\lambda') \times SST_\mathcal{A}(\beta^\pi),}
}
\end{equation*}
which gives the identity
\[ S_{(\lambda, \ell)}(\mathbf{x}_\mathcal{A}) = t^\ell s_{\lambda'}(\mathbf{x}_\mathcal{A}) \sum_{\beta:\rm even} s_\beta(\mathbf{x}_\mathcal{A}) = t^\ell s_\lambda(\mathbf{x}_n) \sum_{\beta:\rm even} s_{\beta'}(\mathbf{x}_n). \]
Here we call a partition $\beta$ even if all of its parts are even.
Also note that we have
\[ s_\mu(\mathbf{x}_\mathcal{A}) =
\sum_{\mathbf{T} \in SST_\mathcal{A}(\mu)} \mathbf{x}_\mathcal{A}^T =
\sum_{\mathbf{T} \in SST_{[n]}(\mu')} \mathbf{x}_{[n]}^T =
s_{\mu'}(\mathbf{x}_n) \]
for $\mu \in \mathscr{P}$ by identifying $x_{i'} = x_i$ for $i \in [n]$.
By Theorem \ref{thm:Cauchy identity}, we also recover the well-known classical identity due to Littlewood \cite{Li} and Weyl \cite{Weyl}. 

\begin{cor}[\cite{Li, Weyl}]\label{cor:Littlewood} 
For $\ell \ge n \ge 1$, we have
\begin{eqnarray*}
\prod_{i=1}^n \prod_{j=1}^\ell (1-x_iz_j)^{-1}(1-x_iz_j^{-1})^{-1} &=& \sum_{\ell(\la)\leq n} sp_\lambda(\mathbf{z}) s_\lambda(\mathbf{x}_n) \prod_{1 \le i < j \le n} (1-x_ix_j)^{-1} \\
 &=& \sum_{\ell(\la)\leq n} sp_\lambda(\mathbf{z}) s_\lambda(\mathbf{x}_n) \sum_{\beta':\rm even} s_\beta(\mathbf{x}_n).
\end{eqnarray*}
\end{cor}

\begin{rem}{\rm
The bijection in Theorem \ref{thm:RSK} even when reduced to the above cases is completely different from the ones in \cite{Tera} and \cite{Sun} for the identities in Corollaries \ref{cor:dual spinors} and \ref{cor:Littlewood}, respectively, where the insertion algorithm in terms of the King tableaux is used.

}
\end{rem}

\end{document}